\documentclass[preprint,aps]{revtex4}
\usepackage{tabularx}
\usepackage[english]{babel}
\usepackage{verbatim, amsfonts, amssymb, amsthm} 
\usepackage{amsmath}
\usepackage{graphicx}
\usepackage{graphics}
\usepackage{dcolumn}
\usepackage{bm}
\usepackage{epsfig}
\usepackage{color}
\usepackage{hyperref}
\usepackage{subfigure}
\usepackage{psfrag}
\usepackage{float}
\usepackage[english]{babel}
\usepackage{array,warpcol} 
\usepackage{varwidth}
\usepackage{tipa}
\usepackage{booktabs,rotating,longtable}
\newcolumntype{Y}{>{\raggedright\arraybackslash}X}
\newcolumntype{W}{>{\raggedleft\arraybackslash}X}
\newcolumntype{Z}{>{\centering\arraybackslash}X}
\newcolumntype{+}{>{\global\let\currentrowstyle\relax}}
\newcolumntype{^}{>{\currentrowstyle}}
\newcommand{\rowstyle}[1]{\gdef\currentrowstyle{#1}}%
\newcommand{\otoprule}{\midrule[\heavyrulewidth]}

\newcommand{\be}{\begin{equation}}
\newcommand{\ee}{\end{equation}}
\newcommand{\ba}{\begin{eqnarray}}
\newcommand{\ea}{\end{eqnarray}}
\newcommand{\comments}[1]{}

\newcommand{\4}{\frac}

\newcommand{\lb}{\left}
\newcommand{\rb}{\right}

\providecommand{\norm}[1]{\lVert#1\rVert}
\def\bs{\begin{split}}
\def\ess{\end{split}}
\def\be{\begin{equation}}
\def\ee{\end{equation}}
\def\ba{\begin{eqnarray}}
\def\ea{\end{eqnarray}}
\def\bea{\begin{eqnarray*}}
\def\eea{\end{eqnarray*}}
\def\f{\frac}

\begin{document}

\title{\Large \bf Higher order geometric flows on three dimensional locally homogeneous spaces}

\author{ Sanjit Das\footnote{Present address: The Institute of Mathematical Sciences, Chennai 600113, India; Electronic address: sanjit@imsc.res.in,}\hskip0.cm, Kartik Prabhu\footnote{Present address: Department of Physics, University of Chicago, 5640 S. Ellis Avenue, Chicago, IL 60637, USA; Electronic address: kartikp@uchicago.edu}, Sayan Kar}
\email{sayan@phy.iitkgp.ernet.in}
\affiliation{\rm Department of Physics \& Meteorology {\it and} Center for Theoretical Studies \\Indian Institute of Technology, Kharagpur, 721302, India}

\begin{abstract}

We analyse  second order (in Riemann curvature) geometric flows 
(un-normalised) on locally homogeneous 
three manifolds and look for specific features through the solutions
(analytic whereever possible, otherwise numerical) of the 
evolution equations. Several novelties appear in the context of 
scale factor evolution, fixed curves, phase portraits,
approach to singular metrics, isotropisation and  
curvature scalar evolution.
The distinguishing features linked  
to the presence of the second order term in the flow equation are
pointed out. 
Throughout the article, we compare the results obtained, with the 
corresponding results for un-normalized Ricci flows. 
\end{abstract}
\pacs{04.20.-q, 04.20.Jb}

\maketitle

\section{Introduction and overview} 
Ricci and other geometric flows \cite{chowbook,chowbook1} were introduced in mathematics by
Hamilton \cite{hamilton1} and in physics, by Friedan \cite{friedan},
around almost the same time, though with very different motivations.
More recently, such geometric flows have become popular, largely
because of Perelman's work \cite{perelman1} 
which led to the proof of the well--known Poincare conjecture.
In physics, geometric flows have been investigated in
varied contexts such as general relativity, black hole entropy \cite{gr}
and string theory \cite{oli,tseytlin}.

The unnormalised Ricci flow equation \cite{chowbook,chowbook1} is given as, 
\begin{equation}
\frac{\partial g_{ij}}{\partial t}=-2R_{ij}
\end{equation}
where $g_{ij}$ is the metric tensor, $R_{ij}$ the Ricci tensor and
$t$ denotes a time parameter (not the physical time).
For a normalised Ricci flow, the corresponding equation turns
out to be:
\begin{equation}
\frac{\partial g_{ij}}{\partial t}=-2R_{ij} + \frac{2}{n}\langle R\rangle
g_{ij}
\end{equation}
where $\langle R\rangle =\frac{\int R dV}{\int dV}$. Thus the normalised
flow ceases to be different from the unnormalised one if we consider
non--compact, infinite volume manifolds (with a finite value for $\int R dV$).
In addition, for constant curvature manifolds, the second term in the 
normalised flow equation reduces to $\frac{2}{n} R g_{ij}$.

Friedan's early work showed how one may arrive at the 
vacuum Einstein field equations of General Relativity from the 
Renormalisation Group (RG) analysis of a nonlinear 
$\sigma$-model \cite{friedan}. 
The RG $\beta$-function equations (after a $\alpha'$ dependent 
scaling of the flow parameter)
are given (upto third order in the $\alpha'$ parameter) 
as follows \cite{sen,jack}.
\be\label{RG flow mod}
		\frac{\partial g}{\partial t} = - 2Rc - \alpha'  \widehat{{Rc}}-2\alpha'^2 \widehat{\widehat{Rc}}-  \ldots
	\ee
where $\widehat{Rc}$ is a symmetric $2$-tensor defined as: 
$\widehat{Rc} = R_{iklm}R_{jabc}g^{ka}g^{lb}g^{mc} $
and $\widehat{\widehat{Rc}}=\f{1}{2}R_{klmp}{R_i}^{{}{mlr}}{{R_j}^{kp}}_{r}-\f{3}{8}R_{iklj}R^{kspr}{R^l}_{spr} + ....$. In component form, we have,
\be\label{RG flow mod1}
		\frac{\partial g_{ij}}{\partial \lambda} = - 2R_{ij} - \alpha'  R_{iklm}{R_j}^{klm} -  \ldots
	\ee
which may be considered as a {\em higher order} geometric flow equation
for the metric on a given manifold. We intend to continue our investigations
\cite{pdk} on such 
higher (mainly second) order geometric flow equations, in this article. 
We mention that apart from a
motivational perspective, our work has no connection with RG flows
and, in what follows, we will not discuss any link with them anywhere.
$\alpha'$ is treated as a parameter without any link to 
the RG flow equations. 
In other words, the higher order flows we talk about are treated
as geometric flows in their own right. Restricting to second order
we mention a scaling feature of this geometric flow. For Ricci flow
a scaling of the metric $g_{ij}\rightarrow \Omega g_{ij'}$ can be
compensated by a scaling of $t$ by $t'=\frac{t}{\Omega}$. Since the
Ricci is invariant under such scaling, flows for $g_{ij}$ and $t$ are
equivalent with those for $g'_{ij}$ and $t'$. For the second order
flows, things are a bit different. Here, we can see that flows
for $g_{ij}$, $t$ and $\alpha'$ will be equivalent to flows for
$g_{ij}/\Omega$, $t'=t/\Omega$ and $\alpha''= \alpha'/\Omega$. In our
work, we have retained the $\alpha'$ when we do analytic calculations
--however, for numerics we have chosen values of $\alpha'=0,\pm 1$--
solutions for other values of $\alpha'$ can be obtained using the
scaling property mentioned above.  

In this paper, we investigate higher (second) order flows on
locally homogeneous three manifolds.
After the seminal
work of Isenberg and Jackson \cite{isenberg} on
the behavior of Ricci flow on homogeneous three manifolds,
several interesting papers on this topic have appeared 
in the literature.
Notable among them are the articles in 
\cite{glick}, \cite{knopf}, \cite{lauret}, \cite{cgsc}, \cite{csc}. 
In particular, the authors in \cite{glick} discuss a 
new approach wherein the Ricci flow of left-invariant metrics on 
homogeneous three manifolds is equivalently viewed as a
flow of the structure
constants of the metric Lie algebra w.r.t. an evolving
orthonormal frame.  
Consequences for 
some other geometric flows on such three manifolds are available in 
\cite{cotton},\cite{cnsc}.

Homogeneous spaces are also known to play an important role in cosmology 
in the context of the so--called Bianchi models (\cite{ryan},\cite{landau}). 
It is worth mentioning that the Misner-Wheeler minisuperspace  
deals with three dimensional homogeneous spaces in cosmology where 
different paths
in superspace correspond to different evolution profiles for metrics. 
For more details in this context see 
\cite{isenberg}, \cite{ryan}, \cite{wheeler}, \cite{hosoya}.
Recently, a connection between self-dual gravitational
instantons and geometric flows of all Bianchi types has been 
shown in \cite{bakas}.

Our article is organised as follows.
In Section II, we discuss some preliminaries. Section III reviews the
general theory of left invariant metrics. Section IV contains the
main results of the article--here we analyse the higher order flows
for the various Bianchi classes. Finally, in Section V we conclude
with some remarks.

\section{Preliminaries on homogeneous three manifolds}
Following Isenberg and Jackson\cite{isenberg}
(to which we refer for details concerning the following discussion), 
we take the view 
point that our original interest is in closed Riemannian 3-manifolds 
that are locally 
homogeneous. By a result of Singer \cite{singer}, the universal cover of a 
locally homogeneous manifold is homogeneous, i.e. its isometry group acts 
transitively. Here we begin from the basic definition of homogeneous manifolds 
and give Thurston's\cite{thur} proposition for eight geometric structures.
A homogeneous Riemannian space is a Riemannian manifold $(M,g)$ whose 
isometry group $\text{Isom}(M,g)$ acts transitively on $M$, i.e 
given $x,y\in M$ there exists a $\phi\in \text{Isom}(M,g)$ such that 
$\phi(x)=y$. 

 If $G$ is a connected Lie group and $H$ is a closed subgroup, 
$G/H$ is the space of cosets $\{gH\}$, $\pi:G \rightarrow G/H$ 
is defined by $g\rightarrow[gH]$. $G/H$ is called a homogeneous space.
 
Let $(M^{3},G,G_{x})$ be a maximal model geometry i.e. the isotropy 
group $G_{x}$ ($g\in G\vert gx=x$) is maximal among all subgroups of 
the diffeomorphism group of $M^{3}$ that have compact isotropy groups.  
 Depending on the isotropy group, there are eight possible manifestations 
of the maximal model geometry which are also known as Thurston geometries. 
Here we mention them briefly. If the isotropy group is 
$\boldmath {SO(3)}$ then the model geometry is any one of the three, namely, 
$SU(2)$, $\mathbb{R}^{3}$ or $\mathbb{H}^{3}$. On the other hand, 
if $\boldmath{SO(2)}$ is the isotropy group, then it has four possibilities 
depending on whether $M^{3}$  can be a trivial bundle over a two-dimensional 
maximal model or not. The first possibility where  $M^{3}$ is a trivial bundle 
over a $2$-dimensional maximal model, gives rise to $\mathbb{S}^{2}\times\mathbb{R}$ or $\mathbb{H}^{2}\times\mathbb{R}$ upto diffeomorphism and the 
latter produces $\mathrm{Nil}$ or $\widetilde{SL(2,\mathbb{R})}$ up to 
isomorphism. The last case where the isotropy group is trivial, 
the Lie group is  isomorphic to $\mathrm{Sol}$.

\section{General theory for left invariants metrics on 3D unimodular Lie groups}
In studying curvature properties of left invariant metrics on a Lie group, 
many results have been obtained in the past. Most of them are contained 
in a survey article by Milnor \cite{milnor}. For more details on  Lie groups 
and homogeneous spaces we refer to \cite{bese}, \cite{cheeg} and \cite{arvan}. We review this work here
largely because, later, we will use the results quoted here while writing
down the second order geometric flow equations.

The left invariant metric on $M^3$ is taken as:
\be\label{milnormetric}
         \begin{split}
g(t) & =A(t)\eta^1\otimes \eta^1+B(t)\eta^2\otimes \eta^2+C(t)\eta^3\otimes \eta^3\\
       & =\pi^1\otimes \pi^1+\pi^2\otimes \pi^2+\pi^3\otimes \pi^3
          \end{split}
\ee       
and its inverse as,
\be\label{milnormetricinv}
         \begin{split}
g^{-1}(t) & =\f{1}{A(t)}F_1\otimes F_1+ \f{1}{B(t)}F_2\otimes F_2 + \f{1}{C(t)}F_3\otimes F_3 \\
                       & =e_1\otimes e_1+e_2\otimes e_2+e_3\otimes e_3
          \end{split}
\ee
where $A(t)$, $B(t)$, $C(t)$ are positive.
Here ${\{F_i}\}^3_{i=1}$ are the left invariant frame field (also called Milnor's frame)  with dual coframe field ${\{\eta^i}\}^3_{i=1}$ . The Lie brackets w.r.t the left invariant frames are of this form:\\
\be\label{defframe} 
[F_i,F_j]=\alpha_{ijk}F_k
\ee\\
Let us define an orthonormal frame field ${\{e_i}\}^3_{i=1}$which will be 
of the form: $e_1:= \frac{1}{\sqrt A}F_1$. If we let $\zeta_1:=A, \zeta_2:=B, \zeta_3:=C $ and $\lambda:=\alpha_{231}, \mu:=\alpha_{312}, \nu:=\alpha_{123}$ then Eq.(\ref{defframe}) looks like\\
\be\label{defoframe}
[e_i,e_j]=\frac{\zeta_k \alpha_{ijk}}{\sqrt{\zeta_i \zeta_j \zeta_k}}e_k
\ee

The  relevant components of the sectional curvature turn out to be\\
\be\label{K12}
K(e_1\wedge e_2):=\langle Rm(e_1,e_2)e_2,e_1\rangle=\f{(\lambda A-\mu B)^2}{4ABC}+\f{\nu(2\mu B+2\lambda A-3\nu C)}{4AB}
\ee
\be\label{K23}
K(e_2\wedge e_3):=\langle Rm(e_2,e_3)e_3,e_2\rangle=\f{(\mu B-\nu C)^2}{4ABC}+\f{\lambda(2 \nu C+2\mu B-3\lambda A)}{4BC}
\ee

\be\label{K31}
K(e_3\wedge e_1):=\langle Rm(e_3,e_1)e_1,e_3\rangle=\f{(\nu C-\lambda A)^2}{4ABC}+\f{\mu(2\lambda A+2\nu C-3\mu B)}{4AC}
\ee

Recall that these are  actually $R_{1212}$, $R_{2323}$, $R_{3131}$, respectively, in the  orthonormal frame. We have to write them back in Milnor's frame. 
Thus the Riemann  tensor in the Milnor frame will be
\be\label{RiemMil}
Rm(F_{i},F_{j},F_{i},F_{j})=\zeta_{i}\zeta_{j}Rm(e_{i},e_{j},e_{i},e_{j})
\ee
From now on we would mean:$\langle Rm(e_i,e_j)e_j,e_i\rangle=Rm(e_{i},e_{j},e_{i},e_{j})$. Similarly the Ricci and $\widehat{Rc}$ would  look as :
\be\label{RcMil}
Rc(F_{i},F_{i})=\sum_{k=1}^{3} Rm(F_{i},F_{k},F_{i},F_{k})\4{1}{\zeta_{k}}
\ee

\be\label{2ndRcMil}
\widehat{\text{Rc}(F_{i},F_{i})}=2\sum_{k=1}^{3} Rm(F_{i},F_{k},F_{i},F_{k})^{2}\4{1}{\zeta_{k}^{2}}\4{1}{\zeta_{i}}
\ee

\ba
\text{Rc}(F_1,F_1)~=~\4{(\lambda A)^2-(\mu B-\nu C)^2}{2 B C}\\
\text{Rc}(F_2,F_2)~=~\4{(\mu B)^2-(\nu C-\lambda A)^2}{2 C A}\\
\text{Rc}(F_3,F_3)~=~\4{(\nu C)^2-(\lambda A-\mu B)^2}{2 A B }
\ea
\ba
\widehat{\text{Rc}(F_1,F_1)}~=~\4{1}{8 A B^2}\lb[ \4{(\lambda A-\mu B)^2}{C} +\nu(2\mu B+2\lambda A-3\nu C) \rb]^2\nonumber\\
+\4{1}{8 A C^2}\lb[ \4{(\nu C-\lambda A)^2}{B} +\mu(2\lambda A+2\nu C-3\mu B) \rb]^2\\
\widehat{\text{Rc}(F_2,F_2)}~=~\4{1}{8 B C^2}\lb[ \4{(\mu B-\nu C)^2}{A} +\lambda(2\nu C+2\mu B-3\lambda A) \rb]^2\nonumber\\
+\4{1}{8 B A^2}\lb[ \4{(\lambda A-\mu B)^2}{C} +\nu(2\mu B+2\lambda A-3\nu C) \rb]^2\\\
\widehat{\text{Rc}(F_3,F_3)}~=\4{1}{8 A^2C}\lb[ \4{(\nu C-\lambda A)^2}{B} +\mu(2\lambda A+2\nu C-3\mu B) \rb]^2\nonumber\\
+\4{1}{8 B^2 C}\lb[ \4{(\mu B-\nu C)^2}{A} +\lambda(2\nu C+2\mu B-3\lambda A) \rb]^2
\ea
Using the above expressions, we can calculate the scalar curvature and the 
norm of the Ricci tensor which are given as,\\
\begin{subequations}
\begin{align}
\text{Scal}& = \4{1}{A}\lb(\text{Rc}(F_1,F_1)\rb)~+~\4{1}{B}\lb(\text{Rc}(F_2,F_2)\rb)~+~\4{1}{C}\lb(\text{Rc}(F_3,F_3)\rb) \\
& =-\4{A^2 \lambda ^2+(B \mu -C \nu )^2-2 A \lambda  (B \mu +C \nu )}{2A B C}\label{scalcurvformula}
\end{align}
\end{subequations}
\begin{subequations}
\begin{align}
{\norm{\text{Rc}}}^2& =\4{1}{A^2}\lb(\text{Rc}(F_1,F_1)\rb)^2~+~\4{1}{B^2}\lb(\text{Rc}(F_2,F_2)\rb)^2~+~\4{1}{C^2}\lb(\text{Rc}(F_3,F_3)\rb)^2\\
& =\4{3 A^4 \lambda ^4-4 A^3 \lambda ^3 (B \mu +C \nu )-4 A \lambda  (B \mu -C \nu )^2 (B \mu +C \nu )+2 A^2 \lambda ^2 (B \mu +C \nu )^2}{4A^2~B^2~C^2}\nonumber \\
&+\4{(B \mu -C \nu )^2 \left(3 B^2 \mu ^2+2 B C \mu  \nu +3 C^2 \nu ^2\right)}{4A^2~B^2~C^2}
\label{normrcformula}
\end{align}
\end{subequations}
The higher order flow equations therefore turn out to be,\\
\ba
\4{dA}{dt}~=~-\left(2\text{Rc}(F_1,F_1)~+~\alpha'~ \widehat{\text{Rc}(F_1,F_1)}\right)\label{hofeqngeneralone}\\
\4{dB}{dt}~=~-\left(2\text{Rc}(F_2,F_2)~+~\alpha'~ \widehat{\text{Rc}(F_2,F_2)}\right)\label{hofeqngeneraltwo}\\
\4{dC}{dt}~=~-\left(2\text{Rc}(F_3,F_3)~+~\alpha'~ \widehat{\text{Rc}(F_3,F_3)}\right)\label{hofeqngeneralthree}
\ea
In the next section we analyze the flow equations in five different 
homogeneous spaces.
\section{Examples on Bianchi classes}
\subsection{ \underline{Computation on $\text{SU(2)}$}}

\subsubsection{\bf Flow equations}
The canonical three sphere ($S^3$) is topologically equivalent to the Lie group $SU(2)$, which is algebraically represented by
\be
SU(2)=\Bigg\{ \begin{pmatrix} 
z_1 & -z_2 \\
\bar{z_2} & \bar{z_1}
\end{pmatrix} :z_1, z_2\in \mathbb{C}, ~~|z_1|^2+|z_2|^2=1\Bigg\}.
\ee
All the structure constants are the same in the Milnor frame field. We  
assume 
$ \lambda=\mu=\nu=~-2$.
With these values of the structure constants, 
the $2nd$ order flow equations turn out 
to be:
\begin{eqnarray}\label{equ1su2}
\f{dA}{dt}=\f{4(B-C)^2-4A^2}{BC}-2\alpha'[\f{1}{AB^2}\{\f{(A-B)^2}{C}\nonumber \\ +(2B+2A-3C) \}^2+\f{1}{AC^2}\lb\{\f{(A-C)^2}{B} +\lb(2C+2A-3B \rb) \rb\}^2 ]
\end{eqnarray}
\begin{eqnarray}\label{equ2su2}
\f{dB}{dt}=\f{4(C-A)^2-4B^2}{AC}-2\alpha'[\f{1}{BC^2}\{\f{(C-B)^2}{A} \nonumber \\ +(2C+2B-3A)\}^2+\f{1}{BA^2}\{\f{(B-A)^2}{C} +(2B+2A-3C)\}^2]
\end{eqnarray}
\begin{eqnarray}\label{equ3su2}
\f{dC}{dt}=\f{4(A-B)^2-4C^2}{AB}-2\alpha'[\f{1}{CA^2}\{\f{(A-C)^2}{B} \nonumber \\
+(2A+2C-3B)\}^2+\f{1}{CB^2}\{\f{(C-B)^2}{A} +(2C+2B-3A)\}^2]
\end{eqnarray}

\subsubsection{\bf Analytical and numerical estimates}

We first try and see if we can analytically estimate the nature of
evolution for the scale factors A, B, and C. 
 
\noindent {\bf (a)} Without any loss of generality, we can take the initial values of the scale 
factors in an ordered way i.e. $A_0>B_0>C_0$. Further, 
from  Eqn.\ref{equ1su2} to \ref{equ3su2} we can write the equations for the 
differences (pairwise) of the scale factors as follows.
\begin{eqnarray}\label{aminusb}
\f{d(A-B)}{dt}~=~-4\frac{ (A-B) \left(A^2+2 A B+B^2-C^2\right)}{A B C}-4\alpha'\4{(A-B)~\mathcal{G}_1(A,B,C)}{A^2 B^2 C^2}
\end{eqnarray}
\begin{eqnarray}\label{aminusb}
\f{d(A-C)}{dt}~=~-4\frac{ (A-C) \left(A^2+2 A C+C^2-B^2\right)}{A B C}-4\alpha'\4{(A-C)~\mathcal{G}_2(A,B,C)}{A^2 B^2 C^2}
\end{eqnarray}
\begin{eqnarray}\label{aminusb}
\f{d(B-C)}{dt}~=~-4\frac{ (B-C) \left(B^2+2 B C+C^2-A^2\right)}{A B C}-4\alpha'\4{(B-C)~\mathcal{G}_3(A,B,C)}{A^2 B^2 C^2}
\end{eqnarray}
 where $\mathcal{G}_1(A,B,C), ~\mathcal{G}_2(A,B,C), ~\mathcal{G}_3(A,B,C)$ are, respectively
\begin{eqnarray}\nonumber 
\mathcal{G}_1(A,B,C)&=&\left(A^4-4 A^3 B+6 A^2 B^2-4 A B^3+B^4  \right. \\ 
&& \left. +2 A^2 C^2+12 A B C^2+2 B^2 C^2-8 A C^3-8 B C^3+5 C^4\right) 
\end{eqnarray} 
\begin{eqnarray}\nonumber 
\mathcal{G}_2(A,B,C)&=&\left(A^4+2 A^2 B^2-8 A B^3+5 B^4-4 A^3 C \right. \\
&& \left. +12 A B^2 C-8 B^3 C+6 A^2 C^2+2 B^2 C^2-4 A C^3+C^4\right)
\end{eqnarray} 
\begin{eqnarray}\nonumber 
\mathcal{G}_3(A,B,C)&=&\left(5 A^4-8 A^3 B+2 A^2 B^2+B^4-8 A^3 C \right. \\
&& \left. +12 A^2 B C-4 B^3 C+2 A^2 C^2+6 B^2 C^2-4 B C^3+C^4\right)
\end{eqnarray} 
From the above expressions for the differences, it can be 
inferred that $A(t)>B(t)>C(t)$ holds
throughout the evolution, if the initial values satisfy $A_0>B_0>C_0$.
Note that each of (30)-(32) may be formally solved. For example, 
if we consider (31), we have
$A(t)-C(t) = \left (A_0-C_0\right ) \exp \left[ \int \mathcal{K}(t')dt'\right ]$.

\

\noindent {\bf (b)} We now look at the evolution of $B(t)$ for $\alpha'= +1$. We have, from
 Eqn.\ref{equ2su2},
\be
\4{dB}{dt}=4\4{(A^2+C^2)}{AC} -2\{ \4{B^2}{AC}+\4{P_1}{BC^2}+\4{P_2}{BA^2}+4 \}
\le 4\4{(A^2+C^2)}{AC} -2\{ \4{B}{A}+\4{P_1}{BC^2}+\4{P_2}{BA^2}+4 \}
\ee
where $P_1= \{\f{(C-B)^2}{A} +(2C+2B-3A)\}^2$ and $P_2= \{\f{(A-B)^2}{C} +(2B+2A-3C)\}^2$ are strictly positive. Choosing $A_0=n+2, B_0=n, C_0=n-2$   
as the initial values of the scale factors, it can be shown 
that $\4{dB}{dt}\big|_{t=t_0}$ is always negative (from the  
inequality above, it may be checked that the negative term always 
dominates over the positive term). Further, it is easy to see that for 
$n>4$ and $\alpha'=(0,+1)$, all the scale factors are decreasing. 

\

\noindent {\bf (c)} On the other hand, for $\alpha'=-1$ the upper bound 
on $B(t)$ can be written as 
\be
\4{dB}{dt}\le-4\4{B}{A} +\4{4(A-C)^2}{AC}+2\4{P_1}{BC^2}+2\4{P_2}{BA^2}
\ee
Here we may note a difference with the $\alpha'=+1$ case. 
As before, let us choose 
$A_0=n+2, B_0=n, C_0=n-2$ to estimate $F_1~=~-4\4{B}{A} +\4{4(A-C)^2}{AC}+2\4{P_1}{BC^2}+2\4{P_2}{BA^2}$ which
turns out to be, $F_1= \frac{4 (256+n (-64+n (24+n (20-(-3+n) n))))}{n \left(n^2-4 \right)^2}$. It is easy to see that 
$F_1$ is bounded below (i.e.$\lim_{n \to \infty}F_1=-4$ and for any finite 
$n$, $F_1$ is larger than $-4$.).
Thus, the scale factor $B(t)$ can either decrease or increase depending on the 
choice of initial values.  

\

\begin{figure}[h]
\centering
\subfigure[$A_{0}=7, B_{0}=5, C_{0}=3, \alpha' = 1, T_s=0.702$]{\includegraphics[width=0.4\textwidth]{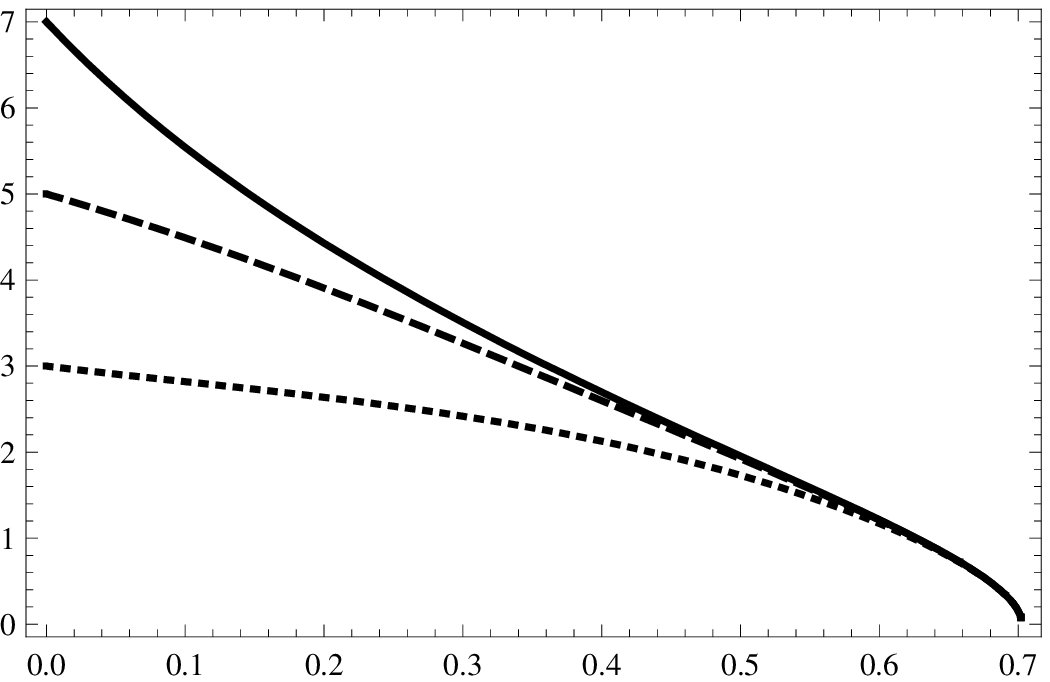}\label{su2one}}
\subfigure[$A_0=11, B_{0}=9, C_{0}=7,\alpha' =-1,T_s=-1.09$]{\includegraphics[width=0.4\textwidth]{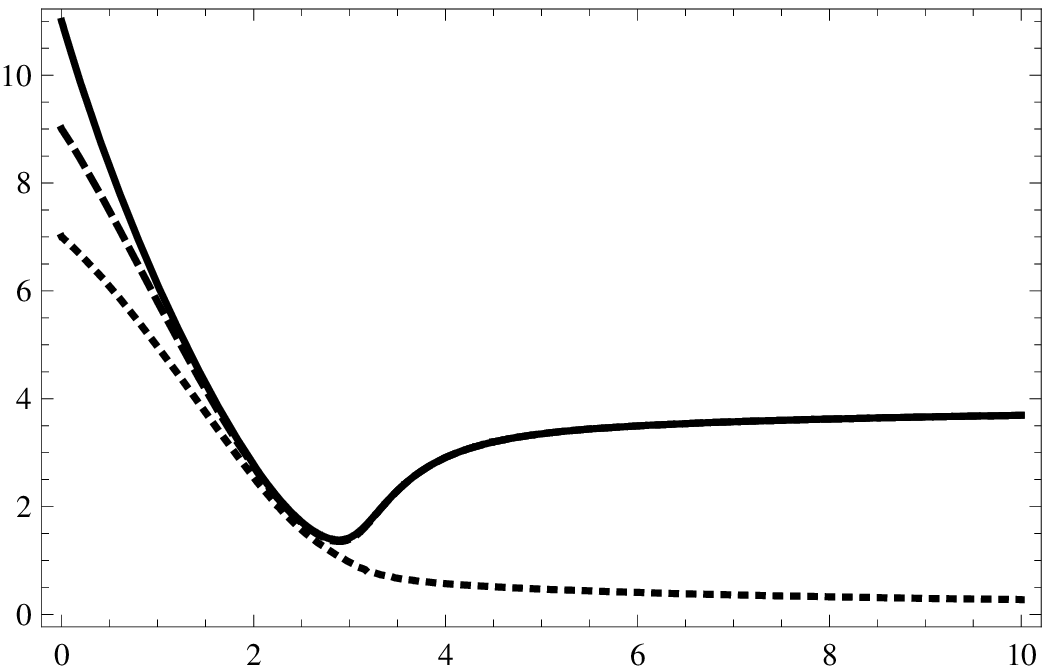}\label{su2three}}
\subfigure[$A_0=5.5, B_{0}=3.5, C_{0}=1.5,\alpha' =-1,T_s=-0.04$]{\includegraphics[width=0.4\textwidth]{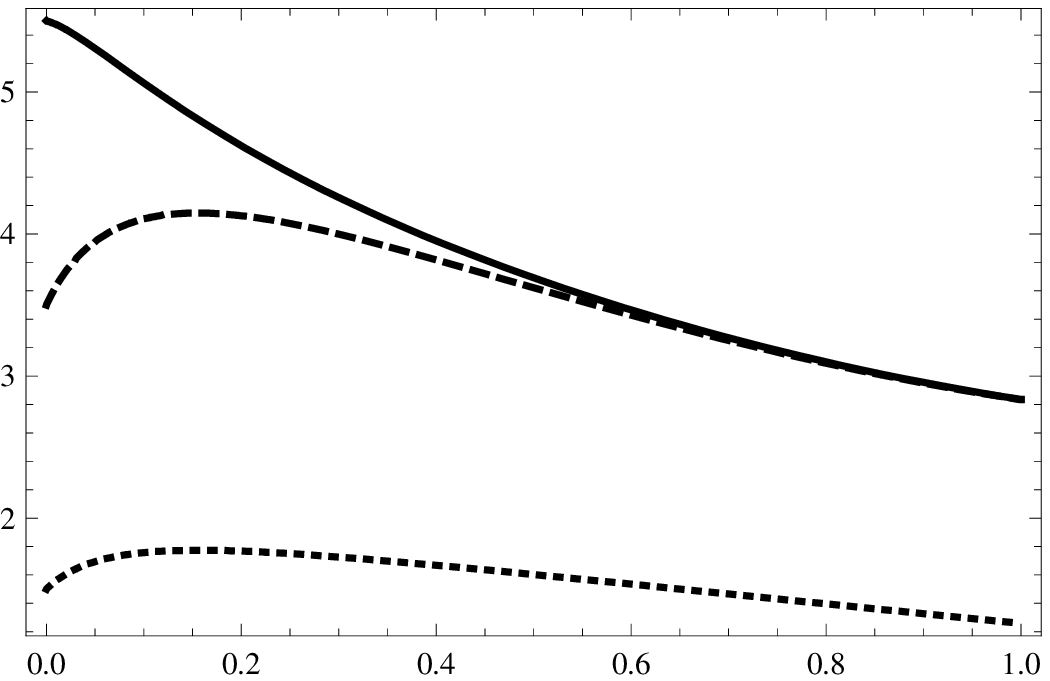}\label{su2three_1}}
\subfigure[$A_0=7, B_{0}=5, C_{0}=3, \alpha' =0,T_s=2.76$]{\includegraphics[width=0.4\textwidth]{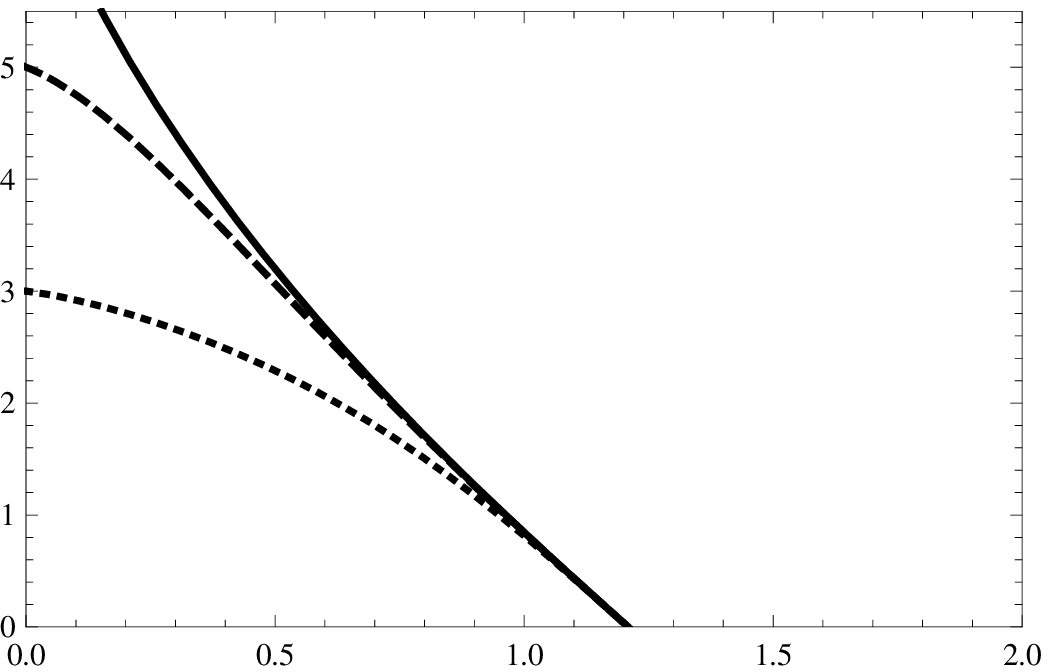}\label{su2two}}
\caption{$A(t),B(t),C(t)$ vs $t$ for  $\alpha'=1$, $\alpha'=-1$ and $\alpha'=0$ for SU(2)}
\label{su2nuevl1}
\end{figure}

\noindent {\bf (d)} It is easy to see from the flow equations that for 
$\alpha'=(0,-1)$, $A(t)$, $B(t)$ and $C(t)$ will 
monotonically decrease 
when $A > (B-C)$. For typical initial values obeying this condition, 
this decreasing feature is shown in Fig. \ref{su2one}and Fig.\ref{su2two}
where
$A(t), B(t)$ and $C(t)$ are obtained by numerically solving the
dynamical system. In the graphs, $A(t)$, $B(t)$ and C(t)) 
are shown in thick, 
dashed and dotted lines, respectively. We follow this convention in our plots,
throughout this article. 

\

\noindent {\bf(e)} We now ask whether all the scale factors can be increasing in $t$. 
For $\alpha'=0$ this is not possible. To see this, 
consider the initial value assumption $A_0>B_0>C_0$. The scale factor $A(t)$ 
will increase for $A<(B-C)$, whereas both $B(t)$ and $C(t)$ will increase 
for $A>(B+C)$. However, both these inequalities cannot hold simultaneously.
Hence all scale factors cannot be increasing in $t$.
 
For $\alpha'=-1$ we may be tempted to believe in a conclusion similar to the
$\alpha'=0$ case mentioned above. However, the higher order contribution 
leads to some unexpected behaviour. 
In Fig.\ref{su2three_1} we see that the scale factor $B(t)$ initially 
increases but eventually decreases after a certain time. This 
is due to the the fact that the $A>(B+C)$ criterion does not hold in the
entire domain of $t$, a fact apparent from the figure. 
In contrast, for Fig.\ref{su2three} we have taken 
initial values 
which do not satisfy the aforesaid condition and though, initially, 
$B(t)$ is decreasing we find that after some time, $B(t)$ increases. 
The reason behind this behaviour is --- the higher order contribution is 
much more positive and the net effect makes $B(t)$ increase after 
an initial decreasing phase.


\

\noindent {\bf (f)}  recall that $SU(2)$ admits round Einstein metrics and
the Ricci flow converges. This feature remains even when we include 
higher order 
terms with $\alpha' = 1$. If we change the value of $\alpha'$(from $1$ to $0$) 
a change appears only in the singularity time, 
though no difference appears in the nature of the graph. 
On the other hand, when $\alpha' = -1$, we notice an expansion of the 
scale factors -- $B(t)$ and $A(t)$ appear to be expanding beyond the minimum in Fig.[\ref{su2three}] at the same rate, whereas $C(t)$ goes to zero 
asymptotically. 

Let us now move on to some special cases.
\subsubsection{\bf Special Case: $A\ne B=C $}
\begin{figure}[htbp] 
\centering
\subfigure[$A_0=7, B_0=5,\alpha' =1,T_s=0.914$]{\includegraphics[width=0.4\textwidth]{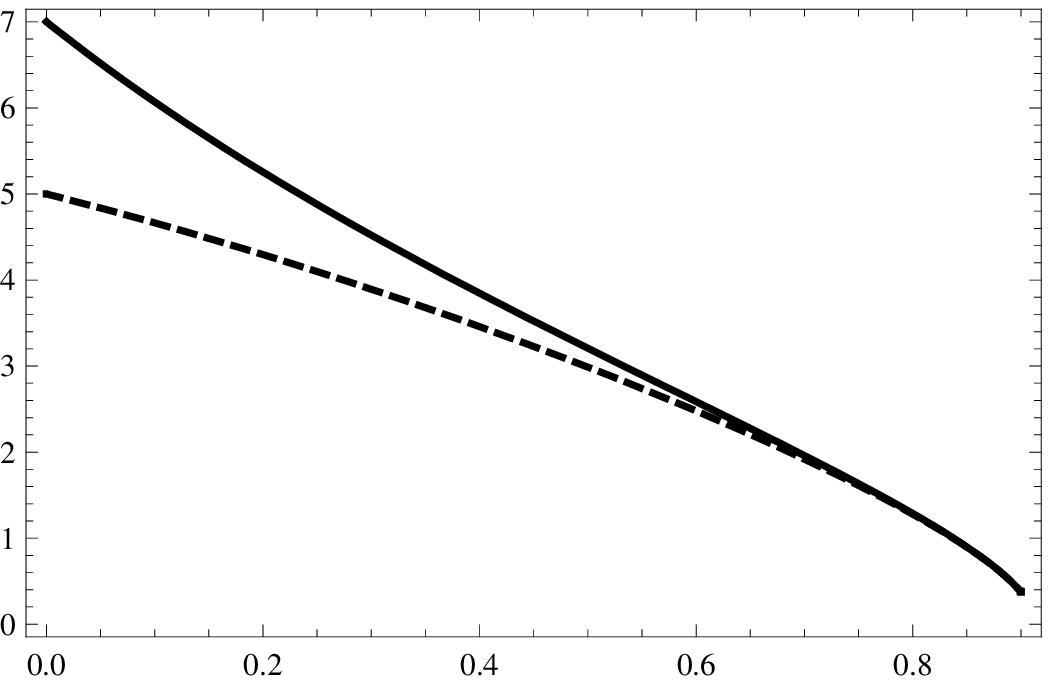} \label{2ndorderSU2one}}
\subfigure[$A_0=7, B_0=5,\alpha' =-1,T_s=-0.602$]{\includegraphics[width=0.4\textwidth]{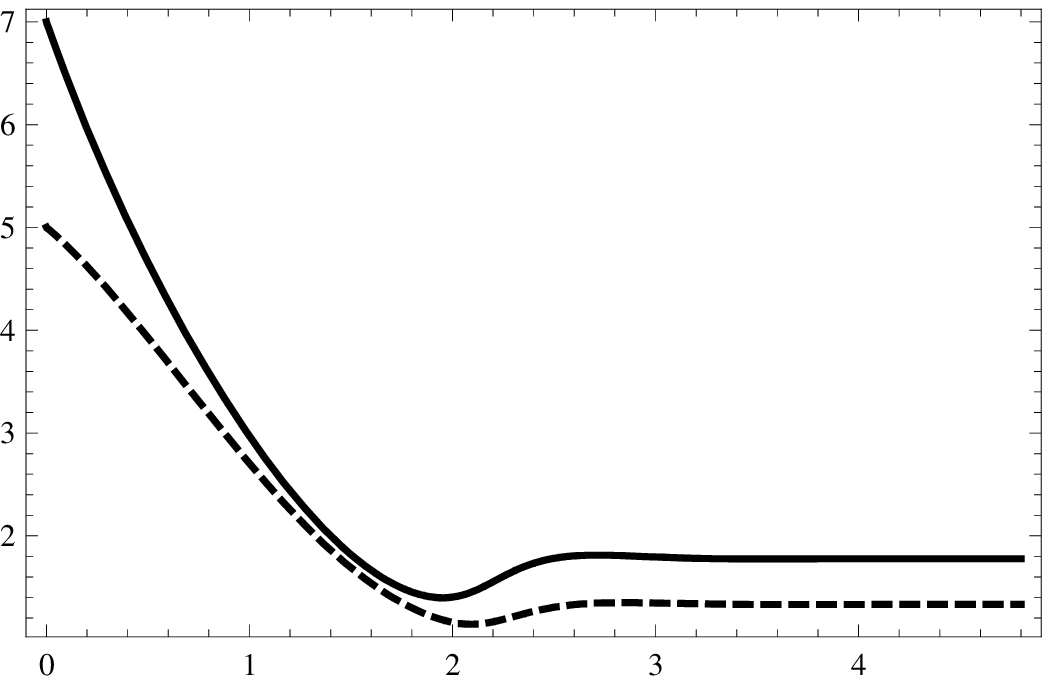} \label{2ndorderSU2two}}
\subfigure[$A_0=1, B_0=0.5,\alpha' =-1,T_s=-0.004$]{\includegraphics[width=0.4\textwidth]{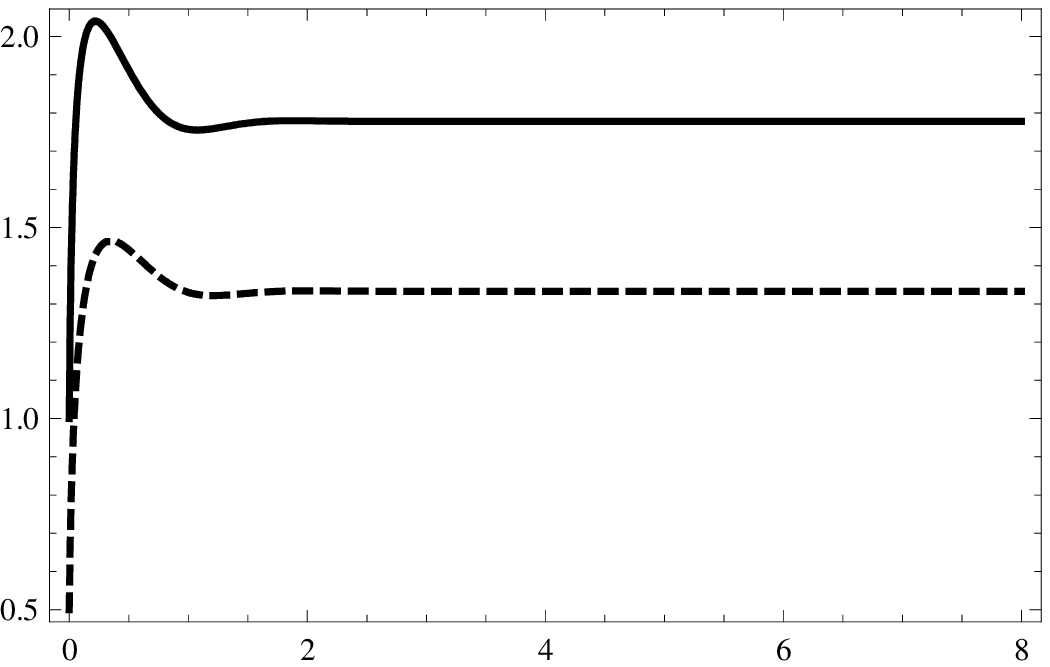} \label{2ndorderSU2twonew}}
\subfigure[$A_0=7, B_0=5,\alpha' =0,T_s=(-0.76)$]{\includegraphics[width=0.4\textwidth]{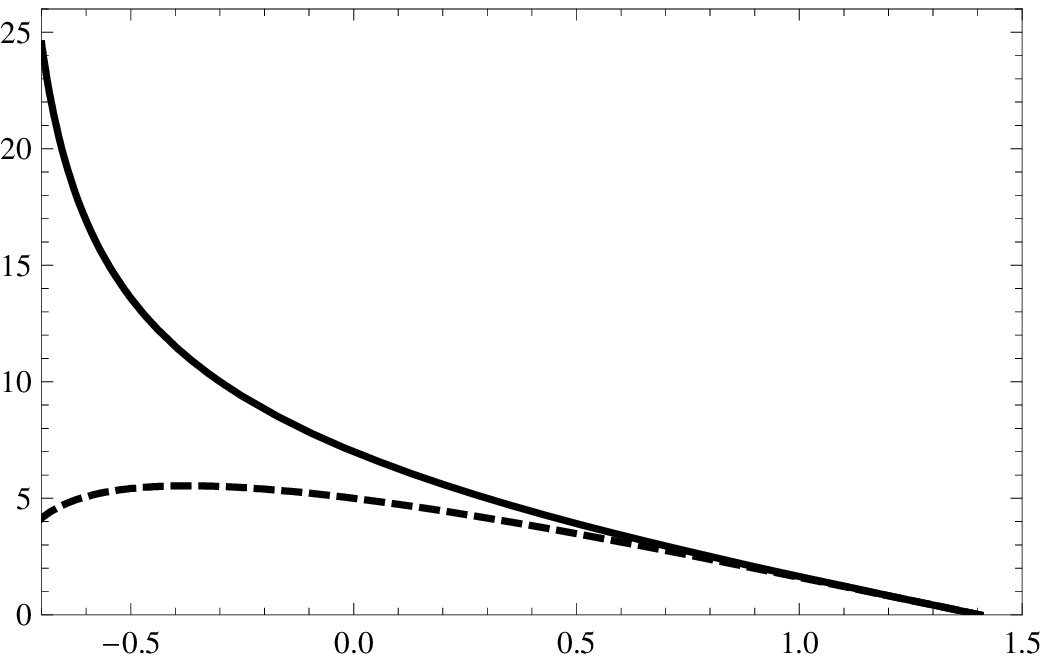} \label{2ndorderSU2three}}
\caption{$A(t),B(t)$ vs $t$ for  $\alpha'=1$, $\alpha'=-1$ and $\alpha'=0$ for SU(2)}
\label{2ndorderSU2}
\end{figure}

The flow equations, with this restriction (leading to the so-called Berger
sphere metrics \cite{chowbook}) become,
\be\label{2equ1su2}
\f{dA}{dt}= -4\frac{A^2}{B^2} - 4\alpha' ~\frac{A^3}{B^4} 
\ee
\be\label{2equ2su2}
\f{dB}{dt}= -8 + 4\frac{A}{B} - 4\alpha' ~\bigl[\4{5}{B}(\4{A}{B}-1)^2+\4{3}{B}(1-\4{2A}{3B})\bigr]
\ee
Let us first look at the case, $\alpha'=0$. It is clear from Eqn. (38) that
$A(t)$ decreases with $t$ and has a upper bound,
$A(t)\leq A_0 - 4t$. However, the same is not true for $B(t)$.
If $ A(t)<2B(t)$ then $B(t)$ decreases with increasing $t$, otherwise it
increases. $B(t)$ does have a lower bound, $B(t)\ge B_0-8t$. 
Depending on the sign of the R. H. S. of (39) $B(t)$ can have an
increasing, stationary or decreasing behaviour during its evolution in $t$. 
From Eqns. (38), (39), one can arrive
at a second order equation for $B(t)$ given as:
\begin{equation}
B(t)\frac{d^2 B}{dt^2} + 2 \left (\frac{d B(t)}{dt}\right )^2
+ 24 \frac{d B(t)}{dt} = -64
\end{equation}
There are trivial solutions to this equation:  $\frac{dB(t)}{dt}=-8$
(which corresponds to $A(t)=0$) and $\frac{dB(t)}{dt}=-4$ (which corresponds to
$A(t)=B(t)$). 
If we define $\frac{dB}{dt}=D(t)$ then we can see the following:
(i) At some $t$ if $D=0$ then $\frac{dD}{dt} <0$ (maximum);
(ii) If, over a range of $t$,  $D>0$ then $\frac{dD}{dt} <0$;
(iii) If, over a range of $t$, $D<0$ then $\frac{dD}{dt}$ can be $<0$, $>0$ or
$=0$. 
We note this behaviour in Fig.\ref{2ndorderSU2three}.
For example, the maximum occurs at the point where $D(t)=0$--substitute
the values of $A(t)$ and $B(t)$ at the 
maximum as found in the graph, into Eqn. (39) (with $\alpha'=0$) and check that it is indeed zero. 

Now we wish to note what happens to the rate of 
change of $(A-B)$ and $(A+B)$. We show that $(A-B)$ decreases faster 
than $(A+B)$ which implies that as $A$ and $B$  decrease, they approach 
each other. From the above equations for $\alpha'=0$
we can write 
\be
\frac{d(A-B)}{dt}=\4{-4(A+2B)}{B^2}(A-B)\\
\le -\frac{4A}{B^2}(A^2+3AB+B^2) \\
\le-20\4{A}{B}
\ee
On the other hand,
\be
\frac{d(A+B)}{dt}=-\frac{4}{B^2}(A^2+2B^2)+4\4{A}{B}\\
\le-4 \4{A}{B}
\ee
It is easy to note that $A-B$ decreases faster than  $A+B$. Thus
$A$ and $B$ approach each other while decreasing. 
For $\alpha'\ne 0$ the analysis can be done following similar logic
but we prefer to solve the equations numerically and demonstrate
our conclusions.

If we plot $A(t)$ (thick line), $B(t)$ (dashed line) treating the evolution equations  as a dynamical system, we see that 
both the scale factors converge to the origin (Fig.\ref{2ndorderSU2one},\ref{2ndorderSU2three}) for $\alpha'=0,1$. However, for $\alpha'=-1$ the evolution of the scale factors depend entirely on the initial values of $A,B,C$-- a fact
shown in the two figures, Fig.\ref{2ndorderSU2two},\ref{2ndorderSU2twonew}.

\subsubsection{\bf Phase plots}

\begin{figure}
	\centering
	\subfigure[$\alpha' = 1$]{\includegraphics[width=0.3\textwidth]{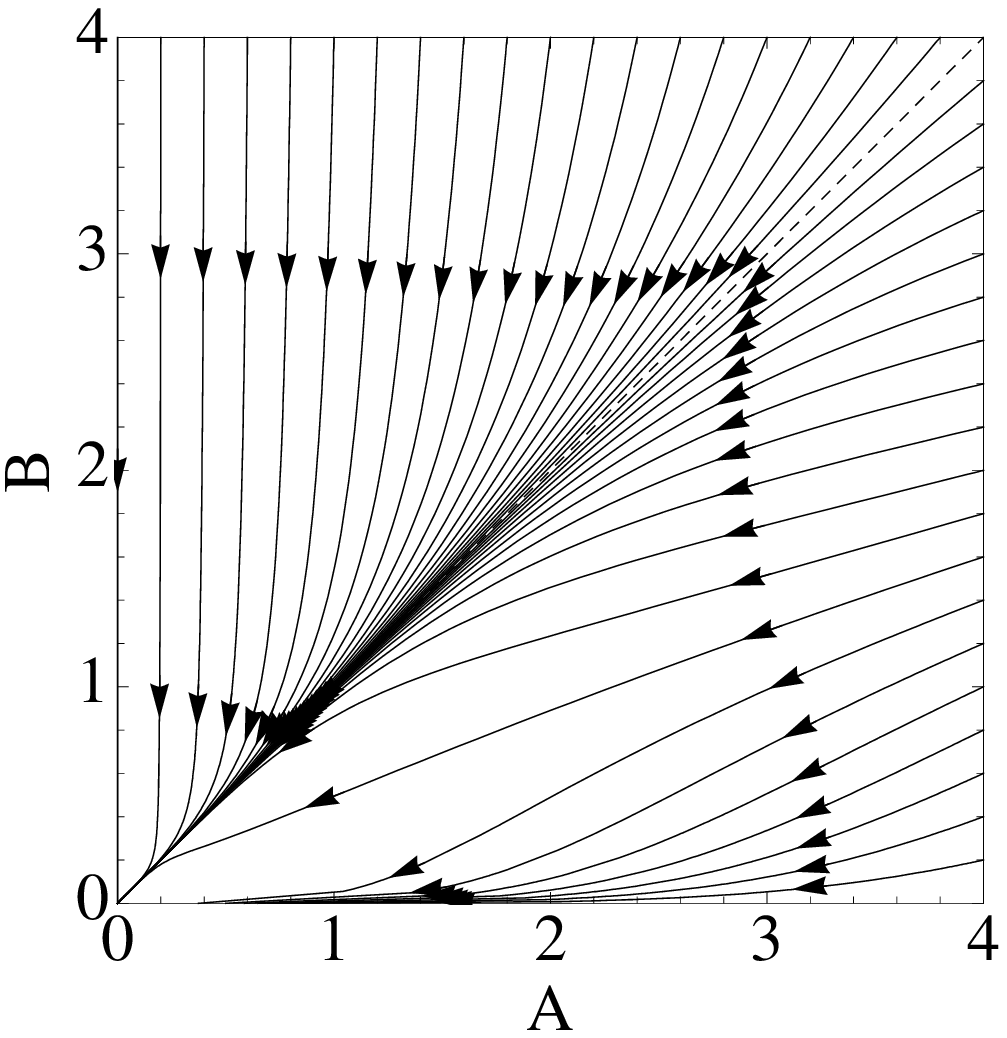}\label{su2 1}}
	\subfigure[$\alpha' = 0$]{\includegraphics[width=0.3\textwidth]{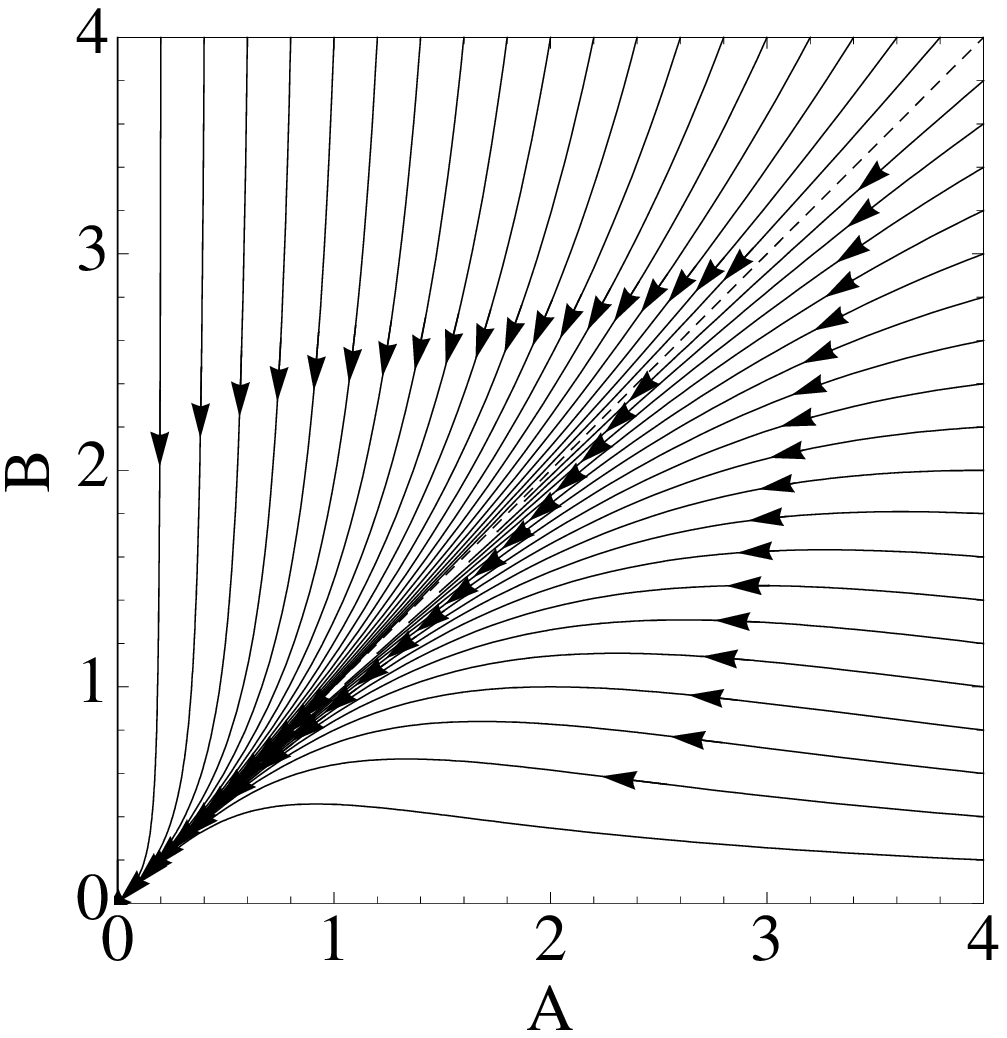}\label{su2 0}}
	\subfigure[$\alpha' = -1$]{\includegraphics[width=0.3\textwidth]{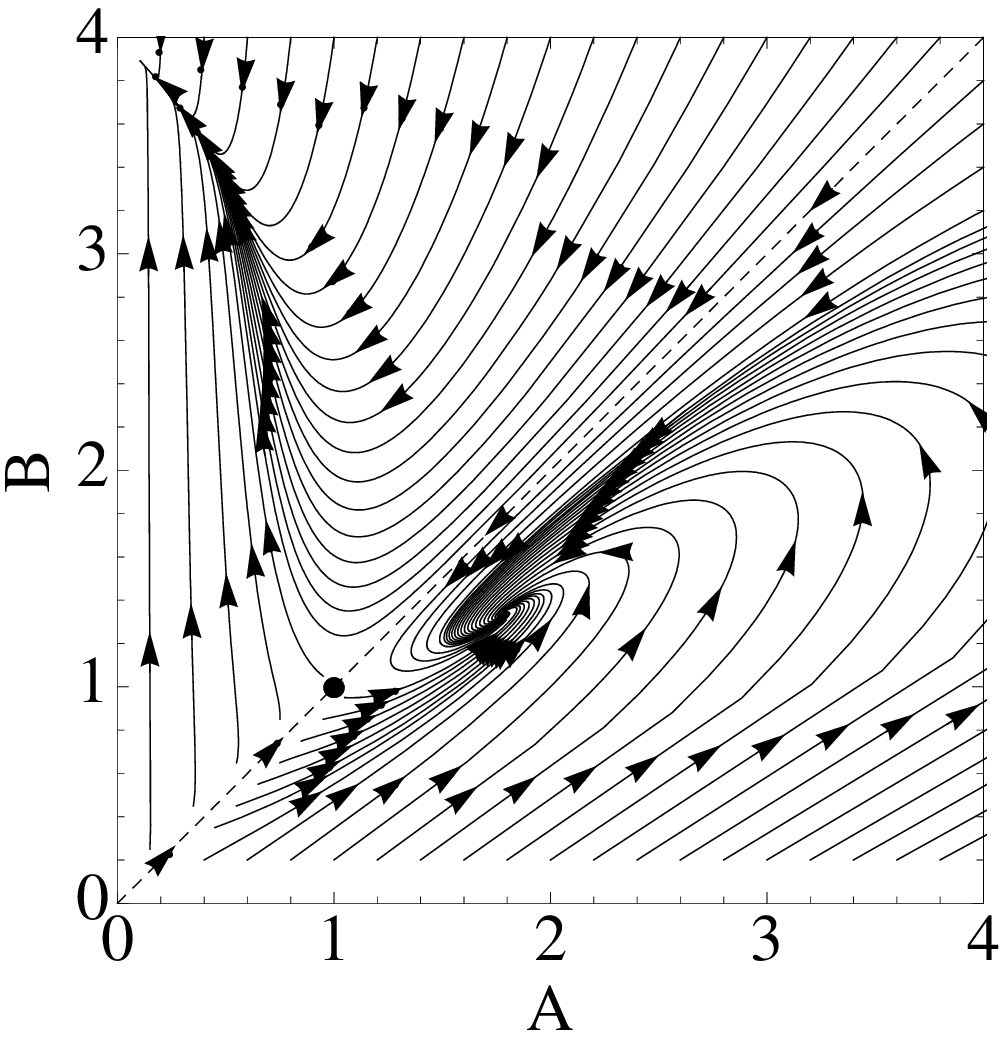}\label{su2 -1}}
	\caption{$2nd$ order  flow on SU(2)for $B=C$}
	\label{su2}
	\end{figure}

We obtain phase plots of the above dynamical system in Fig.\ref{su2}. 
The curves are the trajectories of the system for varying initial conditions, 
and the arrows indicate the direction of increasing time. From the phase 
plots we can deduce the qualitative and quantitative behaviour of the flow.
We can see that the $A=B$ line is a critical curve demarcating the regions 
in phase space with different behaviour. Also, for $\alpha' = 0$ 
all flow trajectories converge towards the $A=B$ line. This leads to isotropization 
of the manifold. However, for $\alpha' = 1$, trajectories from a region in 
the phase space with $A>B$ converge to the singularity $B=0$ before 
isotropization. The behaviour for $\alpha' = -1$ is markedly different 
due to the existence of fixed points given by $(A,B)= (0,4), (\frac{16}{9},\frac{4}{3}), (1,1)$. 
In this case, trajectories from $A>B$ converge to $(\frac{16}{9},\frac{4}{3})$ and those from $A<B$ converge to $(0,1)$. When $A=B$, we note a convergence 
to $(1,1)$. Therefore, in this case there is no 
isotropization, with the asymptotic manifold being non-isotropic -
a fact which is manifest 
in the values of the fixed points.

\subsubsection{\bf Curvature evolution}
	
	Finally, we show the evolution of the scalar curvature. From Eqn.[\ref{scalcurvformula}] and Eqn.[\ref{normrcformula}]
we can immediately compute the  scalar curvature or norm of the Ricci tensor, 
for locally homogeneous three manifolds in 
the different Bianchi classes. These quantitites depend on 
the scale factors. Thus, if we know the behaviour of the scale factors,
under a given flow, we can easily obtain the profile of any 
curvature invariant as a function of the flow parameter. 
However, in most cases, we do not have exact solutions for the scale factors. 
Therefore, the numerical evaluations of the scale factors 
(for various initial values) are used to obtain the curvature invariants. From 
Eqn.\ref{scalcurvformula} we can see, for $SU(2)$, the scalar curvature is as follows
\be
\text{Scal}=\4{-2\left[A^2+(B+C)^2-2A(B+C)\right]}{ABC}=\4{2~ F(A,B,C)}{ABC}
\ee
Therefore, it can be positive, negative or zero depending on the values of $A,B,C$. More specifically  
if $A=B=C$ then $F=3A^2>0$. Similarly if we take $A\ne B=C$ then $F=A(4B-A)$, which can be zero, negative or positive for $A=4B$, $A>4B$ or $A<4B$ respectively.\\
Fig. 4 shows the profiles for the scalar curvature. 
For different initial values and three different $\alpha'$ we have shown the 
evolution. Both $\alpha'=0$ and $\alpha'=1$ show a monotonic nature of the 
scalar curvature. 
For $\alpha'=-1$  the scale factors are not always monotonic
-- the behavior depends on the initial conditions. The scalar 
curvature for $\alpha'=-1$ rises initially but asymptotically reaches 
a constant value. The figures in the first row, second column and
in the third row, second column show the appearance of negative 
scalar curvature.  
We note here that $Scal$ remains entirely negative or entirely positive 
during evolution, for 
$\alpha'=0,1$, but may grow from negative to positive value for
$\alpha'=-1$.
	\begin{figure}[htbp]
	{\includegraphics[width=0.72\textwidth]{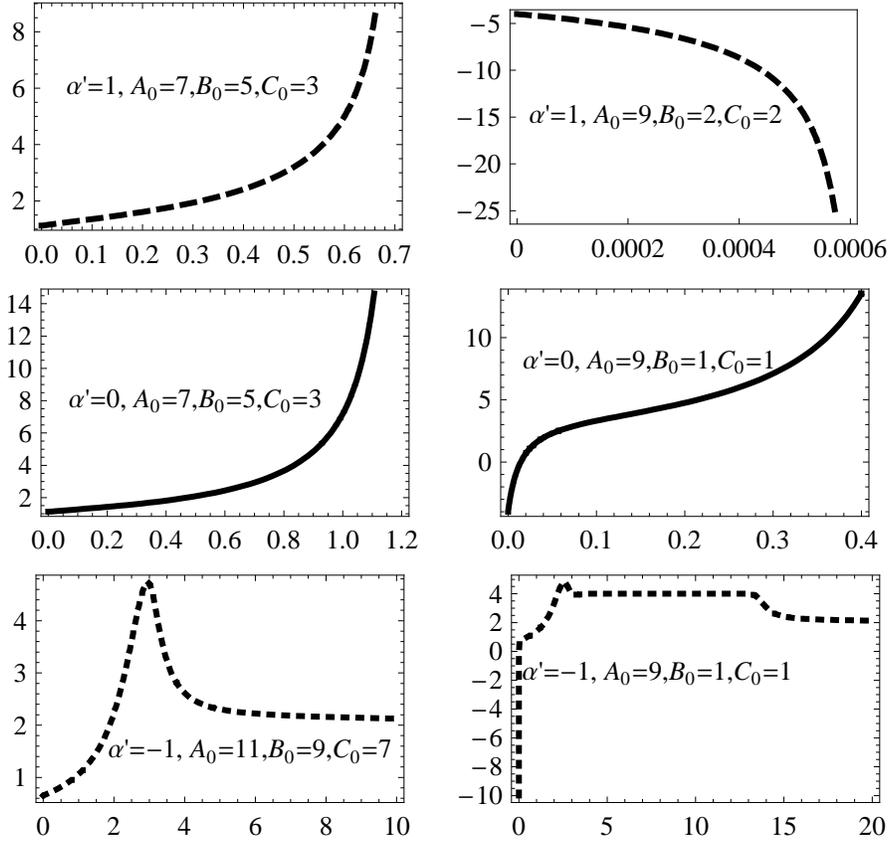}\label{curvevolsu2one}}
	\caption{Evolution of Scal  for $SU(2)$.}
	\end{figure}
	
In summary, several distinguishing features of scale factor evolution
and curvature evolution arising exclusively due to the presence of the
higher order term in the flow, have been pointed out in the above analysis 
and summarised in Table \ref{table1chap03}.

 \begin{table}[htbp]\small
\caption{Comparison of the evolution of different scale factors and Scal. 
curvature for $SU(2)$}\label{table1chap03}\centering
\renewcommand{\tabularxcolumn}[1]{>{\arraybackslash}m{#1}}
\begin{tabularx}{\textwidth}{+Y^Z^Z^Z^Z}
\toprule\rowstyle{\bfseries}
 & $A(t)$ & $B(t)$ & $C(t)$ & $Scal$\\
\otoprule
$\alpha'=0$  & decreasing & decreasing & decreasing & increases, can be -ve.\\
\midrule
$\alpha'=1$& decreasing much faster & decreasing much faster & decreasing much faster& increase/decrease (entirely +ve or -ve.)\\
\midrule
$\alpha'=-1$&  depending on initial value &  depending on initial value & decreasing & not monotonic (can be +ve, -ve)\\
\bottomrule
\end{tabularx}
\end{table}


\subsection{\underline{Computation on Nil}}

\subsubsection{\bf Flow equations}

Nil is a three dimensional Lie group consisting of all $3\times 3$ matrices of the form \\
\[ \left( \begin{array}{ccc}
1 & x & z \\
0 & 1 & y \\
0 & 0 & 1 \end{array} \right)\]
under multiplication, also known as the Heisenberg group. The group action 
in $\mathbb{R}^3$ can be written as 
\be
(x,y,z)\circ(x',y',z')=(x+x',y+y',z+z'+xy')
\ee This is a nilpotent group. We can put a 
left invariant frame field on this Lie group. With the inherent 
metric, \dag\footnotetext[2]{ One such left invariant metric is $ds^2=dx^2+dy^2+(dz-xdy)^2$.}it is a line bundle over the Euclidean plane $\mathbb{E}^2$. Following Milnor \cite{milnor} we have the structure constants, $\lambda=-2,\mu=\nu=0$. The Rc and ${\widehat{Rc}}_{}^{2}$ tensors can be obtained by simple computation. The flow equations in this case will be:
\be\label{equ1Nil}
\f{dA}{dt}=-\f{4A^2}{BC}-4\alpha'\f{A^3}{B^{2}C^{2}}
\ee
\be\label{equ2Nil}
\f{dB}{dt}=\f{4A}{C}-20\alpha'\f{A^{2}}{BC^{2}}
\ee
\be\label{equ3Nil}
\f{dC}{dt}=\f{4A}{B}-20\alpha'\f{A^{2}}{B^{2}C}
\ee

 \begin{figure}[htbp] 
\centering
\subfigure[$(A_0, B_{0},C_{0})=(7,5, 3), \alpha' = 1, T_s=0.134$]{\includegraphics[width=0.4\textwidth]{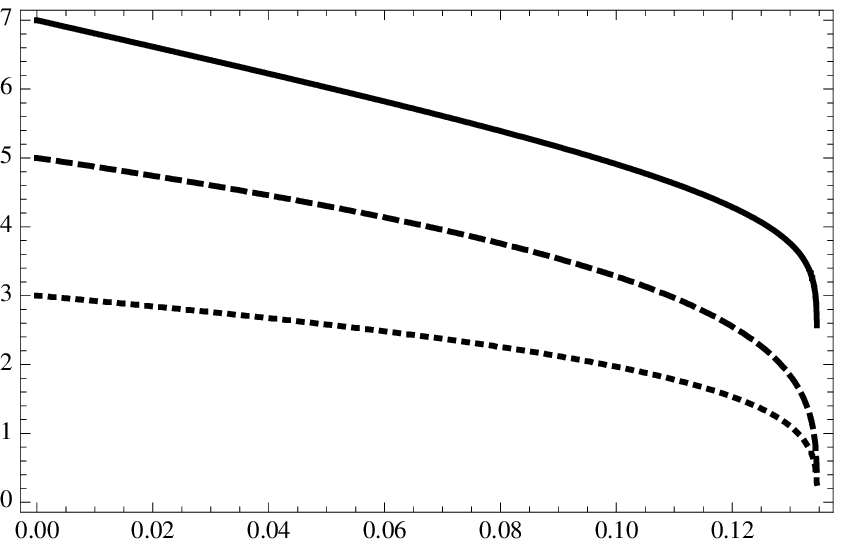}\label{2ndorderNilone} }
\subfigure[$(A_0, B_{0},C_{0})=(7,5, 3),\alpha' =-1,T_s=-0.043$]{\includegraphics[width=0.4\textwidth]{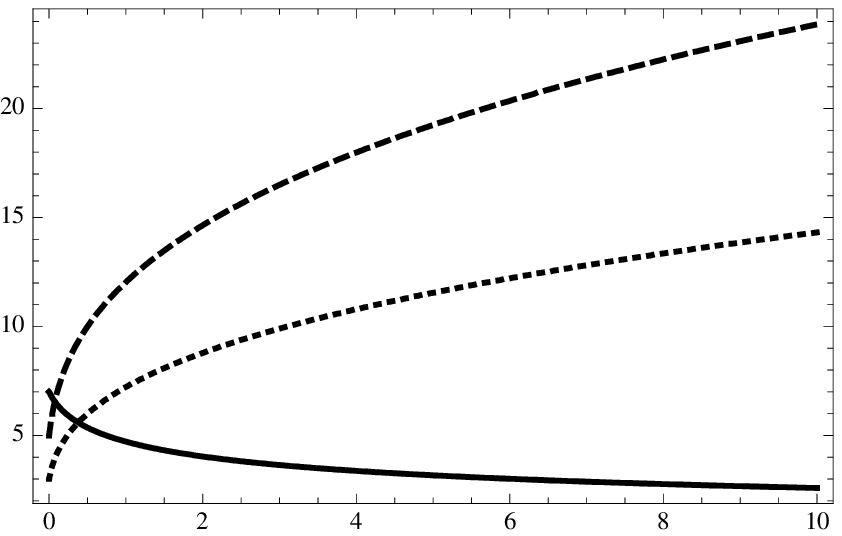}\label{2ndorderNiltwo} }
\subfigure[$(A_0, B_{0},C_{0})=(7,5, 3),\alpha' =0,T_s=-0.17$]{\includegraphics[width=0.4\textwidth]{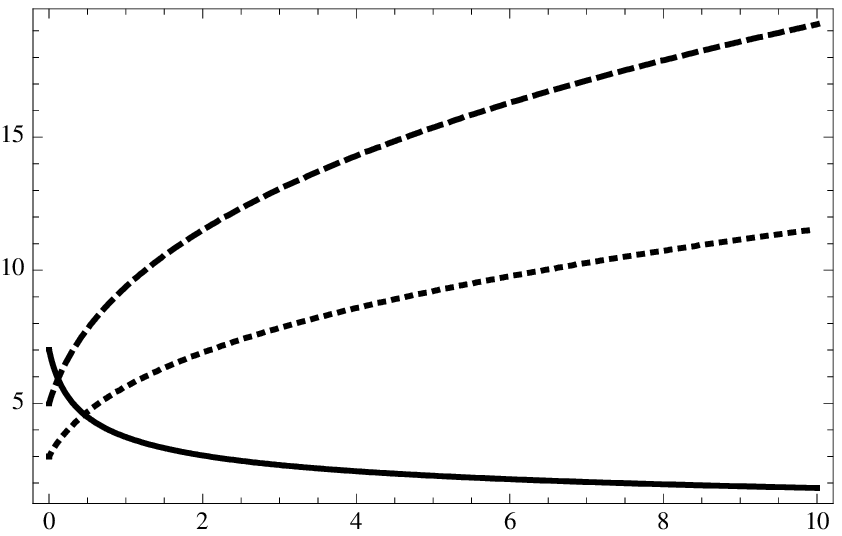}\label{2ndorderNilthree} }
\subfigure[$(A_0, B_{0},C_{0})=(3,5, 7),\alpha' =1$]{\includegraphics[width=0.4\textwidth]{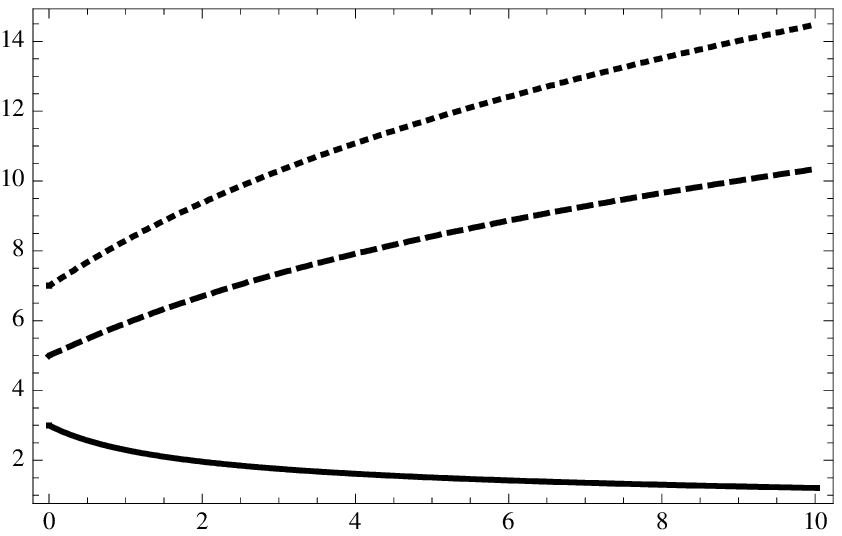} \label{2ndorderNilfour}}
\caption{$A(t),B(t),C(t)$ vs $t$ for  $\alpha'=1$and $\alpha'=0$ for Nil}
\label{2ndorderNil}
\end{figure}

\subsubsection{\bf Numerical estimates}

First we solve the flow equations numerically and look for
specific characteristics of the flow for different types of 
initial conditions. \comments{In Fig.\ref{2ndorderNil} the evolution of the scale 
factors has been depicted where $A(t),B(t)$ and $C(t)$ are plotted, as before,
in thick, 
dashed and dotted lines respectively 
} 
In Fig.\ref{2ndorderNilone} we see that for 
$\alpha' = 1$, $A_{0}\ge B_{0}\ge C_{0}$ and $A_{0},B_{0},C_{0}\ge 1$ the 
flow converges in every direction, contrary to 
 unnormalized Ricci flow (Fig.\ref{2ndorderNilthree}) where the
$B$, $C$ are expanding and diverge diverge while 
($A(t)$) converges i.e. a pancake degeneracy exists (the same happens 
for normalized Ricci flow as 
mentioned in\cite{isenberg}). 
Further, for $\alpha' = 0$ (Fig. 5(b)) and any 
value of $A(t),B(t),C(t)$, the flow has a past singularity.
For $\alpha' = 1$, if $ A_{0}\le B_{0}\le C_{0}$ the evolution of the scale factors resemble the unnormalized Ricci 
flow but have no singularity. On the other hand, for $\alpha'=-1$ (Fig.\ref{2ndorderNiltwo}) and any value of $A(t)$, $B(t)$, $C(t)$, the flow has a past singularity, similar to $\alpha'=0$ case.

\subsubsection{\bf Analytical solution}

It is possible to solve the flow equations analytically. Before
working out the solutions,
we note from the flow equations that $B/C=\text{constant}$. Let us now
choose a new variable $\xi = BC/A$. It can be shown easily \cite{glick2} 
that, with this
choice and appropriate scaling of the coordinates, the nature of the
flow can be determined entirely by finding the evolution of $\xi$. 
In other words, Nil has a one dimensional family of left invariant
metrics which is parametrised by $\xi$. 
The flow equations now take the form
	\begin{subequations}\label{nil gen xi}
		\begin{align}  			
			\dot{\xi} & = 12 - \alpha' ~\frac{36}{ \xi}\label{nil gen xione} \\[20pt] 
			\frac{\dot{A}}{A} & = - \frac{4}{\xi}\lb( 1 + \alpha' ~\frac{1}{ \xi} \rb)\label{nil gen xitwo}  \\[20pt]
			\frac{\dot{B}}{B} & = \frac{4}{\xi}\lb( 1 - \alpha' ~\frac{5}{ \xi} \rb)
		\end{align}
	\end{subequations}
where the $\xi$ equation replaces the $C$ equation and A, B are functions of
$\xi$. So, effectively, there is only one scale factor to worry about, i.e.
$\xi$.
The $\xi$ equation is readily solved to give two solutions-- 
	\begin{subequations}\label{nil soln xi}
		\begin{align}  			
			\xi + 3 \alpha' ~ \ln\lb\vert  \xi - 3\alpha'  \rb\vert & = 12 t + k \\
			& \mathrm{OR} \nonumber \\ 
			\xi & = 3 \alpha' 
		\end{align}
	\end{subequations}
The equations for $A$ and $B$ can then be solved to obtain $A(\xi)$ and $B(\xi)$given as, 
	\begin{subequations}\label{nil soln A B}
		\begin{align}  			
			\frac{A}{A_0} = \lb( \frac{\xi}{\xi_0} \rb)^\frac{1}{9}\lb( \frac{\4{\xi}{3} - \alpha'}{\4{\xi_0}{3} - \alpha'} \rb)^\frac{-4}{9} ~~~~~ & \& ~~~~~  \frac{B}{B_0} = \lb( \frac{\xi}{\xi_0} \rb)^\frac{1}{45}\lb( \frac{\4{\xi}{3} - \alpha'}{\4{\xi_0}{3} - \alpha'} \rb)^\frac{-2}{225} \\
						& \nonumber \\
			A = A_0 \exp \lb( -\frac{16}{3 \alpha'} t \rb) ~~~~~ &  \& ~~~~~ B = B_0 \exp \lb( -\frac{8}{3 \alpha'} t \rb)
		\end{align}
	\end{subequations}
where (49a) is for the $\xi$ given in (48a) and (49b) corresponds to the
solution (48b). The evolution of $C$ can also be found using $C=\frac{A\xi}{B}$.
In the case $\alpha'=1$, and $\xi= 3\alpha'$ we can easily see that  
the scale factors $A(t)$, $B(t)$ and $C(t)$ decrease.  When $\alpha'=0$, $\xi=12 t$ and $A(t)\sim t^{-\frac{1}{3}}$, 
$B(t)\sim t^{\frac{1}{3}}$, $C(t)\sim t^{\frac{1}{3}}$. 
For $\alpha'=-1$, the solution $\xi=3\alpha'$ is not valid as long as we
are dealing with Riemannian manifolds and the evolution is determined 
from the other solution.

Next we move on to a special case where $B=C$.
 
\subsubsection{\bf Special Case: $A\ne B= C $}
The flow equations, under this assumption are 
\begin{equation}\label{listeqn2Nil}
\left\{ \begin{array}{ll}
\f{dA}{dt}=-\f{4A^2}{B^{2}}-4\alpha'\f{A^3}{B^{4}}\\
\f{dB}{dt}=\f{4A}{B}-20\alpha'\f{A^{2}}{B^{3}}
\end{array} \right.
\end{equation}
It is obvious that the same analytical solutions mentioned above are
valid here as long as we define $\xi=\frac{B^2}{A}$. 
However, we discuss some alternative ways of arriving at the nature
of scale factor evolution.

Let us look at the difference between the scale factors given as
\be\label{}
\f{d(A - B)}{dt} = -4\f{A}{B^2}(A+B)-4\alpha' (\f{A}{B^2})^2 (A - 5B)\\
\le-4\f{A}{B^2}(A-B)-4\alpha' (\f{A}{B^2})^2 (A - 5B)
\ee
It is easily seen that, for $\alpha' = 0$ (i.e. unnormalized Ricci flow) 
$(A - B)$ decreases with increasing $t$, but for other values of $\alpha'$ 
different from zero the 
scale factors evolve differently. 
The flow equation for $B(t)$ is
\be
\4{dB}{dt}=\lb(4\4{A}{B^2}-20(\4{A}{B^2})^{2}\rb)B\label{nilBevolsp}
\ee
If $\4{B^2}{A}<5$ (or $-\4{A}{B^2}<-\4{1}{5}$) we have 
\be\label{nilspdynsys1}
\4{dB}{dt}<0
\ee 
which explicitly shows that $B(t)$ is decreasing in this region. 
We then move on to the case where $\4{B^2}{A}>5$. Assume $p= \4{A}{B^2}$.
The quantity inside the bracket in Eq.\ref{nilBevolsp} is a polynomial in $p$, 
$f(p)=4p-20p^2$. $f(p)$ represents a parabola (upside down) and its 
upper bound i.e $ max[{f(p)}]=\4{1}{5}$ occurs at $p=\4{1}{10}<\4{1}{5}$. 
The parabola intersects the abscissa at $p=0,\4{1}{5}$ and
progressively increases in the negative direction of the ordinate
(see Fig.\ref{2ndorderNilspparabola}). Thus, $f(p)$ is positive in the 
region $0<p<\4{1}{5}$ and, therefore $B(t)$ increases with increasing $t$.
Outside this domain of $t$, $B(t)$ decreases with increasing $t$. 
We now show the above-mentioned facts numerically. 
\begin{figure}[htbp] 
\centering
\subfigure[$f(\4{A}{B^2})=4(\4{A}{B^2})-20(\4{A}{B^2})^2$ and $(\4{A}{B^2})=\4{1}{3}$]{\includegraphics[width=0.4\textwidth]{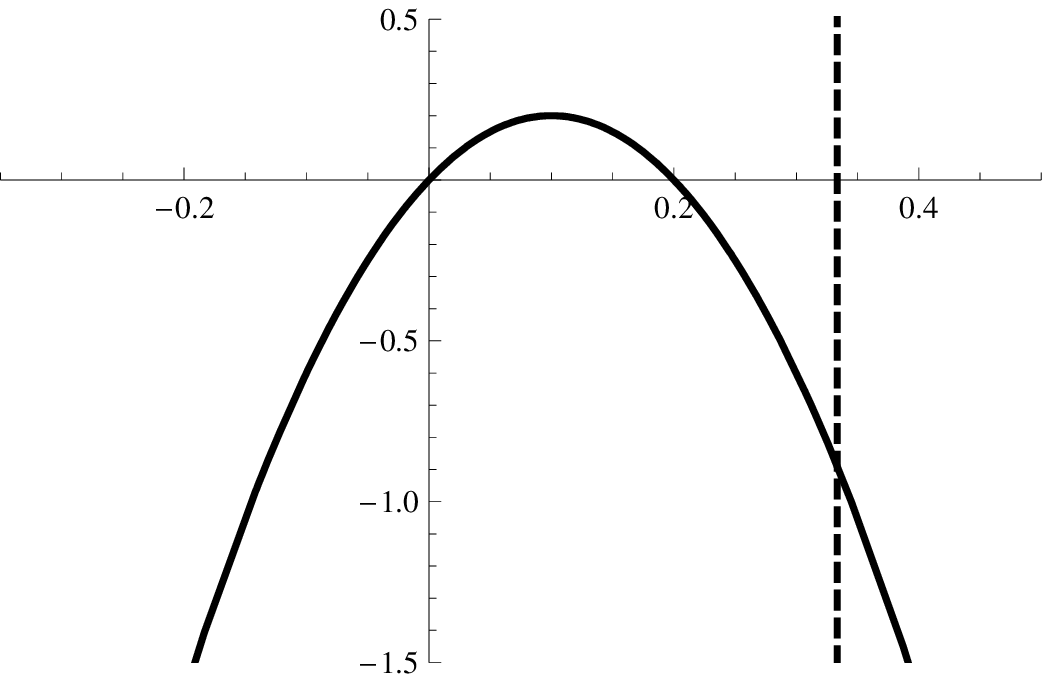} \label{2ndorderNilspparabola}}
\caption{$f(\4{A}{B^2})$ vs $(\4{A}{B^2})$ for  $\alpha'=1$ }
\label{2ndorderNilspparabola}
\end{figure}
In Fig[.\ref{2ndorderNilsp}] we have shown the behaviour of the scale 
factors for different genres of $\4{A}{B^2}$. The first figure 
(Fig.\ref{2ndorderNilspone}) has been plotted for such initial values of 
$A(t)$ and $B(t)$ which correspond to $\4{B^2}{A}<5$, when both the 
scale factors are decreasing (see  Eqn.[\ref{nilspdynsys1}], 
Fig. [\ref{2ndorderNilsptwo}]). In Fig.\ref{2ndorderNilsptwo} we have 
shown an interesting turning behavior where $B(t)$ initially decreases but
eventually increases. The initial condition has been chosen  to satisfy 
$\4{A}{B^2}>\4{1}{5}$ but during the evolution $A(t),B(t)$ decreases 
so that $\4{A}{B^2}$ eventually becomes less than $\4{1}{5}$ leading to 
growing nature of $B(t)$. On the other hand, for $\alpha'=-1,0$, 
the evolution is easy to comprehend from the equations (Eqn.[\ref{listeqn2Nil}]). The corresponding evolution of the scale factors is shown, respectively
 in Fig[\ref{2ndorderNilspthree}] and Fig[\ref{2ndorderNilspfour}].

\subsubsection{\bf Phase plots}

Let us now turn to the phase portraits(Fig.\ref{nil}).
 \begin{figure}[htbp] 
\centering
\subfigure[$(A_0, B_{0})=(7,3), \alpha' = 1,T_s=0.03$]{\includegraphics[width=0.4\textwidth]{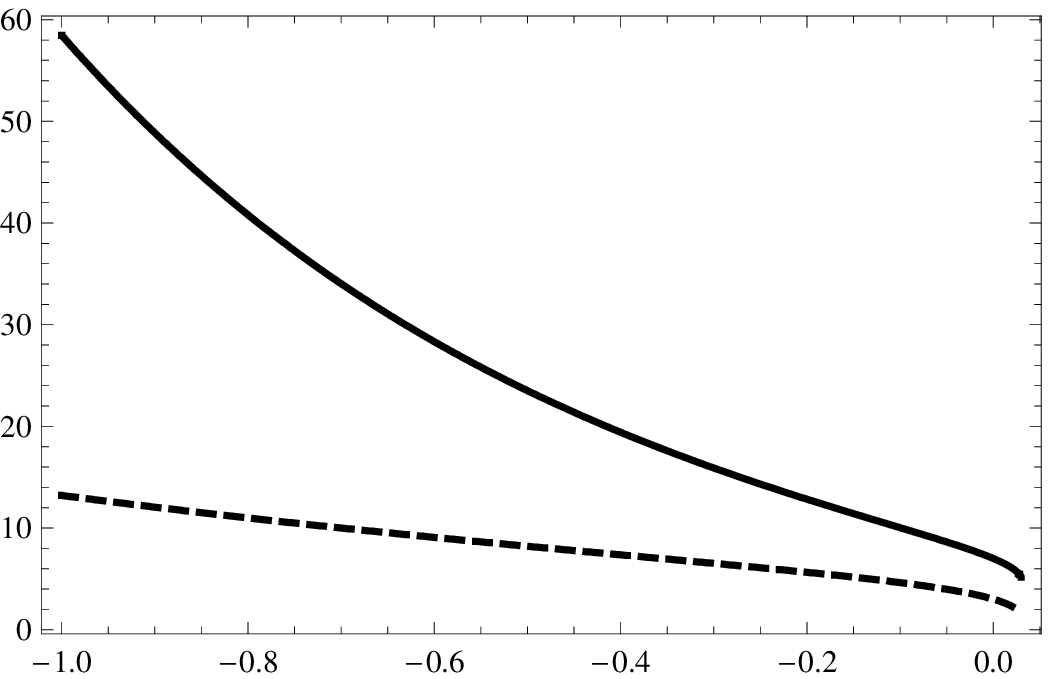} \label{2ndorderNilspone}}
\subfigure[$(A_0, B_{0})=(2,2.45), \alpha' = 1$]{\includegraphics[width=0.4\textwidth]{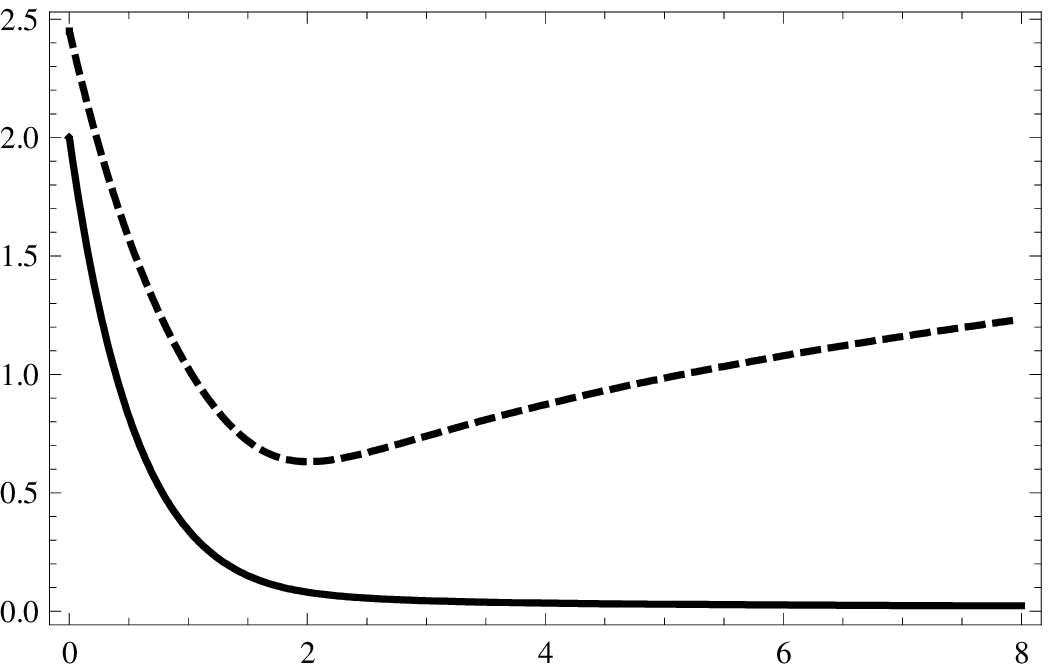}\label{2ndorderNilsptwo} }
\subfigure[$(A_0, B_{0})=(7,1),\alpha' =-1,T_s=-0.0003$]{\includegraphics[width=0.4\textwidth]{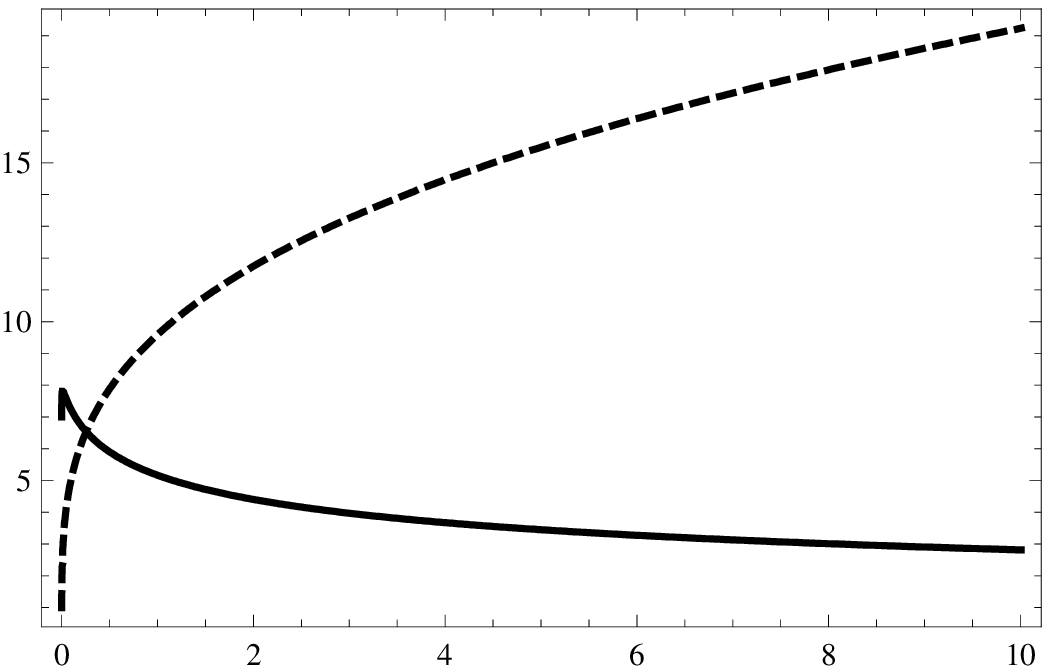}\label{2ndorderNilspthree} }
\subfigure[$(A_0, B_{0})=(5,3),\alpha' =0,T_s=-0.15$]{\includegraphics[width=0.4\textwidth]{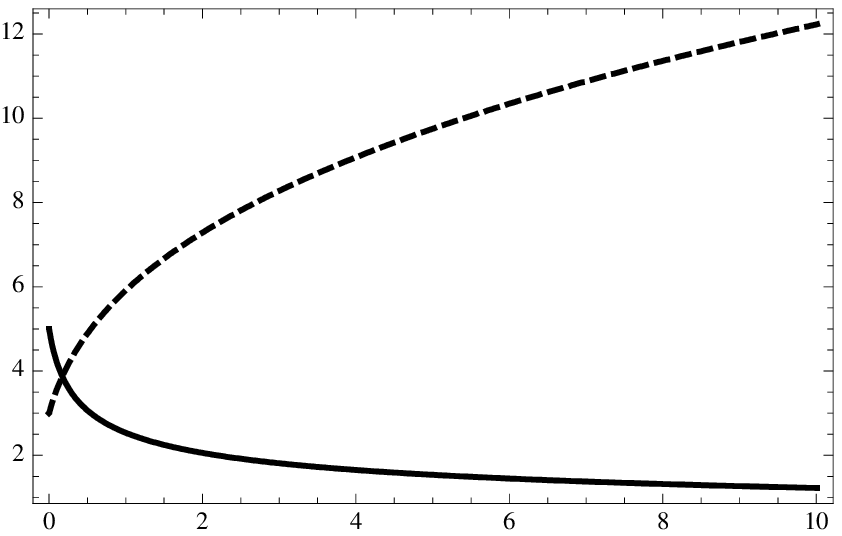}\label{2ndorderNilspfour} }
\caption{$A(t),B(t)$ vs $t$ for  $\alpha'=1$, $\alpha'=-1$ and $\alpha'=0$}
\label{2ndorderNilsp}
\end{figure}

The phase portraits in Fig.\ref{nil} are for $C=B$. From these plots and 
the equations, it is easy to see that $A=0$ is a fixed point of the 
system. For $\alpha' = 1$, trajectories from one region in phase space 
converge to the singularity $B=0$ whereas others move towards the 
fixed point $A=0$. 
The critical curve demarcating this behaviour is $\xi = BC/A = 3$ (plotted in dashed line). For $\alpha' = 0$ and $\alpha'=-1$ the flow show similar behavior with the trajectories proceed towards larger values of $B$ and $A=0$.

	\begin{figure}
	\centering
	\subfigure[$\alpha' = 1$]{\includegraphics[width=0.3\textwidth]{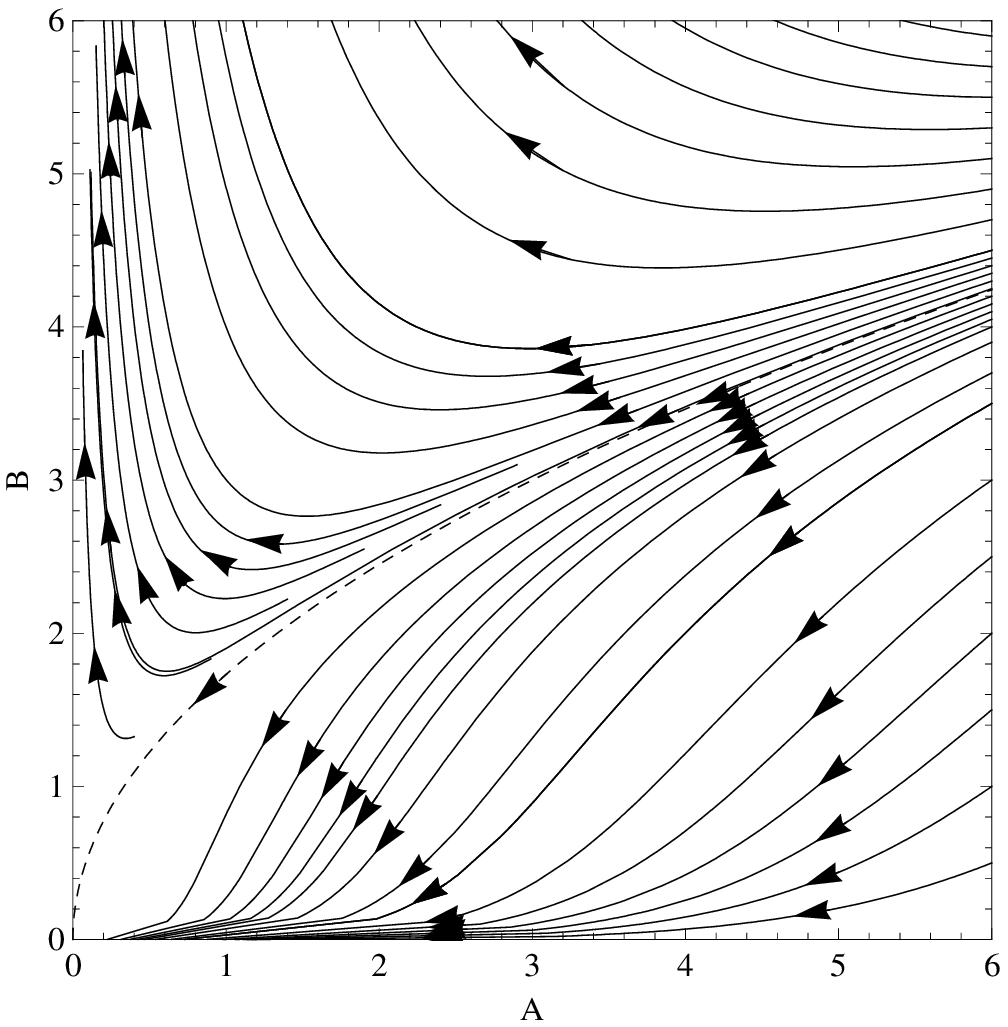}\label{nil 1}}
	\subfigure[$\alpha' = 0$]{\includegraphics[width=0.3\textwidth]{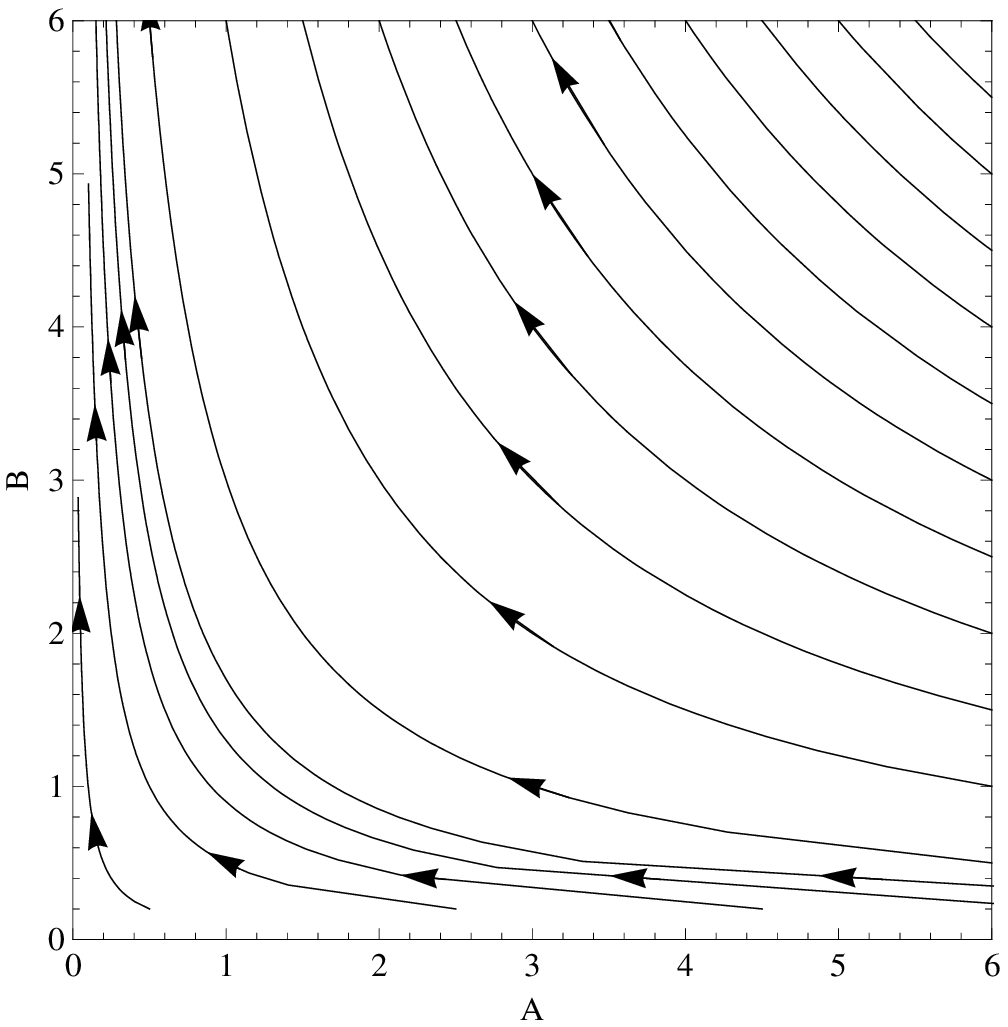}\label{nil 0}}
	\subfigure[$\alpha' = -1$]{\includegraphics[width=0.3\textwidth]{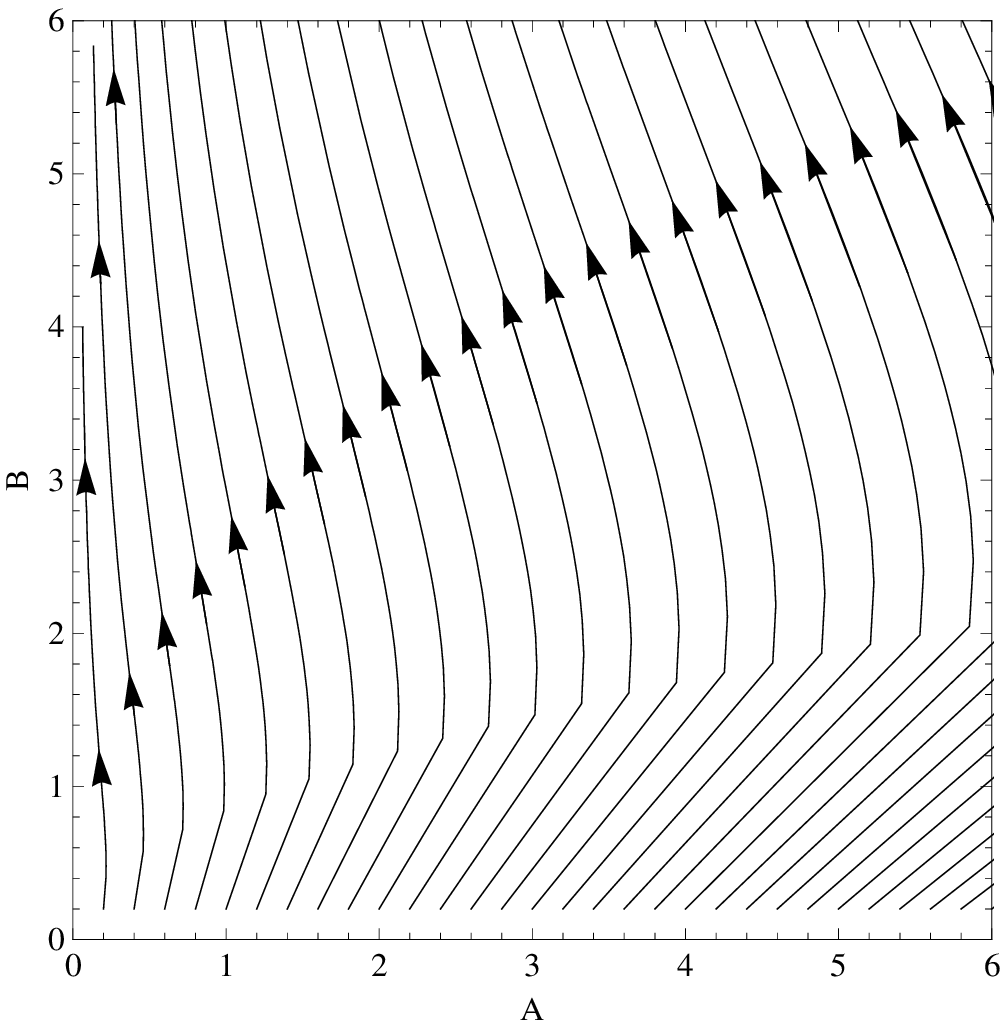}\label{nil -1}}
	\caption{$2nd$ order flow on Nil for $B=C$}
	\label{nil}
	\end{figure}

\subsubsection{\bf Curvature evolution}

\begin{figure}[h]
	\centering
	\subfigure[Evolution of scalar curvature for various $\alpha'$,  and ${A_0,B_0,C_0}$]{\includegraphics[width=0.62\textwidth]{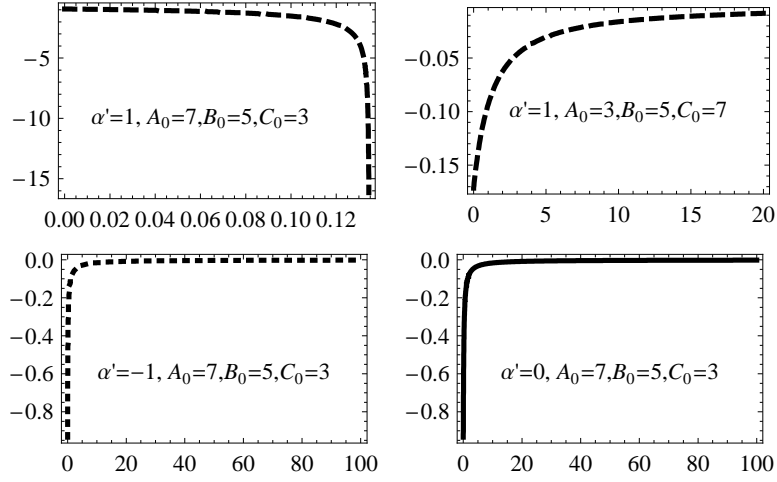}\label{curvevolnilone}}
	\caption{Evolution of scalar curvature for Nil.}
	\label{curvevolnil}
	\end{figure}

As before, here also we analyze the evolution of the scalar curvature. 
We have shown two different regimes of initial conditions, namely 
$A_0>B_0>C_0$ and $A_0<B_0<C_0$ in Fig.[\ref{curvevolnil}] for 
$\alpha'=1$, where we note the increasing or decreasing nature of the 
scalar curvature (for different initial conditions). 
For $\alpha'=0$ and $\alpha'=-1$ the scalar curvature monotonically 
converges to zero asymptotically after beginning from a negative value. 

We conclude this section by providing a comparison of the evolution of 
different scale factors and Scal. curvature for Nil manifolds in table\ref{table2chap03}.
	
 \begin{table}[h]\small
\caption{Comparisons of the evolution of different scale factors and Scal. curvature for Nil}\label{table2chap03}\centering
\renewcommand{\tabularxcolumn}[1]{>{\arraybackslash}m{#1}}
\begin{tabularx}{\textwidth}{+Y^Z^Z^Z^Z}
\toprule\rowstyle{\bfseries}
 & $A(t)$ & $B(t)$ & $C(t)$ & $Scal$\\
\otoprule
$\alpha'=0$  & decreasing & increasing & increasing & asymp. flat starting from -ve.\\
\midrule
$\alpha'=1$& decreasing & depends on i.c & depends on i.c & increase/decrease 
entirely -ve, past/future asymp. flat\\
\midrule
$\alpha'=-1$& decreasing, asymp. flat & increasing & increasing & asymp. flat starting from -ve.\\
\bottomrule
\end{tabularx}
\end{table}
	
\newpage
\subsection{\underline{Computation on Sol}}

\subsubsection{\bf Flow equations}

Sol is a solvable group which has minimum symmetry 
among all the eight geometries. The group action 
in ${\mathbb R}^3$ can be written as \\
\be
\lb(x,y,z \rb)\circ\lb(x',y',z' \rb)=\lb(x+e^{-z}x',y+e^{z}y',z+z' \rb)
\ee
We can put left invariant vector $\text{fields}^{\P}$  \footnotetext[5]{it is easy to check that one such left invariant metric will be $ds^2=e^{2z} dx^2+e^{-2z} dy^2+ dz^2$ } on the manifold, for which the structure constants in the 
Milnor frame will be $\lambda=-2,\mu=0,\nu=+2$. Sol really has a
two parameter family of metrics upto diffeomeorphism (obtained by
scaling and redefinition of the coordinates, see \cite{glick2}).
Using the values of the structure constants, we can find the components of Rc and ${\widehat{Rc}}_{}^{2}$ tensors. The $2nd$ order flow equations reduce to,
\be\label{equ1sol}
\f{dA}{dt}=\4{-4(A^{2}-C^{2})}{BC}-4\alpha'\4{(A+C)^{2}(A^{2}-2AC+5C^{2})}{AB^{2}C^{2}}
\ee
\be\label{equ2sol}
\f{dB}{dt}=4\4{(A+C)^{2}}{AC}-4\alpha'\4{(A+C)^{2}(5A^{2}-6AC+5C^{2})}{A^{2}BC^{2}}
\ee
\be\label{equ3sol}
\f{dC}{dt}=-4\4{(C^{2}-A^{2})}{AB}-4\alpha'\4{(A+C)^{2}(5A^{2}-2AC+C^{2})}{A^{2}B^{2}C}
\ee

\subsubsection{\bf Numerical and analytical estimates}

\begin{figure}[htbp] 
\centering
\subfigure[\small{$(A_0, B_{0},C_{0})=(7,5, 3), \alpha' = 1,T_s=0.078$}]{\includegraphics[width=0.4\textwidth]{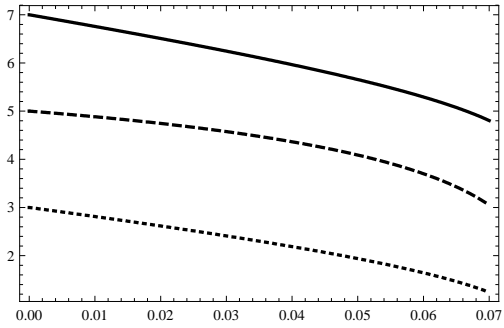}\label{2ndorderSolone} }
\subfigure[$(A_0, B_{0},C_{0})=(5,7, 3), \alpha' = 1$]{\includegraphics[width=0.4\textwidth]{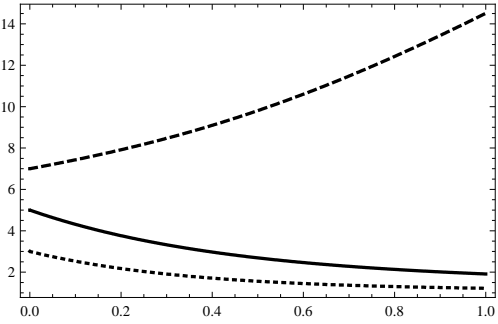}\label{2ndorderSoloneone} }
\subfigure[$(A_0, B_{0},C_{0})=(7,5, 3),\alpha' =-1,T_s=-0.052$]{\includegraphics[width=0.4\textwidth]{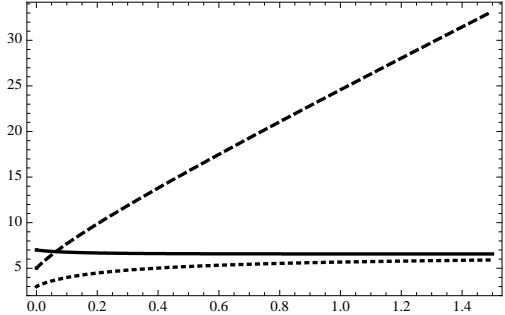}\label{2ndorderSoltwo} }
\subfigure[\small{$(A_0, B_{0},C_{0})=(7,5, 3),\alpha' =0,T_s=-0.154$}]{\includegraphics[width=0.4\textwidth]{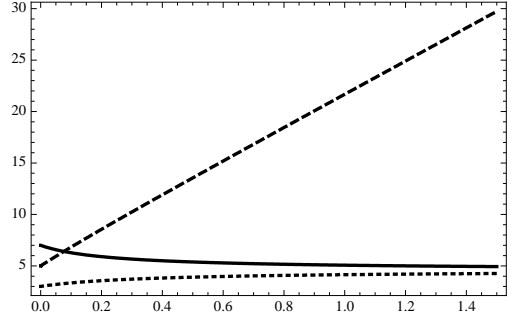} \label{2ndorderSolthree}}
\caption{$A(t),B(t),C(t)$ vs $t$ for  $\alpha'=1$, $\alpha'=-1$ and $\alpha'=0$ in Sol manifold}
\label{2ndorderSol}
\end{figure}

\noindent {\bf (a)} The behaviour of the scale factors for the
above higher order flow are different from 
those in unnormalized Ricci flow. From Fig.\ref{2ndorderSolthree} 
(unnormalized Ricci flow) we note that the scale factors 
do not have a future singularity. In fact, we see the appearance of a cigar degeneracy. 
The above behaviour can be predicted qualitatively from the equations 
themselves. We note that $A$ and $C$ can be interchanged in the equations. 
So instead of three equations it is enough to examine two of them. 
Further, without loss of generality, we may assume $A>C$. 
For $\alpha' = 1$ the scale factors converge (see Fig 10 (a)). 
If $\alpha'=0$ and $A>C$ we 
find that $B(t)$, $C(t)$ increase and $A(t)$ decreases (more detail for 
$\alpha'=0$ can be found in \cite{chowbook}). 
If $\alpha' = -1$ the behaviour of the scale factors bear a 
resemblance with unnormalized Ricci flow, though the singularity time 
(in the past)
changes due to the higher order term ( Fig.\ref{2ndorderSoltwo} illustrates 
one such example for a particular set of initial values). 

\noindent {\bf (b)}  We now consider the case $\alpha'=1$ and analyze the evolution of the 
scale factor $A(t)$ a bit further.
It can be shown that Eqn.\ref{equ1sol} can be written as
\be
\4{dA}{dt}\le -4\4{A}{C}\4{A}{B}\lb(1-\4{C}{A}\rb)^{2}-\4{2}{B}\4{A}{B}\lb(1-\4{C}{A}\rb)^{2}\lb(1+\4{A}{C}\rb)^{2}-\4{8}{B}\4{A}{B}\lb(1+\4{C}{A}\rb)^{2}
\ee
which shows that $A(t)$ is decreasing in forward time. 
This is also true for $C(t)$ where the second term in Eq.\ref{equ3sol} 
dominates over the first term and, therefore, the net effect is a 
decreasing $C(t)$. 
In the same way, we can argue that $B(t)$ is increasing (decreasing) 
for $A<B$ ($A>B$). 
Since $A(t)$ and $C(t)$ are both decreasing,
they may approach each other. To show this, consider 
the difference $A-C$, which satisfies 
\begin{subequations}
\begin{align}
\4{d(A-C)}{dt}& = -4\4{(A-C)}{ABC}(A+C)^2~-4\alpha' \4{(A+C)^2(A-C)}{A^2B^2C^2}(A^2-6AC+C^2)
\\ & \le-4\4{(A+C)^2}{ABC}C~-16\alpha' \4{(A+C)^2}{B^2C}
\end{align}
\end{subequations}
and demonstrates that the rate of change of $A-C$ is always negative. 
This feature is also clearly visible in Fig. 10(b). 
If $\alpha'=0$, we may note that $\4{dA}{dt}\le -4A$ or $A(t)\le (A_0-4t)$.

Next we move on to a special case where $A = C$, which is analytically
solvable. 


\subsubsection{\bf Special Case : $A= C\ne B $}
With this assumption, the flow equations simplify, and we have
\be\label{2equ1sol}
\f{dA}{dt}=-32\alpha'\4{A}{B^{2}}
\ee
\be\label{2equ2sol}
\f{dB}{dt}=8-32\alpha'\4{1}{B}
\ee

When $\alpha'=0$, $A(t)$ is a constant and $B(t)$ increases linearly.
If $\alpha'=1$, $A(t)$ decreases but $B(t)$ may increase/decrease
depending on initial conditions (as shown in Fig 11(a),(b)). 
For $\alpha'=-1$, A and B both increase in forward time (Fig. 11(c)).

The flow equations can be easily solved to obtain explicit solutions as --
 \begin{figure}[htbp] 
\centering
\subfigure[$(A_0, B_{0})=(3, 7), \alpha' = 1$]{\includegraphics[width=0.4\textwidth]{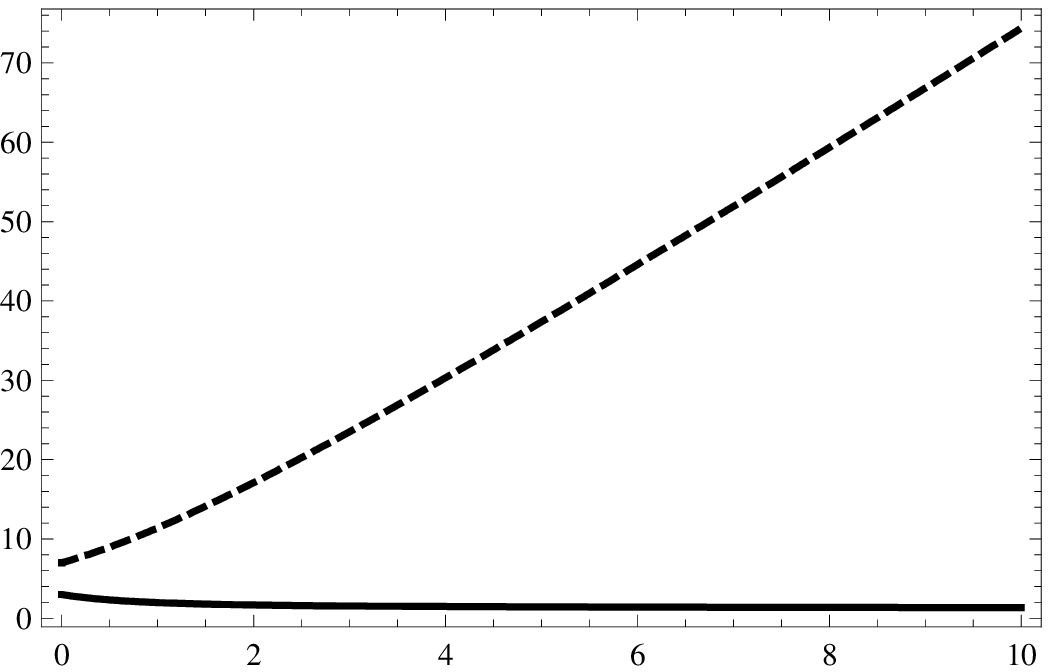} \label{2ndorderSolsp1}}
\subfigure[$(A_0, B_{0})=(7, 3), \alpha' = 1,T_s=0.318$]{\includegraphics[width=0.4\textwidth]{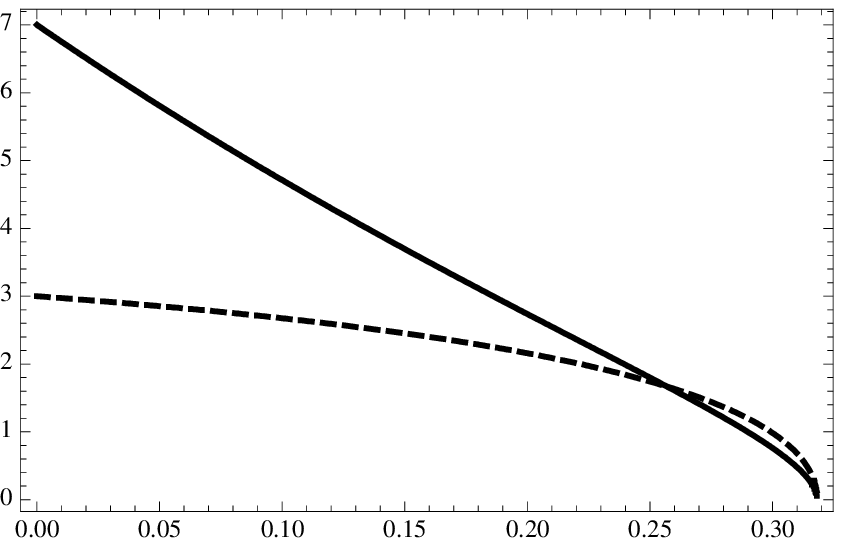} \label{2ndorderSolsp2}}
\subfigure[$(A_0, B_{0},C_{0})=(3,5, 7),\alpha' =-1,T_s=-0.219$]{\includegraphics[width=0.4\textwidth]{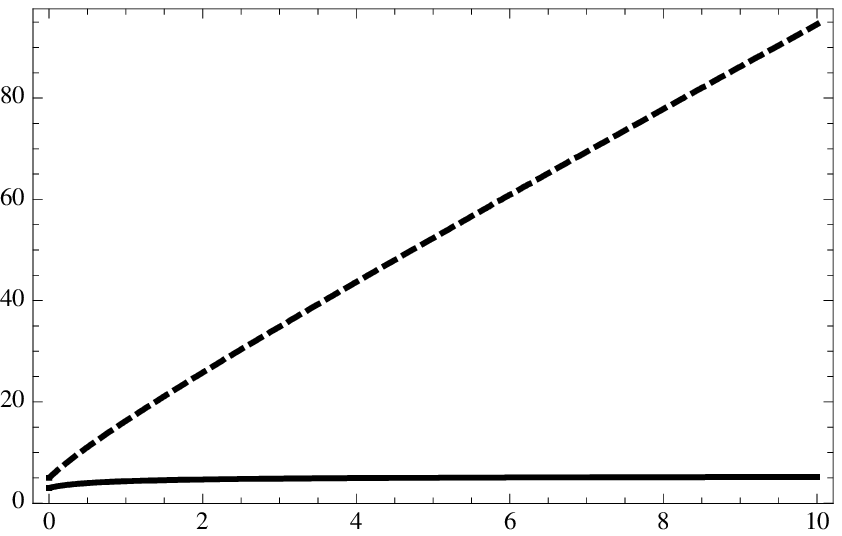} \label{2ndorderSolsp3}}
\caption{$A(t),B(t)$ vs $t$ for  $\alpha'=1$ and $\alpha'=-1$}
\label{2ndorderSolsp}
\end{figure}

	\begin{subequations}\label{sol special soln}
		\begin{align}  			
			A\lb(1- \frac{4\alpha'}{B} \rb) = \mathrm{constant} ~~~~~ & \& ~~~~~  B + 4\alpha' \ln\lb\vert B - 4\alpha' \rb\vert = 8t + k \\
			&\mathrm{OR} \nonumber \\
			A = A_0 \exp \lb( -\frac{2}{\alpha'} t \rb) ~~~~~ &  \& ~~~~~ B = 4\alpha'
		\end{align}
	\end{subequations}
Note that the second solution is not valid as long as we are
concerned with Riemannian manifolds.

\subsubsection{\bf Phase plots}
	
\begin{figure}
	\centering
	\subfigure[$\alpha' = 1$]{\includegraphics[width=0.3\textwidth]{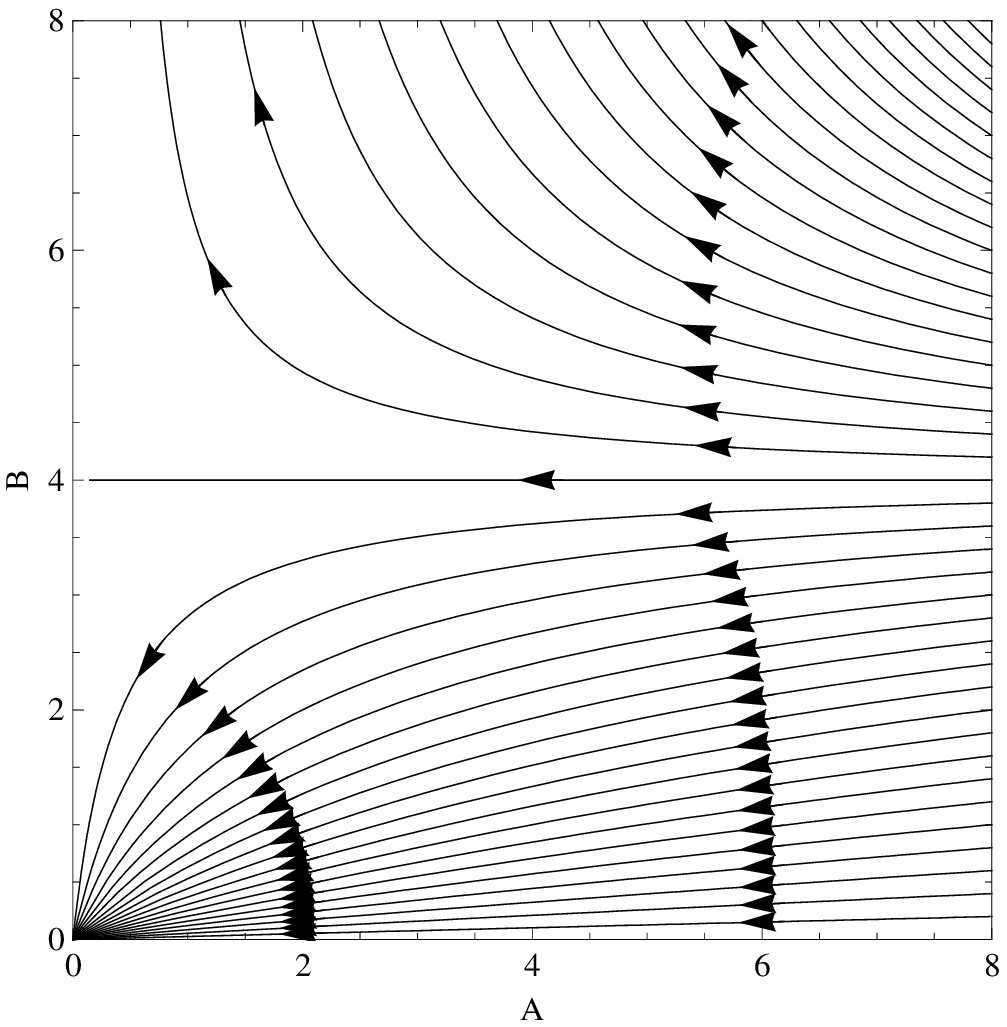}\label{sol 1}}
	\subfigure[$\alpha' = 0$]{\includegraphics[width=0.3\textwidth]{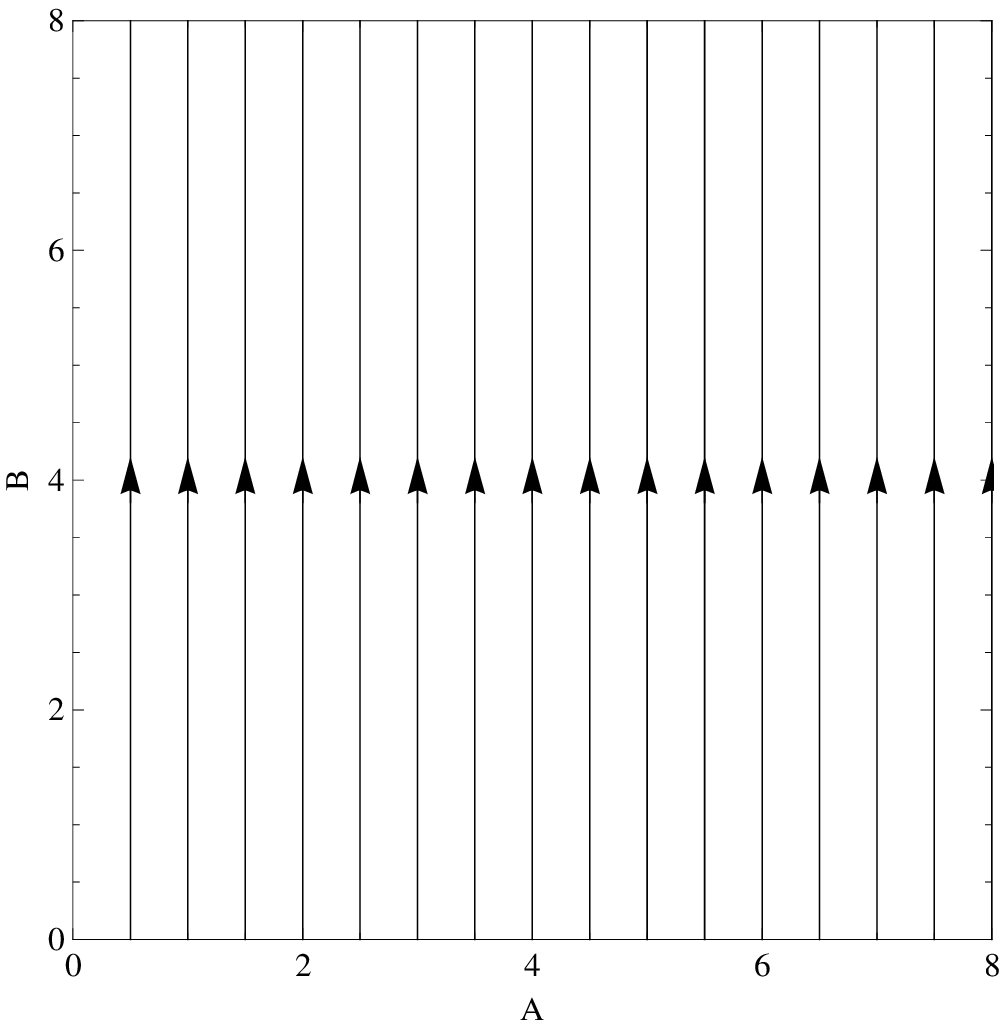}\label{sol 0}}
	\subfigure[$\alpha' = -1$]{\includegraphics[width=0.3\textwidth]{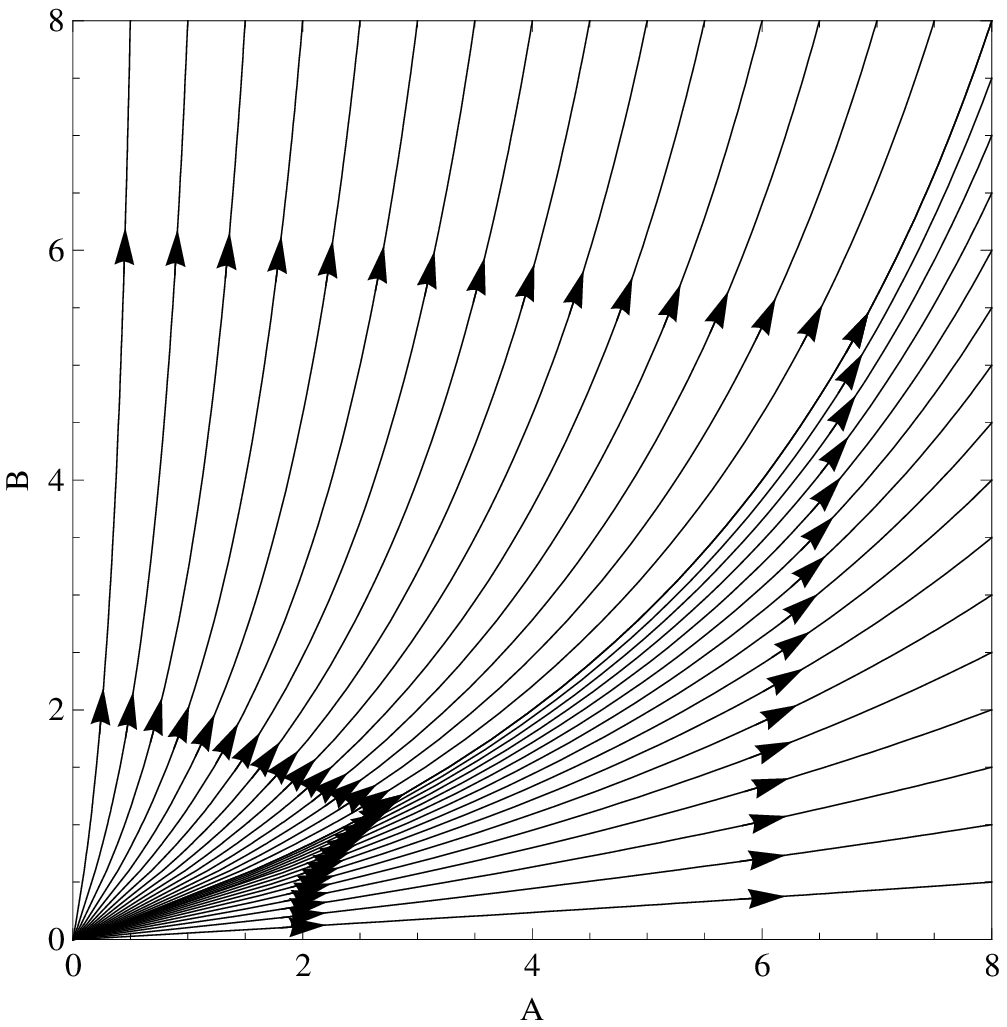}\label{sol -1}}
	\caption{$2nd$ order flow on Sol for $A=C$}
	\label{sol}
	\end{figure}

The phase space is plotted in Fig.\ref{sol}. For $\alpha' = 1$ flow 
trajectories from $B < 4 \alpha'$ converge to the singularity $B=0$ 
whereas others tend towards the fixed point $A=0$. $B= 4\alpha'$ is the 
critical curve in this case. For $\alpha' = 0$ flows we have $A=A_0$ and 
$B = 8t + B_0$, with the trajectories being straight lines flowing towards 
ever increasing values of $B$. In the $\alpha' = -1$ case the trajectories 
expand outwards from $A=B=0$ to larger and larger values.

\subsubsection{\bf Curvature evolution}

Lastly we examine the evolution of scalar curvature. For $\alpha'=0$ and 
$\alpha'=-1$ the behavior is  similar to that of Nil. The only 
difference appears in the $\alpha'=1$ case where we see that the curvature 
increases for $C_0<A_0<B_0$ while it decreases in the other regime.
\begin{figure}
	\centering
	\subfigure[Evolution of scalar curvature for various $\alpha'$, $A_0,B_0,C_0$]{\includegraphics[width=0.73\textwidth]{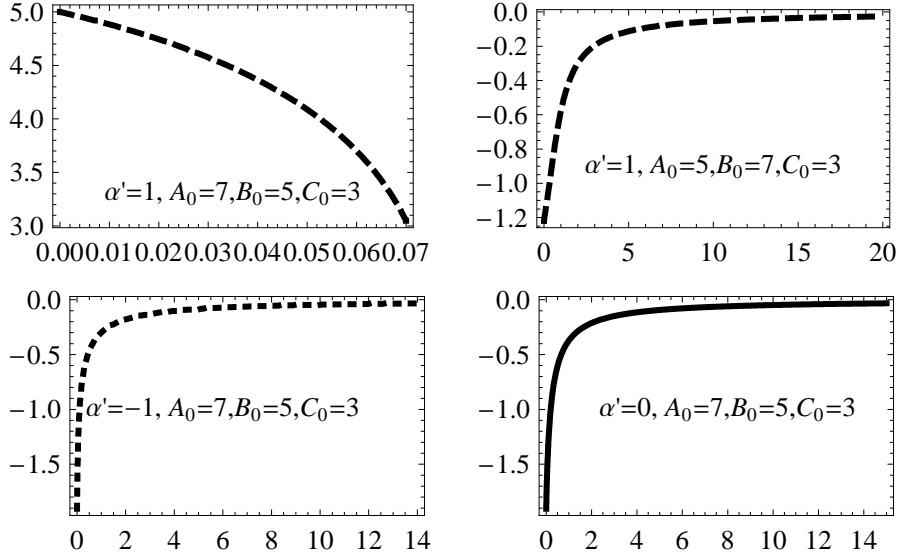}\label{}}
	\caption{Evolution of scalar curvature for Sol.}
	\label{curvevolsol}
	\end{figure}

\subsection{\underline{Computation on $\widetilde{\text{Isom}(\mathbb{R}^{2})}$ }}

\subsubsection{\bf Flow equations}

This case also involves a kind of solvable geometry. It is not
mentioned in Thurston's eight geometries because it is not 
maximally symmetric. According to Milnor's prescription we have the 
structure constants as follows $\lambda=-2,\mu=-2,\nu=0$. 
$\widetilde{\text{Isom}(\mathbb{R}^{2})}$ also, like Sol has
a two parameter family of metrics upto diffeomeorphism (see \cite{glick2}).
As done earlier, we can compute the various curvature quantitites 
and end up with the $2nd$ order flow equations given as
\be\label{equ1IsomR2}
\f{dA}{dt}=\4{-4(A^{2}-B^{2})}{BC}-4\alpha'\4{(A-B)^{2}(A^{2}+2AB+5B^{2})}{AB^{2}C^{2}}\
\ee
\be\label{equ2lsomR2}
\f{dB}{dt}=\4{-4(B^{2}-A^{2})}{AC}-4\alpha'\4{(A-B)^{2}(5A^{2}+2AB+B^{2})}{A^{2}BC^{2}}\\
\ee
\be\label{equ3lsomR2}
\f{dC}{dt}=\4{4(A - B)^{2}}{AB}-4\alpha'\4{(A-B)^{2}(5A^{2}+6AB+5B^{2})}{A^{2}B^{2}C}
\ee

\subsubsection{\bf Analytical estimates}

\noindent {\bf (a)} The $\widetilde{\text{Isom}(\mathbb{R}^{2})}$ class does admit Einstein 
metrics (flat metrics for this case) and converges under unnormalized Ricci 
flow(Fig.\ref{IsomRnullone}). It is easily observable that $A=B$ is the 
fixed point of the flow irrespective of the presence of 
the  higher order term.
For $\alpha' = 1$ each of the scale factors attain 
constant value asymptotically and they are seen to decrease initially. 

\begin{figure}[htbp] 
\centering
{\includegraphics[width=0.4\textwidth]{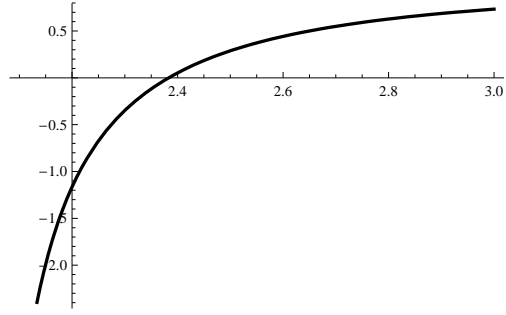}\label{} }
\caption{$F_1(n)$ vs $n$ }
\label{IsomRBevolufig}
\end{figure}

\noindent {\bf (b)} Let us check the evolution of each of the scale factors for $\alpha'=1$. 
Using the fact that $(A^2 + 2AB+5B^2)>(A+B)^2$ and Eqn.\ref{equ1IsomR2}  we can write the evolution of $A(t)$ as
\be
\4{dA}{dt}\le-4\4{(A-B)^2}{BC}-4\4{(A+B)^2 (A-B)^2}{AB^2C^2}
\ee
So the scale factor $A(t)$ always decreases. Using the same kind of argument as earlier $(5A^2 + 2AB+B^2)>(A+B)^2$ we can rewrite Eq.\ref{equ2lsomR2} as
\be
\4{dB}{dt}\le4\4{(A+B)^2}{AC}-4\4{(A+B)^2 (A-B)^2}{A^2BC^2}\\
\le 4\lb( 1-\4{(A-B)^2}{ABC}\rb)\4{(A-B)^2}{AC}
\ee
It can be easily shown that $\lb( 1-\4{(A-B)^2}{ABC}\rb)$(say $F_1$) can be 
negative as well as positive. To illustrate this we assume  
$A=n+2, B=n, C=n-2$ and calculate $F_1$ which turns out to be
$1-\frac{4}{(n-2) n (2+n)}$. From the plot of $F_1$
(Fig.\ref{IsomRBevolufig}), we note that depending on the value of $n$, 
$F_1$ can be 
positive as well as negative. Obviously this holds for $B(t)$ as well. 
Next, we estimate  the behavior of $\4{dC}{dt}$. From Eqn.\ref{equ3lsomR2} 
it is easy to anticipate that $(5A^2 + 6AB+5B^2)$ will be the deciding factor. 
Using the bound $-(5A^2 + 6AB+5B^2)<-(A-B)^2$ we can recast 
the Eqn.\ref{equ3lsomR2} as
\be
\4{dC}{dt}\le 4(\4{A}{B}-1)^2\lb(\4{B}{A}-\4{16B}{AC}\rb)
\ee
The behaviour for $C(t)$ may be found by choosing A, B, C and calculating
a quantity similar to the $F_1$ mentioned earlier. The evolution of 
$C(t)$ is qualitatively similar to that of $B(t)$.  

\subsubsection{\bf Numerical estimates}

\begin{figure}[htbp] 
\centering
\subfigure[$(A_0, B_{0},C_{0})=(5,8, 1), \alpha' = 1$]{\includegraphics[width=0.4\textwidth]{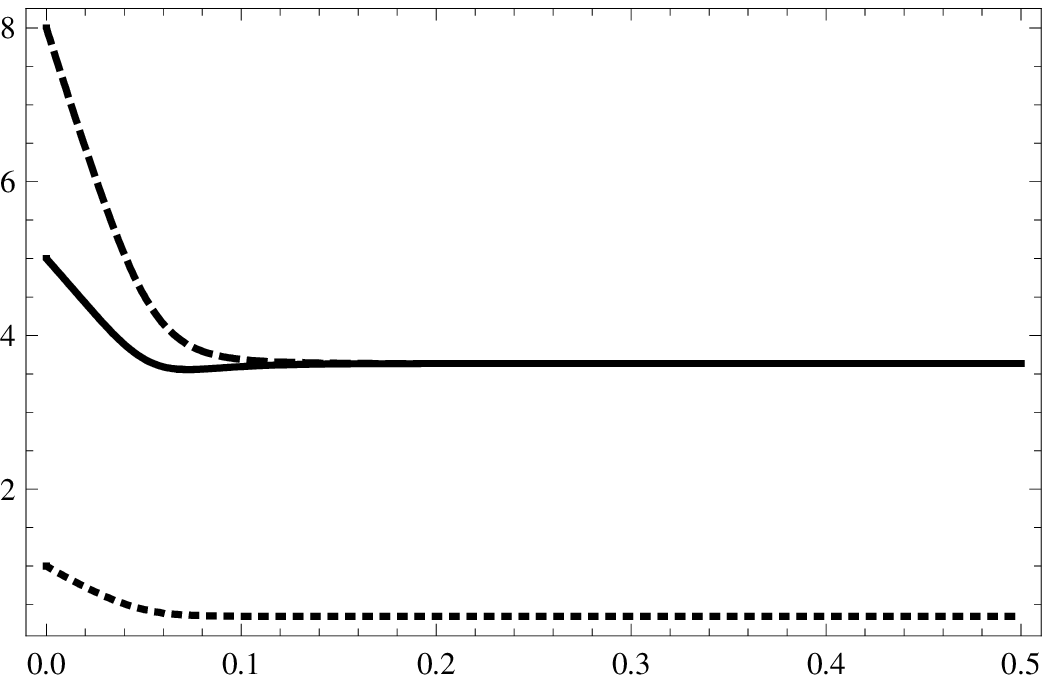}\label{2ndorderIsomRone} }
\subfigure[$(A_0, B_{0},C_{0})=(5,8, 0.8), T_s=0.028,\alpha' = 1$]{\includegraphics[width=0.4\textwidth]{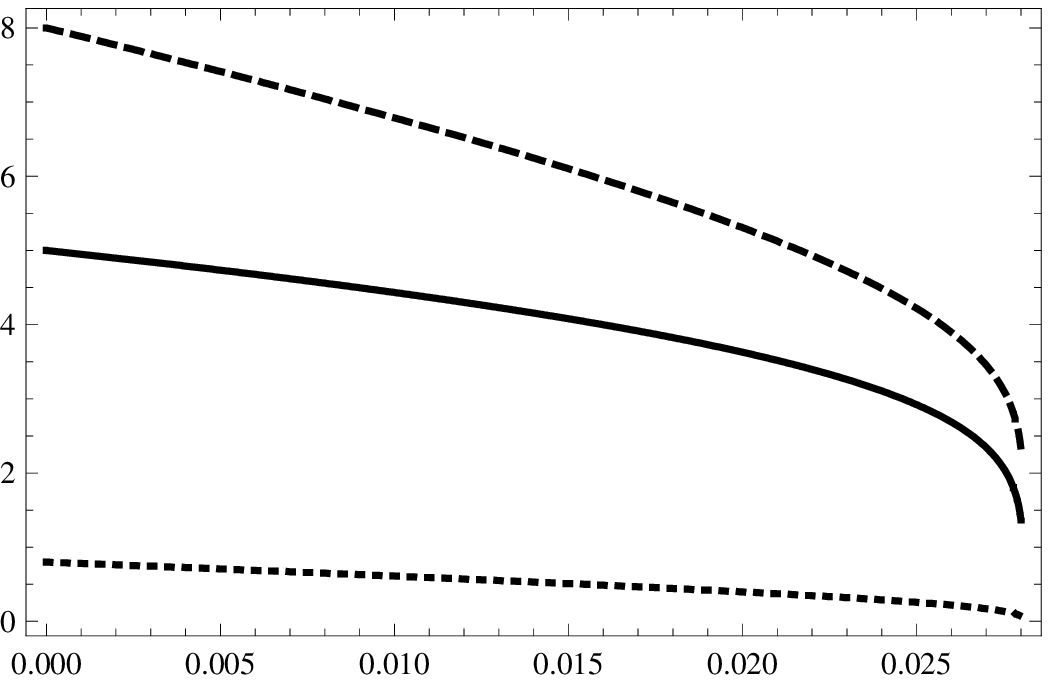}\label{2ndorderIsomRone} }
\subfigure[$(A_0, B_{0},C_{0})=(7,5,3), \alpha' = 1$]{\includegraphics[width=0.4\textwidth]{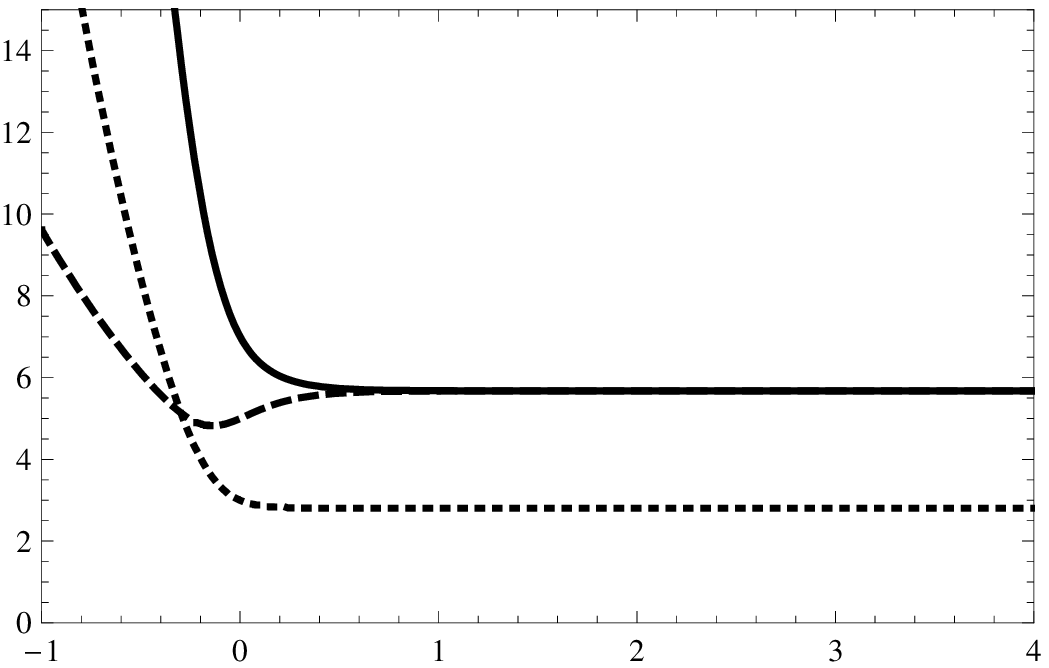}\label{2ndorderIsomRone} }
\subfigure[$(A_0, B_{0},C_{0})=(7,5, 3),\alpha' =-1$,$T_s = -0.128$]{\includegraphics[width=0.4\textwidth]{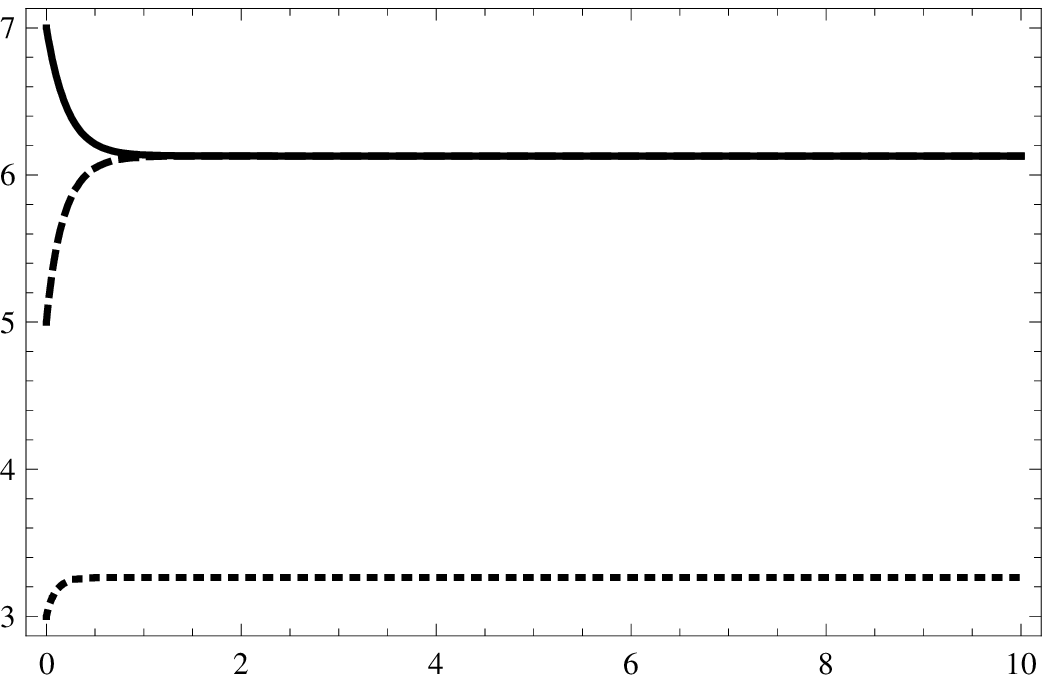}\label{2ndorderIsomRtwo} }
\caption{$A(t),B(t),C(t)$ vs $t$ for  $\alpha'=1$ and $\alpha'=-1$ in $\widetilde{\text{Isom}(\mathbb{R}^{2})}$}
\label{2ndorderIsomR}
\end{figure}

\begin{figure}[htbp] 
\centering
\subfigure[$(A_0, B_{0},C_{0})=(7,5, 3),\alpha'=0,T_s=-0.3$]{\includegraphics[width=0.4\textwidth]{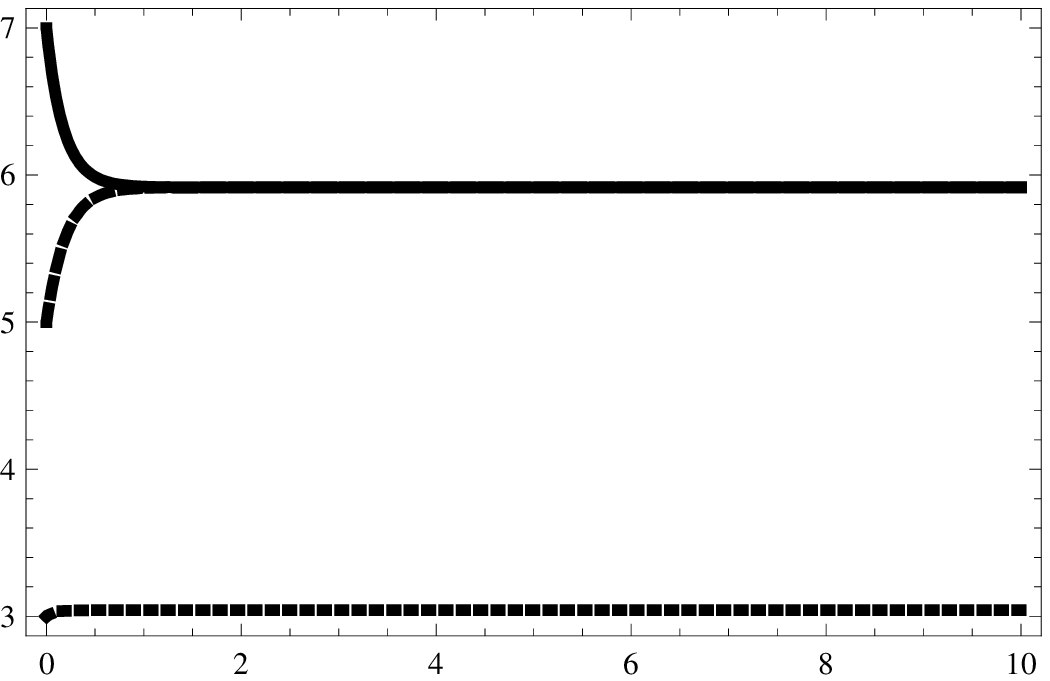}\label{IsomRnullone} }
\subfigure[comparison, $\alpha' =-1$, $\alpha' = 0$(thick)]{\includegraphics[width=0.4\textwidth]{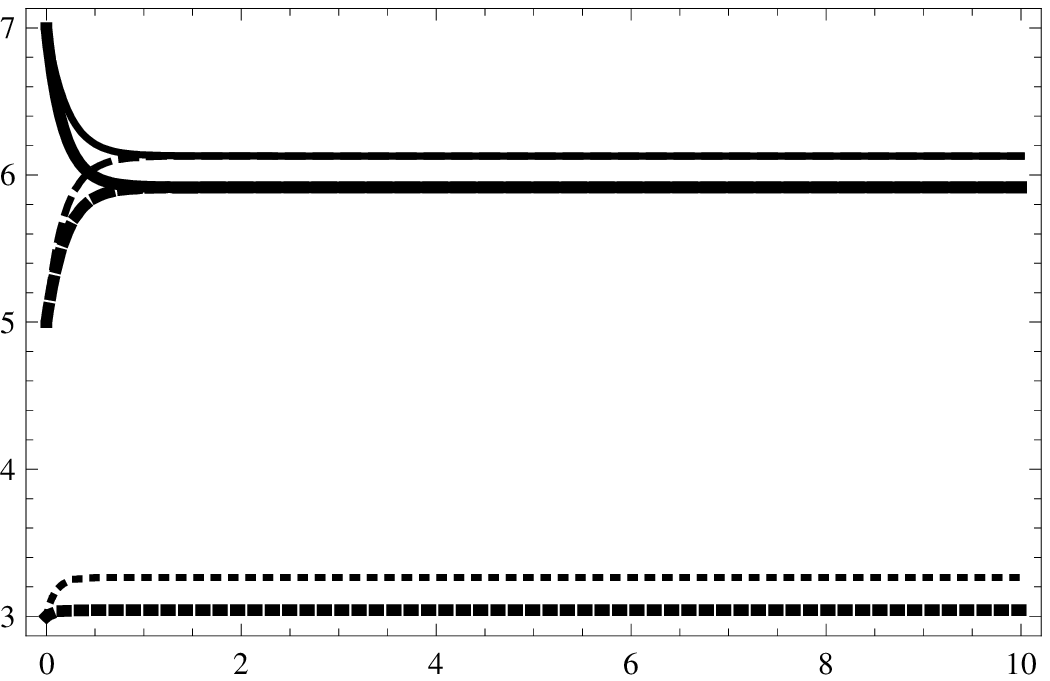} \label{IsomRnulltwo} }
\caption{$A(t),B(t),C(t)$ vs $t$ for $\alpha' = 0$ in $\widetilde{\text{Isom}(\mathbb{R}^{2})}$}
\label{IsomRnullcompar}
\end{figure}

In Fig.[\ref{2ndorderIsomR}] we demonstrate our conclusions for certain 
specific inital values by numericaly solving the dynamical system.
The cases with $\alpha'=1$ and different sets of initial values
are shown   
in Fig. 15(a)-(c). For $\alpha' = -1$ the flow develops a 
past singularity (Fig.\ref{2ndorderIsomRtwo}). The nature of the evolution of scale factors for $\alpha'=0$ is similar to the case $\alpha'=-1$ except for singularity time. These features appear in Fig. \ref{IsomRnullone} and Fig.\ref{IsomRnulltwo}.


\subsubsection{{\bf Special case:} $B = \eta A$}
Let us now choose $B = \eta A$ with $0< \eta < 1$. The flow equations become --

	\begin{subequations}\label{isom special}
		\begin{align}  			
			\f{dA}{dt}   & = 4\frac{A \left(\eta ^2-1\right)}{C \eta } -4 \alpha' ~\frac{A (\eta -1)^2 ( 5 \eta^2 +2 \eta + 1)}{C^2 \eta ^2} \\[20pt]
			\f{d\eta}{dt}& = -8 ~\frac{\eta ^2 - 1}{C} + 16\alpha' \frac{(\eta +1) (\eta -1)^3}{C^2 \eta } \\[20pt]
 			\f{dC}{dt}   & = 4\frac{\lb( \eta - 1\rb)^2}{\eta}  -4\alpha' ~\frac{(\eta -1)^2 ( 5 \eta^2 +6\eta + 5)}{C \eta ^2}
		\end{align}
	\end{subequations}

\begin{figure}
	\centering
	\subfigure[$\alpha' = 1$]{\includegraphics[width=0.3\textwidth]{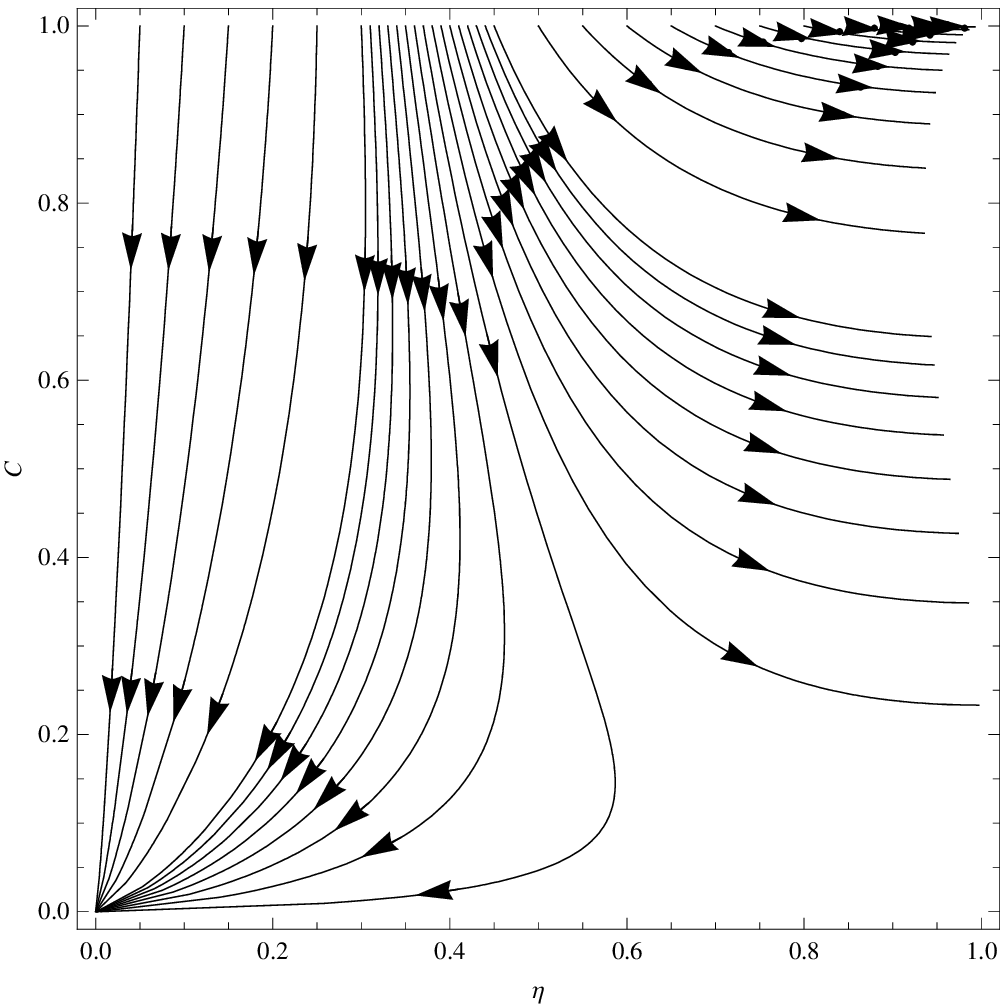}\label{isom 1}}
	\subfigure[$\alpha' = 0$]{\includegraphics[width=0.3\textwidth]{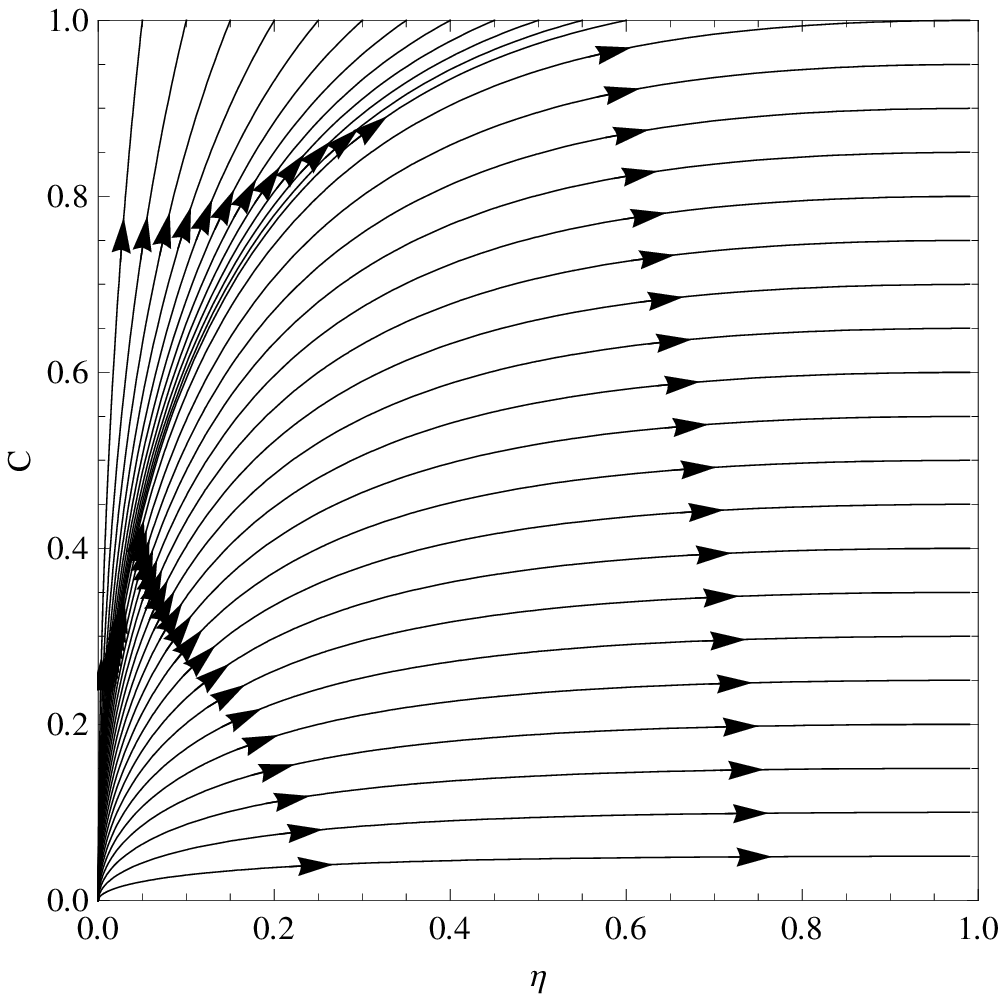}\label{isom 0}}
	\subfigure[$\alpha' = -1$]{\includegraphics[width=0.3\textwidth]{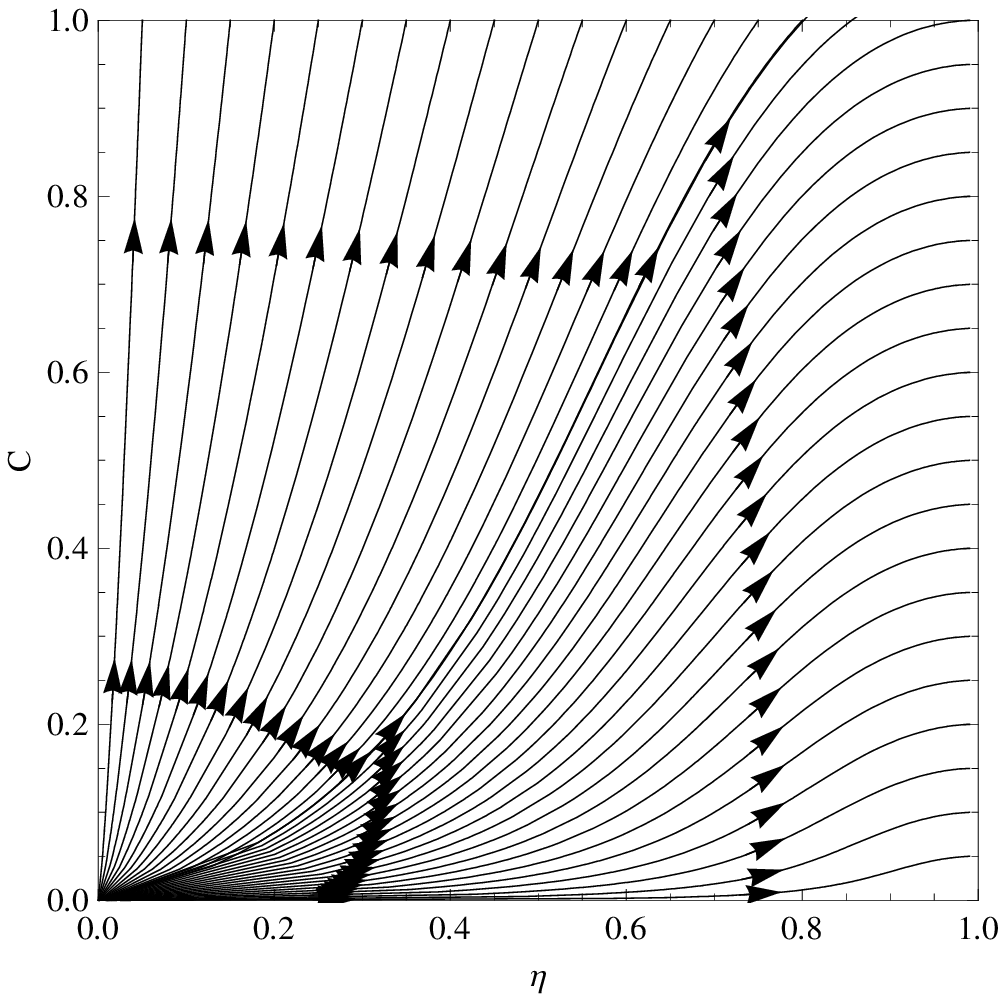}\label{isom -1}}
	\caption{$2nd$ order flow on $\widetilde{\text{Isom}(\mathbb{R}^{2})}$}
	\label{isom}
	\end{figure}

For $\alpha'=0$ we can exactly  solve the  $\eta - C$ system and obtain a relation between $\eta(=\4{B}{A})$ and C which is given as \\
\be
A(t)= B(t)\4{\lb( k-\sqrt{k^2-4C(t)^2}\rb)}{\lb( k+\sqrt{k^2-4C(t)^2}\rb)}
\ee
We substitute the above relation back in Eqn.\ref{isom special}b and 
arrive at an exact solution for $\eta(t)$ as\\
\be
\4{2\sqrt{\eta}}{(1+\eta)}+ \ln\4{(\sqrt{\eta}-1)}{(\sqrt{\eta}+1)}=-\4{32t}{k}+ C_{1}
\ee
We have analyzed numerically the   $\eta - C$ system in Fig.\ref{isom}. We note that $\eta = 1$ is a fixed line and attracts flow trajectories in all cases. For $\alpha' = 0,-1$ all trajectories in phase space end up on $\eta =1$. 
However, for $\alpha' = 1$ some flow to $\eta = 1$ whereas others 
flow towards the singularity at $\eta = 0 = C$.

\subsubsection{\bf Curvature evolution}

\begin{figure}
	\centering
	\subfigure[Evolution of scalar curvature for various $\alpha'$, $A_0,B_0,C_0$]{\includegraphics[width=0.71\textwidth]{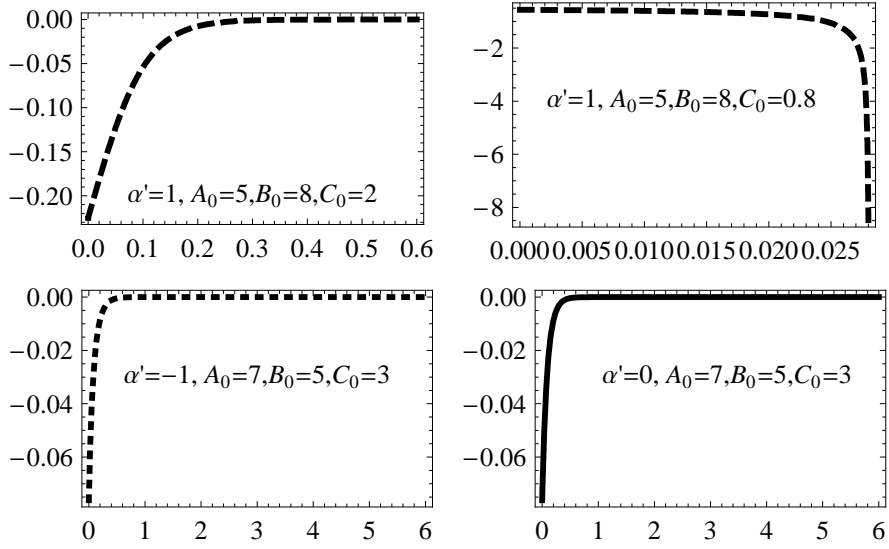}\label{}}
	\caption{Evolution of Scalar curvature for $\widetilde{\text{Isom}(\mathbb{R}^{2})}$.}
	\label{curvevolisom}
	\end{figure}

Before we end our discussion, it is worthwhile to mention 
the evolution of scalar curvature for this case. Like the Sol and Nil 
manifolds, the scalar curvature, for $\alpha'=0,-1$,
asymptotically reaches zero but 
by starting from a negative value. The case  $\alpha'=1$ has 
different behavior depending on the initial conditions, as seen in the 
Fig.[\ref{curvevolisom}]. For some particular initial values, the 
scalar curvature 
starts from a negative value, decreases first, then 
asymptotically reaches zero. For the same $\alpha'=+1$, but different 
initial conditions, we obtain negative scalar curvature as
before,  but it grows to larger negative values till the flow exists.
	
\subsection{\underline{Computation on ${\bf \widetilde{\text{SL}(2,\mathbb{R)}}}$}}

\subsubsection{\bf Flow equations}

This is a three dimensional Lie group consisting of all $2\times 2$ matrices 
with unit determinant and  its universal cover is denoted by 
$\widetilde{\text{SL}(2,\mathbb{R})}$. It can be shown that $\widetilde{\text{SL}(2,\mathbb{R})}$ is not isometric to $\mathbb{H}^2\times \mathbb{R}$ which is 
also a locally homogeneous space. We can put a left invariant metric\ddag\footnotetext[3]{one such set of left invariant forms are 
\be
\left\{ \begin{array}{ll}
\pi^1=dx+(1+x^2)^2 dy + (x-y-x^2 y) dz \\
\pi^2=2 x dy+(1-2 x y) dz\\
\pi^3=dx+(x^2 - 1) dy + (x + y - x^2 y) dz\nonumber\end{array} \right.
\ee
} 
(or right invariant metric but not bi-invariant) on this manifold. 
Using Milnor's prescription $\lb(\lambda=-2,\mu=-2,\nu=+2 \rb)$ we
find the various curvature-related quantities.
The evaluated  
components of Rc and ${\widehat{Rc}}_{}^{2}$ eventually lead to the
following flow equations:\\

\comments{\be\label{listRcSL2R}
\left\{ \begin{array}{ll}
Rc(F_{1},F_{1})=\4{2A^{2}-2(B+C)^{2}}{BC}\\
Rc(F_{2},F_{2})=\4{2B^{2}-2(A+C)^{2}}{AC}\\
Rc(F_{3},F_{3})=\4{2C^{2}-2(A-B)^{2}}{AB}\\
\end{array} \right.
\ee

\be\label{list2ndRcSL2R}
\left\{ \begin{array}{ll}
{\widehat Rc}(F_{1},F_{1})= \4{2(A^{4} + 2 A^{2} (B + C)^{2} - 
   8 A (B - C) (B + C)^{2} + (B + C)^{2} (5 B^{2} - 6 B C + 5 C^{2}))}{B^{2}C^{2}A}\\
{\widehat Rc}(F_{2},F_{2})=\4{2 (5 A^4 + 4 A C (B + C)^2 + A^3 (-8B + 4 C) + 
   2 A^2 (B^2 - 4 BC - C^2) + (B + C)^2 (B^2 - 2 BC + 5 C^2))}{A^2C^2B}\\
{\widehat Rc}(F_{3},F_{3})=\4{2 (5 A^4 - 4 A^3 (B - 2 C) - 4 A B (B + C)^2 - 
   2 A^2 (B^2 + 4 BC - C^2) + (B + C)^2 (5 B^2 - 2 BC + C^2))}{A^2B^2C}
\end{array} \right.
\ee
}
\be\label{SL2RevolAfull}
\begin{split}
\4{dA}{dt}= -4\4{A^{2}-(B+C)^{2}}{BC}-2\alpha'\4{2(A^{4} + 2 A^{2} (B + C)^{2} - 
   8 A (B - C) (B + C)^{2})}{B^{2}C^{2}A}\\
   -2\alpha' \4{2(B + C)^{2} (5 B^{2} - 6 B C + 5 C^{2})}{B^{2}C^{2}A}
\end{split}
\ee
\be\label{SL2RevolBfull}
\begin{split}
\4{dB}{dt}= -4\4{B^{2}-(A+C)^{2}}{AC}-2\alpha'\4{2 (5 A^4 + 4 A C (B + C)^2 + A^3 (-8B + 4 C))}{A^2C^2B}\\
-2\alpha'\4{ 2(2 A^2 (B^2 - 4 BC - C^2) + (B + C)^2 (B^2 - 2 BC + 5 C^2))}{A^2C^2B}
\end{split}
\ee
\be\label{SL2RevolCfull}
\begin{split}
\4{dC}{dt} = -2\4{2C^{2}-2(A-B)^{2}}{AB}-2\alpha'\4{2 (5 A^4 - 4 A^3 (B - 2 C) - 4 A B (B + C)^2)}{A^2B^2C}\\ 
   -2\alpha'\4{2( -2 A^2 (B^2 + 4 BC - C^2)+(B + C)^2 (5 B^2 - 2 BC + C^2))}{A^2B^2C}
\end{split}
\ee

\subsubsection{\bf Analytical and numerical estimates}

\noindent {\bf (a)} The abovementioned equations
are very hard to solve, analytically, when $\alpha'\ne0$. For $\alpha'=0$ 
we can see from Eqn.\ref{SL2RevolAfull}  and Eqn.\ref{SL2RevolBfull} that the 
flow is symmetric in $A$ and $B$. Assuming $A_0\ge B_0$, 
without any loss of generality, we can show that $A(t)\ge B(t)$ 
throughout the $\alpha'=0$ flow. This may be inferred (for $A>B>C$) 
from the evolution of the 
difference of $(A-B)$ given as\\
\be\label{evola-bsl2rURF1}
\4{d(A-B)}{dt}=-\4{4}{ABC}\lb(A-B\rb)\lb(A+B-C\rb)\lb(A+B+C\rb) \le 0
\ee

\begin{figure}[htbp] 
\centering
\subfigure[$(A_0, B_{0},C_{0})=(9,7,5), \alpha' = 1$]{\includegraphics[width=0.4\textwidth]{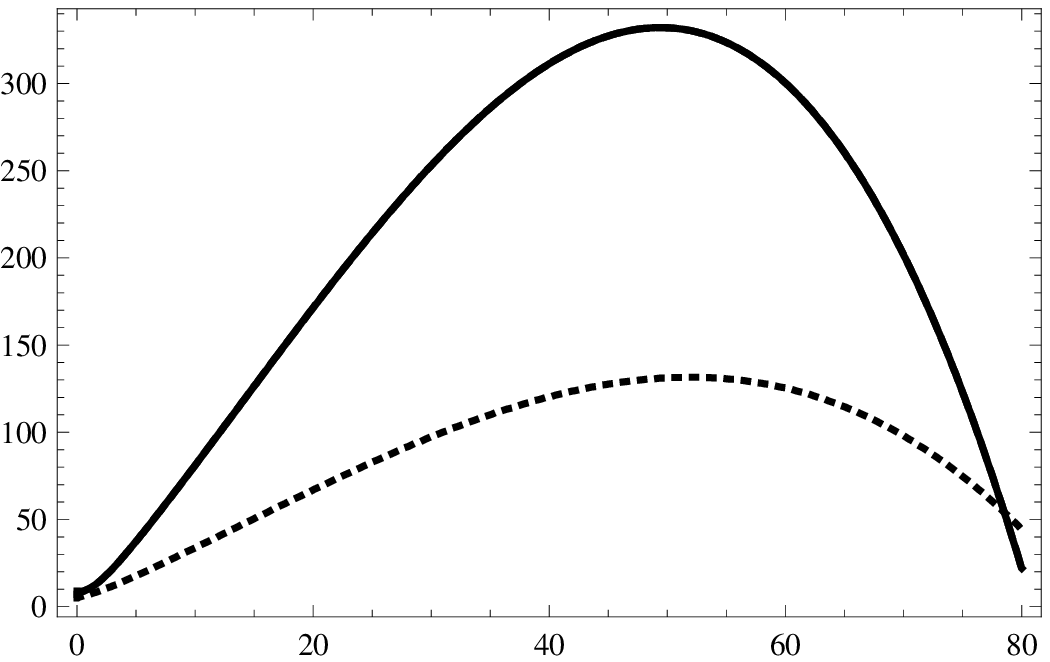}\label{2ndorderSL2Rzero} }
\subfigure[$(A_0, B_{0},C_{0})=(3,5, 7), \alpha' = 1, T_s=0.033$]{\includegraphics[width=0.4\textwidth]{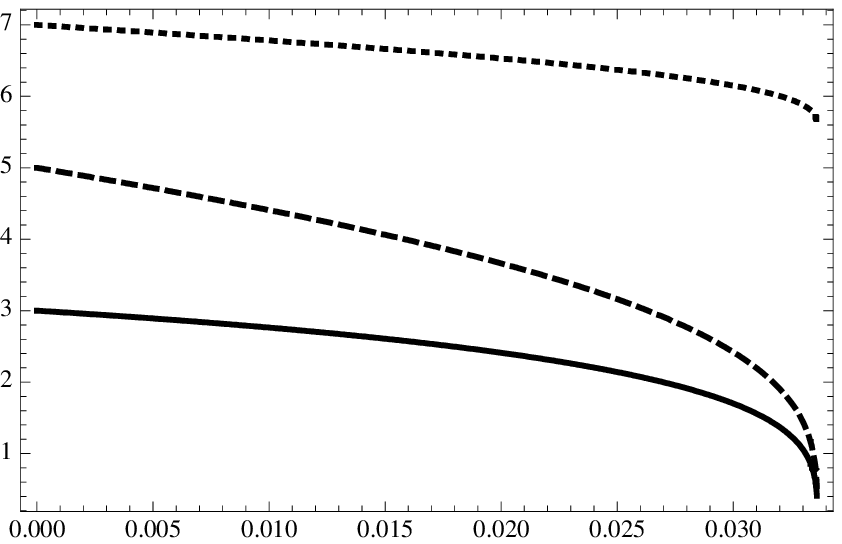}\label{2ndorderSL2Rone} }
\subfigure[$(A_0, B_{0},C_{0})=(7,5, 3), \alpha' = 1, T_s=0.6$]{\includegraphics[width=0.4\textwidth]{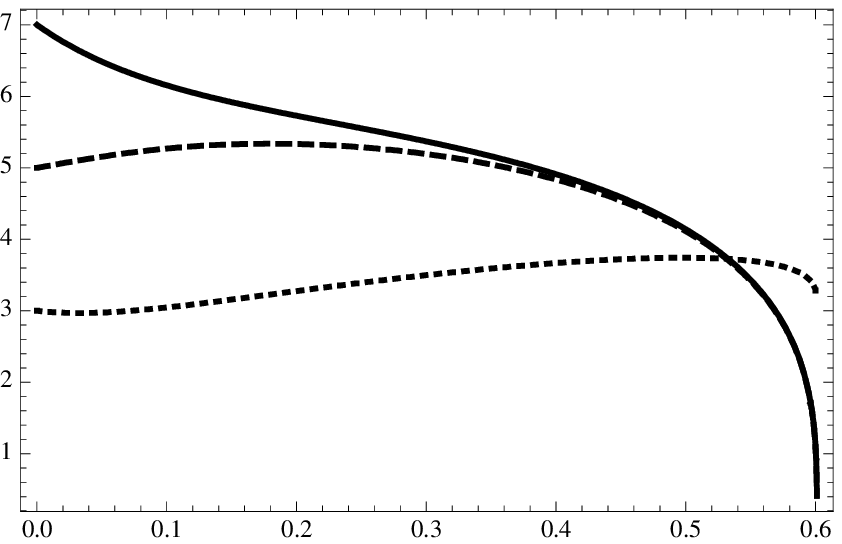}\label{2ndorderSL2Rtwo} }
\subfigure[$(A_0, B_{0},C_{0})=(7,5, 3),\alpha' =-1,T_s=-0.07$]{\includegraphics[width=0.4\textwidth]{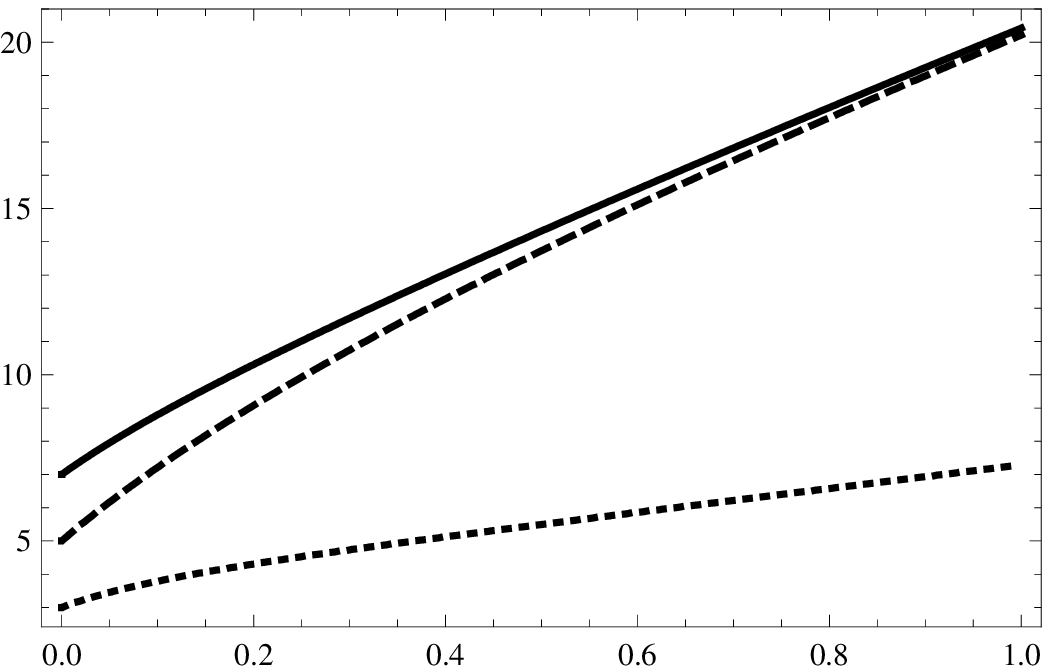} \label{2ndorderSL2Rthree}}
\subfigure[$(A_0, B_{0},C_{0})=(7,5, 3),\alpha' = 0,T_s=-0.24$]{\includegraphics[width=0.4\textwidth]{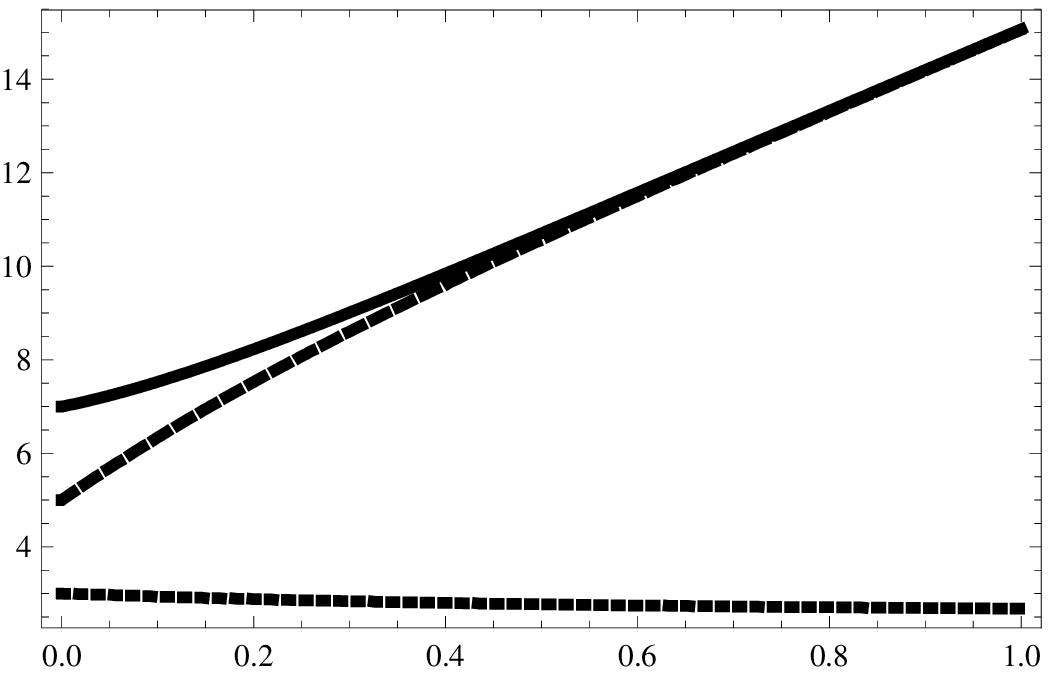}\label{2ndorderSL2Rfour}}
\subfigure[$\alpha' = -1$ and $\alpha' = 0$(thick)]{\includegraphics[width=0.4\textwidth]{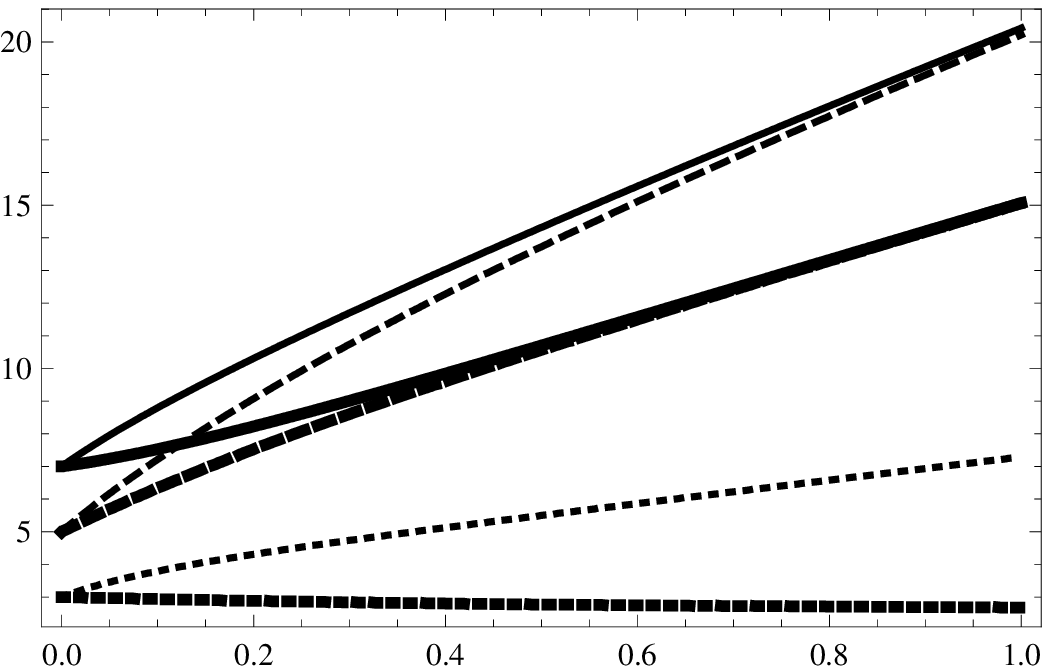} }
\caption{$A(t),B(t),C(t)$ vs $t$ for  $\alpha'=1$, $\alpha'=-1$and $\alpha'= 0$ for $\widetilde{\text{SL}(2,\mathbb{R})}$}
\label{2ndorderSL2R}
\end{figure}

The evolution of $C(t)$ is straightforward. We note that
there exists a lower bound, in the following sense,\\
\begin{align}\label{evolCsl2rURF}
\4{dC}{dt}&=\4{4}{AB}\bigl(\lb(A-B\rb)^2-C^2\bigr)\nonumber\\
&=4\bigl(\lb(\4{A}{B}+\4{B}{A}\rb)^2-\4{C^2}{AB}-2\bigr)\ge-4
\end{align}
Eqn.\ref{evolCsl2rURF} shows that $C(t)(=C_0-4t)$ is monotonically decreasing 
--a point of difference from what we find for normalised Ricci flow. 
Similarly, it is not difficult to show that $B(t)$ is also 
 monotonically increasing, if we assume $A_0>B_0>C_0$. However,
 the nature of the evolution of $A(t)$ shows an increase though it is
 not monotonic. 
This can be justified as follows. We have\\
\be
\4{dA}{dt}=\4{4}{BC}\lb(A+B+C\rb)\lb(B+C-A\rb)
\ee
Thus, $A(t)$ will increase monotonically provided $A<(B+C)$, 
otherwise it will decrease initially and then increase. 
We further note that as $A(t)$ and $B(t)$ increase they approach each other.
This feature follows
from Eqn.\ref{evola-bsl2rURF1}, assuming $A_0>B_0>C_0$. We can write\\
\begin{align}\label{evola-bsl2rURF2}
\4{d(A-B)}{dt}&=-\4{4}{ABC}\lb(A-B\rb)\lb(A+B-C\rb)\lb(A+B+C\rb)\nonumber\\
&\le-\4{4}{AB}\lb(A-B\rb)\lb(A+B-C\rb)\nonumber\\
&\le-12\4{C}{AB}(A-B)
\end{align}
All the above stated features for $\alpha'=0$ are shown in Fig.\ref{2ndorderSL2Rfour}.

\noindent {\bf (b)} For $\alpha'\ne0$ it is difficult to understand even 
the qualitative nature 
of the evolution of scale factors, so we tried numerical solutions for some 
particular initial values. Before we start discussing the numerical 
results, let us recall that the $\widetilde{\text{SL}(2,\mathbb{R})}$ class 
does not contain Einstein metrics. The unnormalized Ricci flow does not 
converge-- instead it evolves towards the 
pancake degeneracy. But the inclusion of higher order terms 
seems to oppose this. For $\alpha' = 1$ the scale factors converge for 
different initial values -a fact depicted in Figs. (19)(a-c). 
On the other hand, the $\alpha' = -1$ 
case resembles the unnormalized Ricci flow where two of the expanding 
scale factors approach each other, though the past singularity time 
is larger (Fig.\ref{2ndorderSL2Rthree}). The generic 
behavior of scale factors 
for $\alpha'=0$  are retained for $\alpha'=-1$, except for the evolution 
of $C(t)$--a fact depicted in the Figs. (19)(e-f). 

Finally, 
we move on the special case where $A=B\ne C$ 
 where we have an exact relation between the scale factors for $\alpha'=0$.
\subsubsection{{\bf Special case: $A=B\ne C$}}
The flow equations  in this case are:
	\begin{subequations}\label{sl2rspecial}
		\begin{align}  		
			\f{dA}{dt} & =  8 + 4\frac{C}{A} -4 \alpha' ~\frac{8 A^2+12 C A+5 C^2}{A^3}\label{sl2rspecial1}\\
			\f{dC}{dt}  &= -4\frac{C^2}{A^2} - 4\alpha' ~\frac{C^3}{A^4}\label{sl2rspecial2}
		\end{align}
	\end{subequations}	
For $\alpha'=0$ we can exactly solve the system and arrive at a relation 
between $A$ and $C$ given as\\
\be
C(t)=\4{A(t)^2 ~k_1~ \text{exp}(-A(t)^2)}{\4{C(t)}{A(t)}+k_1~\text{exp}(-A(t)^2)}
\ee

\begin{figure}[htbp] 
\centering
\subfigure[$(A_0,C_{0})=(7,5), \alpha' = 1$]{\includegraphics[width=0.41\textwidth]{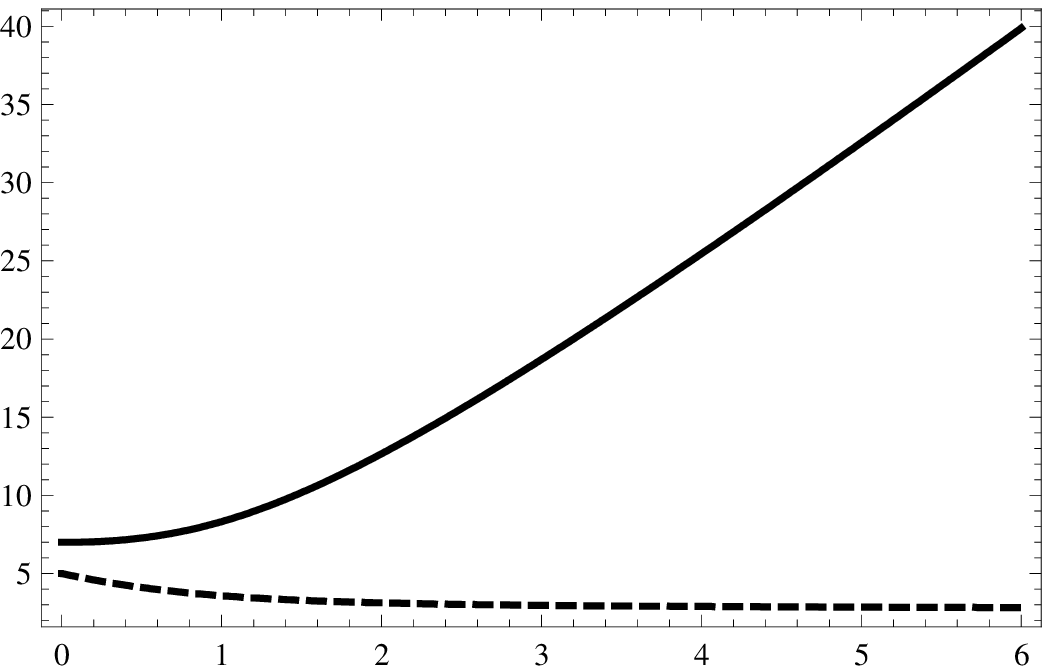}\label{2ndorderSL2Rspone} }
\subfigure[$(A_0,C_{0})=(3.5,2), \alpha' = 1, T_s=0.11$]{\includegraphics[width=0.41\textwidth]{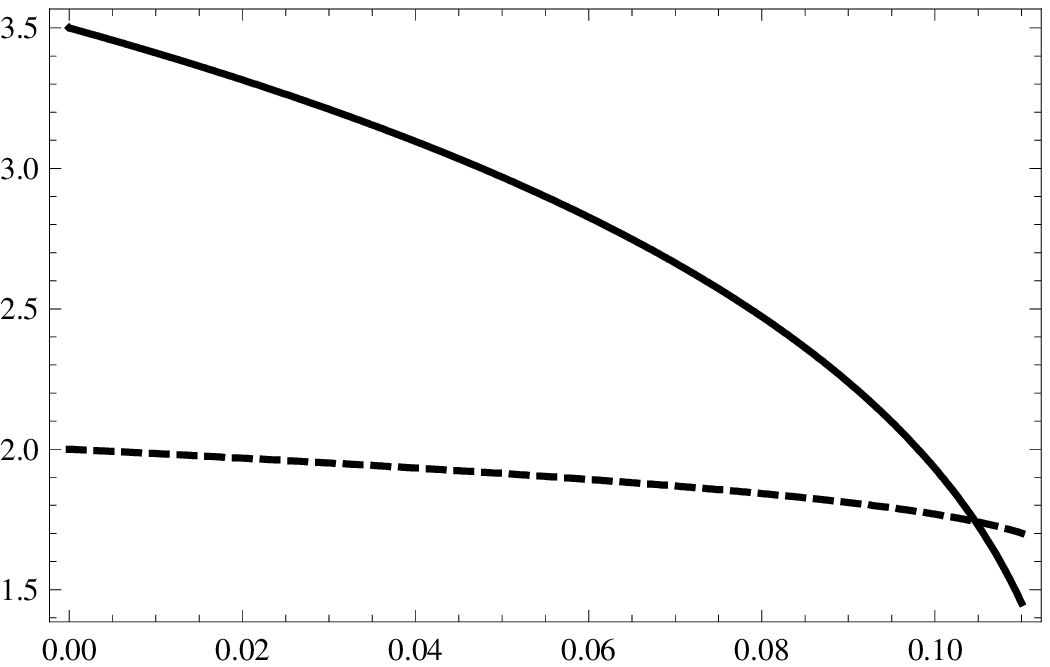}\label{2ndorderSL2Rspthree} }
\subfigure[$(A_0, C_{0})=(7,5), \alpha' = -1, T_s=-0.13$]{\includegraphics[width=0.41\textwidth]{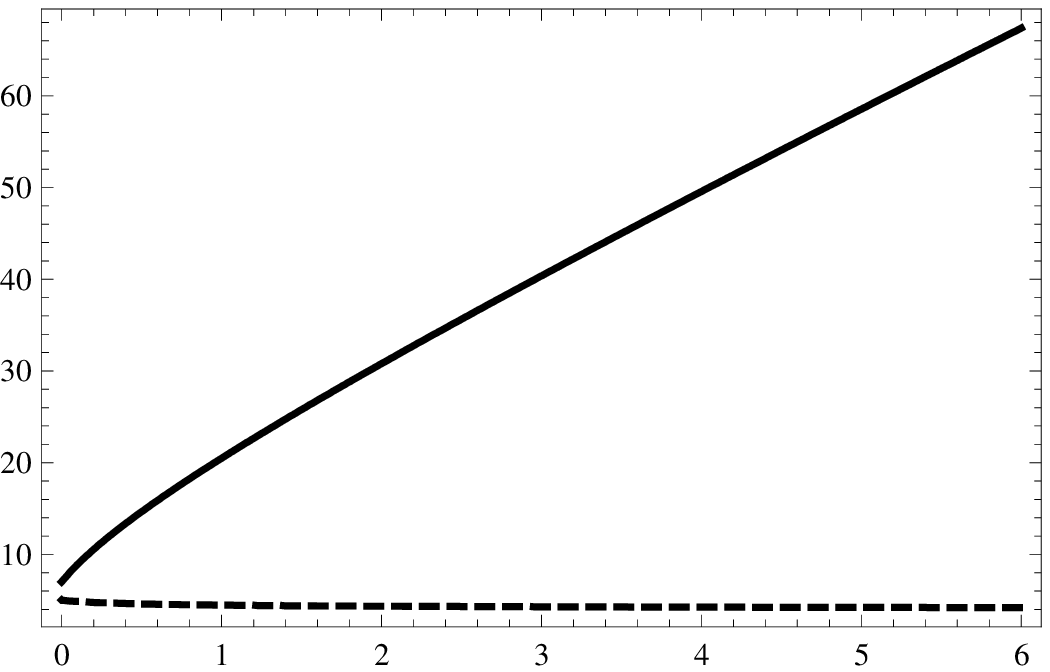} \label{2ndorderSL2Rsptwo}}

\caption{$A(t),C(t)$ vs $t$ for  $\alpha'=1$, and $\alpha'= -1$}
\label{2ndorderSL2Rsp}
\end{figure}

We can also check that for $\alpha'=0$ the scale factors diverge from each other. This may be noted from the rate of evolution of the difference of
scale factors, as given below.
\be
\4{d(A-C)}{dt} =8~+~\4{C}{A}~\lb(1+\4{C}{A}\rb)\ge 8
\ee
For $\alpha'=-1$ we can write down the expressions for the 
differences of the scale factors 
from Eqn.\ref{sl2rspecial}, as follows,

\begin{align}
\4{d(A-C)}{dt} &= 8~+~\4{4\lb( 5C^2+A^2\lb(8+C\rb)+AC\lb(12+C\rb)\rb)}{A^3}-\4{4}{A}\lb(\4{C}{A}\rb)^3 
\\ & \ge 8+\4{4\lb(A+C\rb)\lb(5C+A\lb(C+7\rb)\rb)}{A^3}
\end{align}
where we maintain $A>C$. Here the scale factors are diverging. 
This fact has been numerically verified in Fig.\ref{2ndorderSL2Rsptwo}. 
We do not show the evolution of the scale factors $A(t)$ and $C(t)$ 
for $\alpha'=0$ because their nature are similar to the $\alpha'=-1$
case.  Lastly, if  $\alpha'=1$, the evolution of the scale factor $C(t)$  
can be found from Eqn.\ref{sl2rspecial2}. However, the evolution of $A(t)$ 
will be different for different initial conditions --- features 
which are shown in Fig.\ref{2ndorderSL2Rspone}and inFig.\ref{2ndorderSL2Rspthree}. Let us now turn to the phase potraits.

\subsubsection{\bf Phase plots} 

\begin{figure}[htbp]
	\centering
	\subfigure[$\alpha' = 1$]{\includegraphics[width=0.3\textwidth]{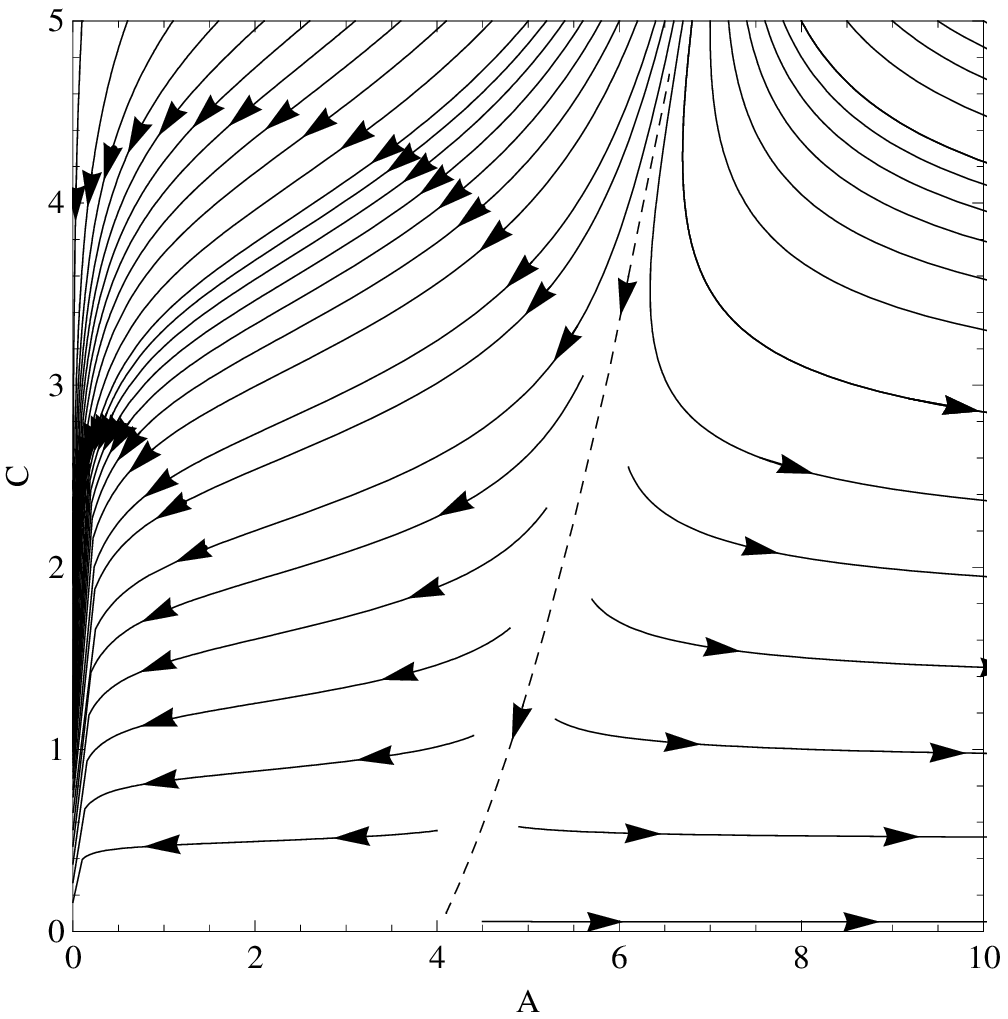}\label{sl2r 1}}
	\subfigure[$\alpha' = 0$]{\includegraphics[width=0.3\textwidth]{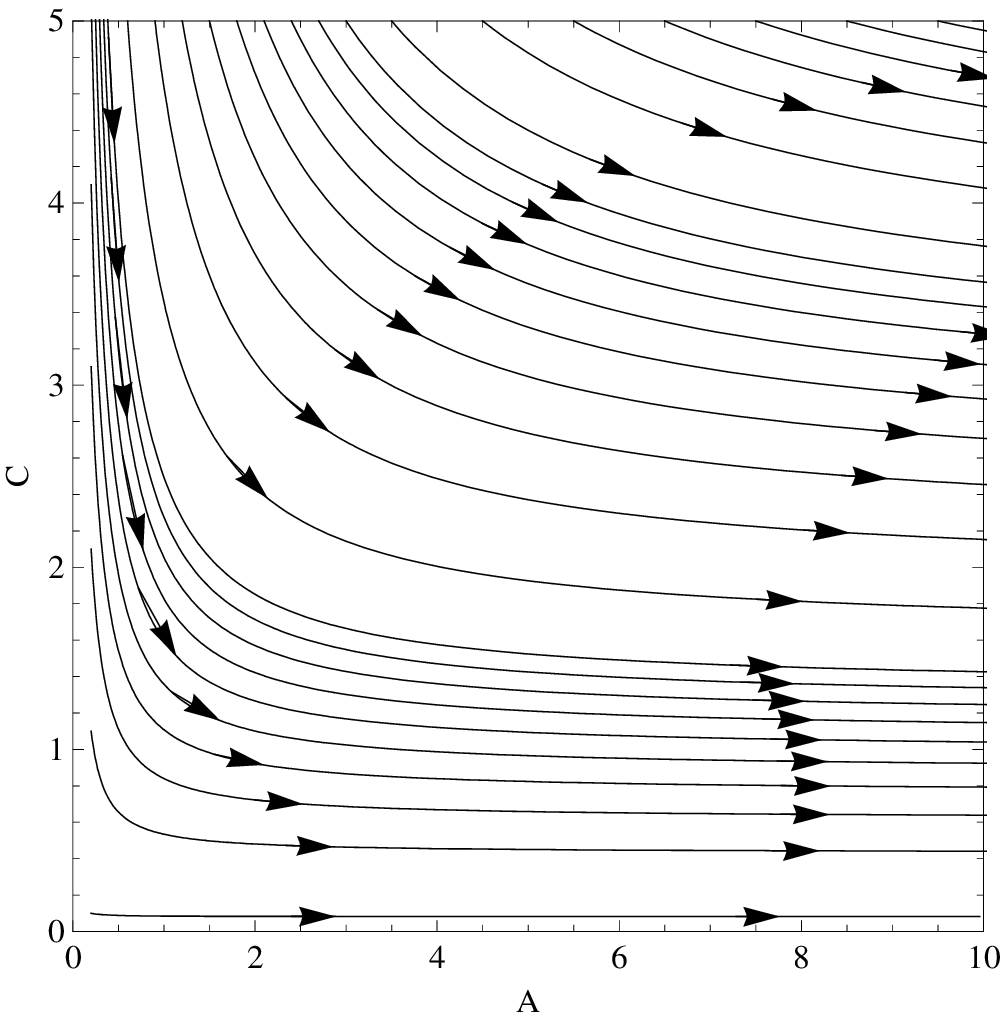}\label{sl2r 0}}
	\subfigure[$\alpha' = -1$]{\includegraphics[width=0.3\textwidth]{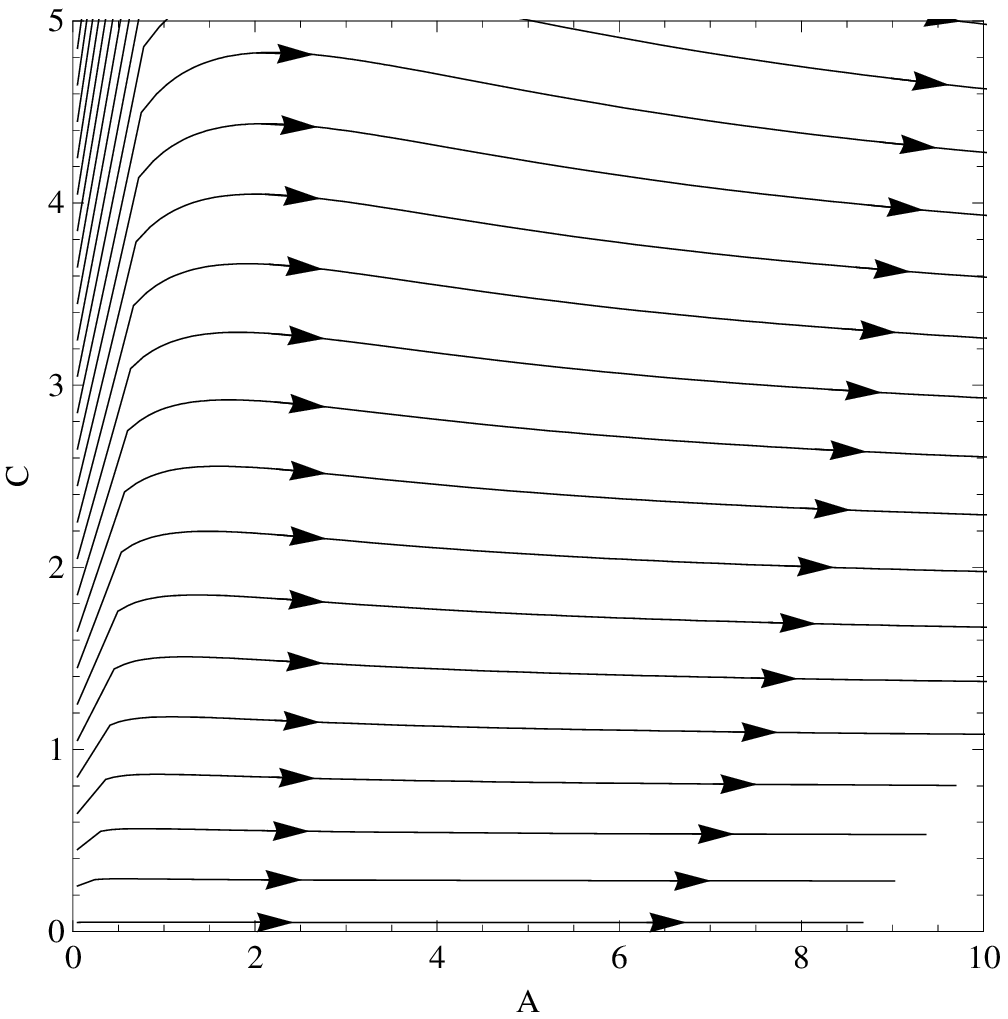}\label{sl2r -1}}
	\caption{$2nd$ order flow on $\widetilde{\text{SL}(2,\mathbb{R})}$ for $A=B$}
	\label{sl2r}
	\end{figure}

As before, we plot the above--mentioned reduced system in Fig.\ref{sl2r}. 
For $\alpha' = 0,-1$ trajectories flow to larger values of $A$. But in the 
case of $\alpha' = 1$, the presence of a fixed point $(A,C)= (4\alpha',0)$ 
makes the flow more interesting. In this case, some trajectories from a certain 
region in the phase space go towards larger $A$ as before, but others 
converge to the singularity at $A=0$. There also exists a critical curve, 
which flows into the fixed point, and demarcates the regions with different 
asymptotics.

\subsubsection{\bf Curvature evolution}

\begin{figure}[h]
	\centering
	\subfigure[Evolution of scalar curvature for various $\alpha'$, $A_0,B_0,C_0$]{\includegraphics[width=0.72\textwidth]{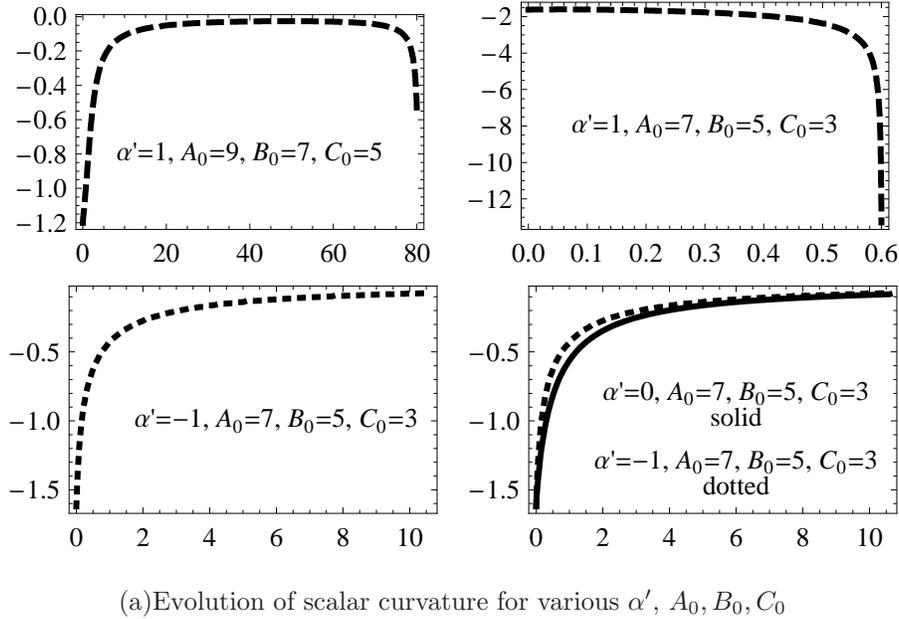}\label{}}
	\caption{Evolution of scalar curvature  for $\widetilde{\text{SL}(2,\mathbb{R})}$.}
	\label{curvevolisl2r}
	\end{figure}

We end by plotting the scalar curvature evolution. 
This is not very different from the earlier cases (Sol and Nil), for 
$\alpha'=0$ and $\alpha'=-1$ where the curvature asymptotically goes to zero 
starting from a negative value. However, when $\alpha'=+1$ we find some
interesting features which can be seen in Fig.[\ref{curvevolisl2r}]. 
For a particular initial 
condition, the scalar curvature increases from a negative value towards zero, 
stays there for a while and then drops to lower values once again. 
\section{Conclusions}
Our overall aim in this work has been to study in detail the 
consequences of second order
(in Riemann curvature) geometric flows on three dimensional homogeneous
spaces, using analytical, semi-analytical
and numerical methods. 
Through the analysis carried out,
we believe we have been  able to obtain quite a few characteristics 
which seem to arise for second order geometric flows  
on three dimensional homogeneous (locally) geometries.
Here, we briefly summarize our
results and mention a few possibilities for the future.

From our results, we can say that for manifolds which do
not contain Einstein metrics the flow characteristics 
show new features essentially caused by the inclusion of 
the higher order term. On the other hand, the class of group manifolds 
which admit metrics of Einstein class (whether it is flat or round),
results for second order flows appear to be 
refinements over known results for Ricci flows.
However, generic new characteristics do arise with varying sign of $\alpha'$ 
--in particular, a negative $\alpha'$-- and these have been noted in our
work.

In several restricted cases (i.e. where two of the scale factors
are related), we are able to solve the flow equations exactly. 
We have obtained analytical solutions in such restricted cases 
for {\em Nil}, {\em Sol},  $\widetilde{\text{Isom}(\mathbb{R}^{2})}$, 
$\widetilde{\text{SL}(2,\mathbb{R})}$ manifolds. 
The exact solutions are instructive because they help in obtaining
analytical expressions for fixed points (curves) 
as well as in understanding the evolution of the scale factor. 
In addition, we also use them for checking our numerics.
 
A generic observation is the fact that the singularity time changes due 
to the inclusion of
higher orders in the flow equations. This pattern is noticeable throughout
in our numerical work. 

The results for $\alpha'<0$ are, in quite a few cases, strikingly
different from those for $\alpha'=0$ or $\alpha'>0$. Even the evolution 
of the scalar curvature exhibits a different behaviour in many of the cases
studied. 

For $SU(2)$ all the scale factors do not converge for $\alpha'=-1$ which is 
exactly opposite to the  characteristics for $\alpha'=1$ and $\alpha'=0$. 
Here, if  $\alpha'=0,1$, the scalar curvature increases but when  
$\alpha'=-1$ it increases first and then decreases. 
The appearance of negative scalar curvature in $SU(2)$ is also
noted. In the case for
 Sol manifold  the behavior of the scale factors depend on 
different initial conditions for $\alpha'=\pm 1$---so does the 
evolution of the scalar curvature.
  
In the case of $\widetilde{Isom(\mathbb{R}^{2})}$
all the scale factors converge for $\alpha'=1$. When $\alpha'=-1$
one scale factor decreases while the other two increase.  
In the last case, SL(2,R), we find that for $\alpha'=1$ all the scale factors 
may initially increase but they converge towards a singularity. However, 
for $\alpha'=-1$ all the scale factors diverge. 

In all cases we have obtained the phase portraits,
which we feel, helps in visualising the flow features as well as the
fixed points. 
We also provide a summary of all our results in several
tables in the sections as well as, in the end. 

Among possible future directions, we mention a few below.

$\bullet$ It would be interesting to pursue the approach presented in 
\cite{glick} for
such higher order flows. 

$\bullet$ A more systematic and exhaustive analysis of the
stability and classification of fixed points, which is largely an algebraic
problem, can be carried out for such second order flows. This will surely
shed more light on the behaviour of these flows from an analytical
perspective.  
 
$\bullet$ Given the fact that such three manifolds do
arise in various physically relevant contexts, it will be nice to know 
whether our results on higher order flows
can help us understand such scenarios in any meaningful way. 

$\bullet$
The relevance of our results in the context of 
renormalisation group flows of the bosonic nonlinear $\sigma$-model 
deserve some attention. 

$\bullet$ 
Since homogeneous four manifolds are already classified and studied
with reference to Ricci flows \cite{fourmanifolds} it would be
worthwhile extending our results to four manifolds.

We hope to address some of these issues in future articles.  
   
\



\begin{sidewaystable}[p]\small
\caption{Comparisons of scale factors for different cases}\label{table3chap03}\centering
\renewcommand{\tabularxcolumn}[1]{>{\arraybackslash}m{#1}}
\begin{tabularx}{\textheight}{+l^c^Z^Z^Z}
\toprule
Manifold &  Bianchi Type & $\alpha'=0$ & $\alpha'=1$ & $\alpha'=-1$\\
\otoprule
$\text{SU}(2)$ & IX & all scale factors converge to finite time singularity & different convergence rates, singularity times & does not always converge to a singularity\\
\midrule
Nil & II & pancake degeneracy & depends on initial data & pancake degeneracy\\
\midrule
Sol & $\text{IV}_{-1}$ & cigar degeneracy & depends on initial data & depends on initial data\\
\midrule
$\widetilde{\text{Isom}(\mathbb{R})^2}$ & $\text{VII}_{0}$ & one decreases and other two increase, converges to flat & all decrease and converge to flat & one decreases and other two increase, converges to flat\\
\midrule
$\widetilde{\text{SL}(2,\mathbb{R})}$ & VIII & one decreases and two increase & increases or decrease & all increase and diverging\\
\bottomrule
\end{tabularx}
\end{sidewaystable}
\newpage
\setlength{\heavyrulewidth}{0.1em}
\begin{table}[h]\small
\caption{Comparisons of scalar curvature between different cases}\label{table4chap03}\centering
\renewcommand{\tabularxcolumn}[1]{>{\arraybackslash}m{#1}}
\begin{tabularx}{\textwidth}{+Y^Z^Z^Z}
\toprule\rowstyle{\bfseries}
Manifold & $\alpha'=0$ & $\alpha'=1$ & $\alpha'=-1$\\
\otoprule
$\text{SU}(2)$  & Scal. curvature increases can be -ve. & Scal. curvature increases or decreases (entirely -ve. or +ve.) & Scal. curvature increases/decreases, can be -ve. \\
\midrule
Nil & asymptotically goes to flat from -ve. value & increases $(A_0<B_0<C_0)$ or decreases $(A_0>B_0>C_0)$ depending on initial value & asymptotically goes to flat from -ve. value\\
\midrule
Sol & asymptotically goes to flat from -ve. value & increases $(C_0<A_0<B_0)$ or decreases $(A_0>B_0>C_0)$ depending on initial value & asymptotically goes to flat from -ve. value\\
\midrule
$\widetilde{\text{Isom}(\mathbb{R})^2}$ &  starts with -ve. value, asymptotically goes to flat & starts with -ve. value, may develop singularity, asymptotically flat & starts with -ve. value, asymptotically flat\\
\midrule
$\widetilde{\text{SL}(2,\mathbb{R})}$  & asymptotically goes to flat from -ve. value & depending on initial value, curvature may decrease (tend to zero) or increase & asymptotically goes to flat from -ve. value\\
\bottomrule
\end{tabularx}
\end{table}

\begin{acknowledgements}
SD acknowledges H. Seshadri, S. Panda for useful discusions and 
thanks the Institute of Mathematical Sciences, Chennai, India,
for support through a post-doctoral fellowship.
\end{acknowledgements}

\newpage

\end{document}